\documentclass[review,hidelinks,onefignum,onetabnum,reqno]{siamart250211}
\usepackage{amssymb}
\usepackage{xcolor}
\usepackage{amsmath}
\usepackage{lipsum}
\usepackage{amsfonts}
\usepackage{graphicx}
\usepackage{epstopdf}
\usepackage{algorithmic}
\ifpdf
  \DeclareGraphicsExtensions{.eps,.pdf,.png,.jpg}
\else
  \DeclareGraphicsExtensions{.eps}
\fi

\newsiamremark{remark}{Remark}
\newsiamremark{hypothesis}{Hypothesis}
\crefname{hypothesis}{Hypothesis}{Hypotheses}
\newsiamthm{claim}{Claim}
\newsiamremark{fact}{Fact}
\crefname{fact}{Fact}{Facts}

\headers{Semi-Discrete Convexification and Deep Learning for EIT}{Michael V. Klibanov, Kirill V. Golubnichiy, and Benjamin Jiang}

\usepackage{amsopn}

\ifpdf
\hypersetup{
  pdftitle={A Carleman Semi-Discrete Convexification Method Combined With Deep
Learning for Electrical Impedance Tomography},
  pdfauthor={Michael V. Klibanov, Kirill V. Golubnichiy, and Benjamin Jiang}
}
\fi

\begin{document}

\nolinenumbers

\title{A Carleman Semi-Discrete Convexification Method Combined With Deep
Learning for Electrical Impedance Tomography}
\author{Michael V. Klibanov \thanks{
Department of Mathematics and Statistics, University of North Carolina at
Charlotte, Charlotte, NC, 28223, USA, mklibanv@charlotte.edu}, \and Kirill
V. Golubnichiy \thanks{%
Department of Mathematics and Statistics, Texas Tech University, Lubbock,
TX, 79409, USA, kgolubni@ttu.edu,} \and Benjamin Jiang\thanks{%
Department of Mathematical and Statistical Sciences, University of Alberta,
Edmonton, AB, T6G 2J5, Canada, bjjiang@ualberta.ca}}
\date{}
\maketitle

\begin{abstract}
In this paper, a new, two-stage numerical procedure for Coefficient Inverse Problems is proposed. On the first stage a new semi-discrete version of the Carleman estimate-based convexification globally convergent numerical method is applied on a coarse grid. On the second stage, solution obtained on the first stage is used as the input for the deep learning procedure. Our analytical development allows us to obtain an a priory accuracy estimate for the entire procedure, under two reasonable assumptions. Numerical testing demonstrates a good performance.
\end{abstract}

\sloppy
\overfullrule=0pt

\nolinenumbers



\pagestyle{myheadings} \thispagestyle{plain} 
\markboth{M. V. KLIBANOV, K. V.
GOLUBNICHY, AND B. JIANG}{}




\section{Introduction}

\label{sec:1}

We present a new two-stage numerical procedure for a broad class of
Coefficient Inverse Problems (CIPs) for Partial Differential Equations
(PDEs). This idea, combined with two reasonable assumptions, allows us to
obtain a priory accuracy estimates for a Deep Learning Procedure (DNN) for a
CIP. \textquotedblleft Broad class" needs to be specified. To better
illustrate our concept as well as to simplify the presentation, we focus
below on the CIP\ of the Electrical Impedance Tomography (EIT)\ in the 2-D
case, see, e.g. \cite{Borcea,KEIT} for this problem. Thus, all results below
are related only to this CIP. More precisely, we focus on the recent version
of the convexification method for EIT, which uses the viscosity term \cite%
{EITIP2025}. We refer to, e.g. \cite{Harrach1,Harrach2,Harrach3} for some
other approaches to numerical solutions of CIPs for EIT. There are many
publications where deep neural network is used for EIT, see, e.g. \cite%
{Santos}. Now, since the convexification method is applicable to a broad
class of CIPs see, e.g. \cite{KEIT,KL,KLZhyp2,TTTP,Ktransp,EITIP2025} about
this method, then we conjecture that the idea of this paper might be
extended to many other CIPs.

On the first stage of our two-stage procedure, the convexification method is
applied on a coarse grid. This method is based on Carleman estimates, and it
is a numerical development of the initially purely theoretical idea of \cite%
{BukhKlib}. The publication \cite{BukhKlib} is the first one, in which the
tool of Carleman estimates was introduced in the field of CIPs. The goal of 
\cite{BukhKlib} was to prove uniqueness and stability theorems for CIPs,
see, e.g. \cite{Isakov,KT,KL} for some follow up publications.

The convexification method was originated in \cite{Klib95,Klib97}. The
reason why the convexification is used is that conventional numerical
methods for CIPs are based on the minimization least squares cost
functionals, which are non convex. Hence, those techniques face the well
known phenomenon of multiple local minima and ravines of corresponding cost
functionals, see, e.g. \cite{Beilina1,Beilina2,Chavent,Giorgi,Gonch1,Gonch2}
for such functionals. In addition, we refer to \cite{Scales} for a numerical
example of multiple local minima. Hence, any minimization procedure for such
a functional can be trapped at any local minimum. The main advantage of the
convexification method is that it is free from this phenomenon.

The key point of any version of the convexification method is a construction
of a weighted Tikhonov-like regularization functional $F$. The main feature
of $F$ is the presence of a weight function, which is the Carleman Weight
Function (CWF). The CWF is used as the weight in the Carleman estimate for
the corresponding PDE operator. The convergence analysis establishes then
that the functional $F$ is strongly convex on an appropriate convex bounded
set $G\subset H$ of the diameter $d>0,$ where $H$ is an appropriate Hilbert
space. Since restrictions on the smallness of $d$ are not imposed, then that
functional is globally strongly convex. In addition, the existence and
uniqueness of the minimizer of $F$\ on the set $G$\ are proven. Furthermore,
an accuracy estimate of that minimizer of $F$ on $G$ implies that the
distance between the minimizer and the true solution of the corresponding
CIP is small for an appropriate choice of parameters. These facts lead to
the conclusion, which is proven in, e.g. above cited publications about the
convexification method \cite{KEIT,KL,KLZhyp2,TTTP,Ktransp,EITIP2025}: that
it converges globally, i.e. its convergence to the true solution is
rigorously guaranteed without an a priory knowledge of a good starting point
for the solution.

\subsection{Two-stage procedure}

\label{sec:1.1}

The convexification method performs slowly. \sloppy Therefore, \fussy it is
reasonable to use it only as the first stage of our numerical procedure, in
which a coarse grid is used. This speeds up computations quite significantly.
This is certainly advantageous from the computational standpoint because of
the time-consuming constrained minimization procedure within the
convexification method. Indeed, the same computer has required:

\begin{equation}
\left.
\begin{array}{c}
57\text{ minutes on the fine grid for each image of \cite{EITIP2025}}\\
\text{and only 4 minutes on the coarse grid of the current paper.}
\end{array}
\right.
\label{1.1}
\end{equation}

However, since only continuous versions of the convexification method were
considered in the past, the use of a coarse grid requires the development of
the theory of a semi-discrete version of that method. In this case, the PDE
operators involved in $F$ are written in finite differences, whereas the
penalty term is written in the continuous form.

The rapid solution obtained by the semi-discrete convexification method on a
coarse grid is quite blurry. Hence, to improve its accuracy, we use the
second stage: deep learning.

The roles of the two stages are the same during training and testing, but they
should be distinguished explicitly. During the training stage, for every
training sample we first compute its own reconstruction by the semi-discrete
convexification method on the coarse grid. This reconstruction is used as the
input image $x$ of the DNN, while the corresponding true image is used as the
training target $y$. Thus, the DNN is trained to learn a mapping from a
coarse-grid convexification reconstruction to a refined reconstruction.

During the testing stage, the convexification procedure is applied again,
independently, to every new testing sample. In particular, a new coarse-grid
convexification reconstruction is computed for each testing sample and is
then supplied as the input to the already trained DNN. No convexification
reconstruction from the training set is reused during testing. The true image
of a testing sample is not used by the reconstruction algorithm and is used
only afterward for quantitative evaluation and comparison.

Therefore, the complete testing pipeline for each new sample is
\[
\begin{aligned}
\text{measurement data}
&\longrightarrow
\text{coarse-grid convexification reconstruction} \\
&\longrightarrow
\text{trained DNN}
\longrightarrow
\text{refined reconstruction}.
\end{aligned}
\]

The resulting two-stage numerical procedure works quite rapidly and
accurately.

\subsection{Four main new elements of this paper}

\label{sec:1.2}

These elements are:

\begin{enumerate}
\item The so-called \textquotedblleft $h-$strong convexity property" is
introduced. This property is established for that semi-discrete analogue of
the CWF-weighted Tikhonov-like objective functional, see Theorem 5.1 and
Remark 5.1 below. Here $h\in \left( 0,1\right) $ is the grid step size.

\item Let $F_{sd}$ be that semi-discrete functional
defined on a coarse grid. Suppose that a minimizer of $F_{sd}$ on
that grid exists (see subsection 1.3). Then the result of item 1 allows us
to estimate the distance between that minimizer of $F_{sd}$ and
the minimizer of the above continuous functional $F$, see Theorem
5.2 below. In addition, the distance between that minimizer of $F_{sd}$%
 and the true solution of our CIP is estimated in Theorem 5.2.
Finally, the distances between the terms of a minimizing sequence of $F_{sd}$%
 and that minimizer of $F$ as well as the true solution
of our CIP are estimated without the assumption of the existence of a
minimizer of $F_{sd},$ see Theorem 5.3.

\item In the proposed deep-learning procedure, the solution
obtained by the semi-discrete version of the convexification method on a
coarse grid is used as the input image to the DNN. This construction is used
for both the training and testing samples.

\item Therefore, items 2 and 3 result in an a priory
accuracy estimate for the convexification reconstruction used as the input
image to the DNN throughout the deep-learning procedure. This estimate is 
\begin{equation}
O\left( \sqrt{\alpha }+\sqrt{h}\right) ,\text{ as }\sqrt{\alpha }+\sqrt{h}%
\rightarrow 0^{+},  \label{1.3}
\end{equation}%
where $\alpha \in \left( 0,1\right) $ is a small regularization parameter
used in $F_{sd}.$
\end{enumerate}

\textbf{Remark 1.1.} We recall that accuracy estimate (\ref{1.3}) is
obtained only under the assumption of the existence of a minimizer of $F_{sd}
$, also see subsection 1.3, Theorem 5.2 and Remark 5.2. 

\subsection{Some considerations}

\label{sec:1.3}

The issue of extensions of various results for the continuos version 
$F$of the CWF-weighted Tikhonov-like objective functional to its
semi-discrete analog $F_{sd}$ is a quite challenging one. It is
natural, therefore, that there are two important properties for the
continuous case, which we are yet unable to extend analytically to the
semi-discrete case. These analytical extensions are challenging tasks, and
they are outside of the scope of the current publication. We hope to address
them in the future. 

Those properties are:

\begin{enumerate}
\item We cannot prove the existence of a minimizer of $F_{sd}$ 
on a certain set. Therefore, we assume that existence: compare item
3 of Theorem 3.2 for the continuous case with point 3 of Theorem 5.1 for the
semi-discrete case. 

\item While the global convergence of the gradient descent method of
the minimization of the continuous version $F$ of that functional
was proven in Theorem 5.4 of \cite{EITIP2025}, we cannot prove a similar
result for the semi-discrete case of the functional $F_{sd}$.

\item In addition, we cannot establish analytically the accuracy of
a certain interpolation, see subsection 7.1. For this reason, we impose a
certain quite natural assumption there.
\end{enumerate}

Nevertheless, the $h$ $-$ strong convexity property
enables us to establish the following three features of that semi-discrete
functional $F_{sd}$:

\begin{enumerate}
\item Assuming the existence of a minimizer of \ the functional $%
F_{sd}$, we estimate the accuracy of this minimizer. More
precisely, we estimate in Theorem 5.2 the distance between this minimizer
and the minimizer of the continuos version of that functional. In addition,
we estimate the distance between that minimizer and the true solution of the
continuous case. Both these estimates are as in (\ref{1.3}). 

\item Therefore, it follows from (\ref{1.3}) that even if the
semi-discrete functional $F_{sd}$ has several minimizers, still
they are close to each other as well as to the minimizer of $F$
and to the true solution, as long as $\sqrt{\alpha }+\sqrt{h}$ is
sufficiently small. This is good for computations. In addition, for the sake
of comparison, we refer to Theorem 3.3 for accuracy estimates of the
minimizer for the continuous case.

\item Second, for any fixed value of the regularization parameter $%
\alpha $, there exists a minimizing sequence for $F_{sd}$ on an
appropriate set, which is specified in subsection 5.3. This sequence moves
values of the functional $F_{sd}$ to its infimum value on that
set.\ So, we estimate in Theorem 5.3 the distances between some terms of
that sequence and the minimizer of $F.$\ In addition, we estimate
the distances between those terms and the true solution of our CIP.
Estimates in both cases are as in (\ref{1.3}). In other words, those
distances are small as long as $\sqrt{\alpha }+\sqrt{h}$\ is
small. 
\end{enumerate}

With respect to the previous paragraph, we recall that in the
classical regularization theory convergence of a minimizing sequence of a
general Tikhonov regularization functional is proven only if the
regularization parameter tends to zero, which is unlike our case of any
fixed value of $\alpha $ $.$ In addition, estimates of
distances between terms of that sequence and the true solution are not
provided in the classical regularization theory \cite{T}.

\subsection{An additional assumption}

\label{sec:1.4}


Recall that it is stated in subsection 1.1 that we use the solution
computed by the semi-discrete version of the convexification method on the
coarse grid as the initial input image for the DNN. We have consistently
observed in our computational experiments that the DNN refines images
obtained on that first stage. Therefore, based on our computational
experience, we impose the following quite reasonable assumption:

\textbf{Assumption 1.1.} We assume that accuracy estimate (\ref{1.3}%
) for the solution of our CIP\ via the minimization of the above
semi-discrete functional $F_{sd}$ on a coarse grid, which we use
as the starting point for the follow up DNN procedure, does not become worse
due to this procedure.

In addition, we formulate in subsection 7.1 below our second
assumption of this sort, which is Assumption 7.1. As soon as Assumptions 1.1
and 7.1 are accepted, we conclude that (\ref{1.3}) provides a priory
accuracy estimate for our entire two-stage procedure outlined in subsection
1.1.

All functions considered below are real valued ones. The paper is organized
as follows. Section \ref{sec:2} formulates the forward and inverse problems
of EIT. In section \ref{sec:3} we outline the convexification method for the
continuous case. In section \ref{sec:4} we introduce the semi-discrete
version of the convexification method. In section \ref{sec:5} we carry out
convergence analysis of that semi-discrete version. In \ref{sec:6} we
describe the minimization procedure for our semi-discrete functional.  In
section\ref{sec:7} we describe the deep learning procedure and discuss a
priory accuracy estimates. Section \ref{sec:8} is devoted to numerical
studies. A summary of results is presented in section \ref{sec:9}.

\section{Statements of Forward and Inverse Problems}

\label{sec:2}

Below $\mathbf{x}=\left( x,y\right) $ denotes points in $\mathbb{R}^{2}.$ We
now construct such a domain $\Omega \subset \mathbb{R}^{2},$ which is
convenient for our computational purpose. Let $B,D$ be two numbers such that 
$0<B<D.\mathbb{\ }$\ Let $a,b$ be two other numbers, where $a>0$. Consider
two concentric disks $P_{D}\left( a,b\right) $ and $P_{B}\left( a,b\right)
\subset $ $P_{D}\left( a,b\right) $ with the center at $\mathbf{x}=\left(
a,b\right) $, and let $E_{B}\left( a,b\right) =\partial P_{B}\left(
a,b\right) $ be the circle, which is the boundary of the disk $P_{B}\left(
a,b\right) ,$ 
\begin{equation}
\left. 
\begin{array}{c}
P_{D}\left( a,b\right) =\{\mathbf{x}=\left( x,y\right)
:(x-a)^{2}+(y-b)^{2}<D^{2}\}, \\ 
P_{B}\left( a,b\right) =\{\mathbf{x}=\left( x,y\right)
:(x-a)^{2}+(y-b)^{2}<B^{2}\}\subset P_{D}\left( a,b\right) , \\ 
E_{B}\left( a,b\right) =\{\mathbf{x}=\left( x,y\right)
:(x-a)^{2}+(y-b)^{2}=B^{2}\}.%
\end{array}%
\right.  \label{1.00}
\end{equation}%
Let $c\in \left( 0,a\right) $ be a number. We refer to Figure 1. Define the
square $\Omega $ as%
\begin{equation}
\Omega =\left\{ \mathbf{x}:x\in (a-c,a+c),y\in (b-c,b+c)\right\} .
\label{1.0}
\end{equation}%
We choose $c$ such that $\overline{\Omega }\subset P_{B}\left( a,b\right) $.
Thus, in addition, by (\ref{1.0})%
\begin{equation}
\overline{\Omega }\cap \left\{ x=0\right\} =\varnothing .  \label{1.2}
\end{equation}

In our setting point sources are 
\begin{equation}
\mathbf{x}_{0}=\mathbf{x}_{0}(\theta )\in E_{B}\left( a,b\right) ,\theta \in %
\left[ 0,2\pi \right] .  \label{2.01}
\end{equation}%
We model the point sources by a smooth, compactly supported function $g(%
\mathbf{x}-\mathbf{x}_{0})$ which approximates the $\delta -$function, 
\begin{equation}
g\left( \mathbf{x}-\mathbf{x}_{0}\right) =\left\{ 
\begin{array}{c}
C_{\xi }\exp \left( \frac{\left\vert \mathbf{x}-\mathbf{x}_{0}\right\vert
^{2}}{\left\vert \mathbf{x}-\mathbf{x}_{0}\right\vert ^{2}-\xi ^{2}}\right)
,\left\vert \ \mathbf{x}-\mathbf{x}_{0}\right\vert <\xi , \\ 
0,\text{ }\left\vert \mathbf{x}-\mathbf{x}_{0}\right\vert \geq \xi%
\end{array}%
\right. ,  \label{2.02}
\end{equation}
where $\xi \in \left( 0,1\right) $ is a sufficiently small number such that 
\[
\left\{ \left\vert \mathbf{x}-\mathbf{x}_{0}\right\vert <\xi \right\} \cap 
\overline{\Omega }=\varnothing ,\text{ }\left\{ \left\vert \mathbf{x}-%
\mathbf{x}_{0}\right\vert <\xi \right\} \subset P_{D}\left( a,b\right) ,%
\text{ }\forall \mathbf{x}_{0}\in E_{B}\left( a,b\right) , 
\]%
and the number $C_{\xi }$ is chosen such that 
\[
\int\limits_{\left\vert \mathbf{x}-\mathbf{x}_{0}\right\vert <\xi }g\left( 
\mathbf{x}-\mathbf{x}_{0}\right) dx=1,\text{ }\forall \mathbf{x}_{0}\in
E_{B}\left( a,b\right) . 
\]
Let $\sigma \left( \mathbf{x}\right) $ denotes the electric conductivity
coefficient. We assume that satisfies the conditions 
\begin{equation}
\sigma \in C^{2}\left( \overline{P_{D}\left( a,b\right) }\right) ,\quad
\sigma (\mathbf{x})\geq 1\mbox{ in }P_{D}\left( a,b\right) ,  \label{1-1}
\end{equation}%
\begin{equation}
\sigma (\mathbf{x})=1\quad \mbox{for }~\mathbf{x}\in P_{B}\left( a,b\right)
\setminus \Omega .  \label{2-1}
\end{equation}
For each $\mathbf{x}_{0}\in E_{B},$the underlying elliptic PDE is:%
\begin{equation}
\left. 
\begin{array}{c}
\nabla _{\mathbf{x}}\cdot \left( \sigma \left( \mathbf{x}\right) \nabla _{%
\mathbf{x}}v\left( \mathbf{x,x}_{0}\right) \right) =-g\left( \mathbf{x-x}%
_{0}\right) ,\text{ }\mathbf{x}\in P_{D}\left( a,b\right) , \\ 
v\left( \mathbf{x,x}_{0}\right) =0,\text{ }\forall \mathbf{x}\in \partial
P_{D}\left( a,b\right) .%
\end{array}%
\right.  \label{2-2}
\end{equation}%
Here $v\left( \mathbf{x,x}_{0}\right) $ is the voltage generated by the
electric current induced at the point source $\mathbf{x}_{0}\in E_{B}.$
Conditions (\ref{1.0})-(\ref{2-1}) imply that for every $\gamma \in \left(
0,1\right) $ and for every $\mathbf{x}_{0}\in E_{B}\left( a,b\right) $ there
exists unique solution $v(\mathbf{x},\mathbf{x}_{0})\in C^{3+\gamma }\left( 
\overline{P_{D}}\left( a,b\right) \right) $ of problem (\ref{2-2}). Here $%
C^{3+\gamma }\left( \overline{P_{D}\left( a,b\right) }\right) $ is the H\"{o}%
lder space \cite{GT}. By the maximum principle \cite{GT}, there exists a
number $\beta >0$ such that $v(\mathbf{x,x}_{0})\geq \beta ,~\forall \mathbf{%
x}\in \overline{\Omega },$ $\forall \mathbf{x}_{0}\in E_{B}\left( a,b\right) 
$.

Let $\Gamma _{0}\subset \partial \Omega $ be the right side of the square $%
\Omega $ in (\ref{1.0}), 
\begin{equation}
\Gamma _{0}=\left\{ \mathbf{x}=\left( x,y\right) :x=a+c,\text{ }y\in
(b-c,b+c)\right\} ,  \label{2.06}
\end{equation}%
We formulate the coefficient inverse problem as follows:

\textbf{Coefficient Inverse Problem (CIP).} For each $\mathbf{x}_{0}\in
E_{B}\left( a,b\right) $ let $v(\mathbf{x},\mathbf{x}_{0})\in C^{3+\alpha
}\left( \overline{P_{D}}\left( a,b\right) \right) $ be the above solution of
problem (\ref{2-2}). Assume that the following functions $h_{0}(\mathbf{x},%
\mathbf{x}_{0})$ and $h_{1}(\mathbf{x},\mathbf{x}_{0})$ are given:%
\begin{equation}
\left. 
\begin{array}{c}
v(\mathbf{x},\mathbf{x}_{0})=h_{0}(\mathbf{x},\mathbf{x}_{0}),\text{ }%
\forall \mathbf{x}\in \partial \Omega , \\ 
v_{x}(\mathbf{x},\mathbf{x}_{0})=h_{1}(\mathbf{x},\mathbf{x}_{0}),~\forall 
\mathbf{x}\in \Gamma _{0}, \\ 
\forall x_{0}\in E_{B}(a,b).%
\end{array}%
\right.  \label{2.6}
\end{equation}%
Find the coefficient $\sigma \left( \mathbf{x}\right) $ in for $\mathbf{x}%
\in \Omega $.

\begin{figure}[]
\centering
\includegraphics[width=0.5\textwidth]{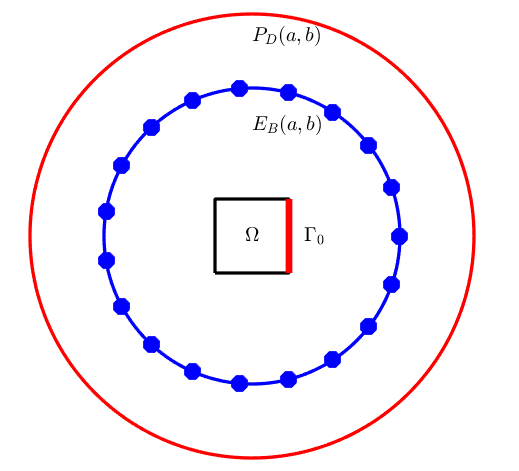}
\caption{A schematic diagram of our measurements. The large disk is the
domain $P_{D}\left( a,b\right) $ where the forward problem is solved. The
smaller circle is the circle $E_{B}\left( a,b\right) $ where point sources
are located$.$ The square is our domain of interest $\Omega ,$ where the CIP
is solved.}
\label{fig:Plot_Domain_02}
\end{figure}

As to the above CIP, two versions of the convexification method for it were
published in \cite{KEIT,EITIP2025}. Below we follow the version of \cite%
{EITIP2025}. First, the convexification method, transforms the inverse
problem in a system of two coupled elliptic PDEs with a viscosity term,
subject to the Cauchy boundary conditions. These PDEs do not contain the
unknown coefficient $\sigma \left( \mathbf{x}\right) .$ Within the framework
of the convexification method, a viscosity term was introduced for the first
time in \cite{KlibHJ} for the numerical solution of the Hamilton-Jacobi
equation. We also refer to \cite{LN} in this regard. In both \cite{EITIP2025}
and the current paper, an extended idea of \cite{KlibHJ,LN} about the
viscosity term in the convexification method is used. The presence of the
viscosity term enables us to avoid the truncation of a certain Fourier-like
series, such as the one in \cite{KEIT}. In addition, the viscosity term was
used in the convexification method for CIPs in \cite{Ktransp,KT3}.

\section{The Carleman Estimate-Based Convexification Method for Coefficient
Inverse Problem (\protect\ref{2-2})-(\protect\ref{2.6})}

\label{sec:3}

For a better understanding of the convexification related material of this
paper, we outline in this section the continuous version of this method,
which we take from \cite{EITIP2025}.

\subsection{Transformation procedure}

\label{sec:3.1}

First, we change variables as $V\left( \mathbf{x},\mathbf{x}_{0}\right) =%
\sqrt{\sigma \left( \mathbf{x}\right) }v\left( \mathbf{x},\mathbf{x}%
_{0}\right) .$ Using (\ref{2-2}) and (\ref{2.6}) we obtain 
\begin{equation}
\Delta V+r\left( \mathbf{x}\right) V=0,~\mbox{in}~\Omega ,  \label{3.1}
\end{equation}%
\begin{equation}
\left. 
\begin{array}{c}
V(\mathbf{x},\mathbf{x}_{0})=h_{0}(\mathbf{x},\mathbf{x}_{0}),\text{ }%
\mathbf{x}\in \partial \Omega , \\ 
V_{x}(\mathbf{x},\mathbf{x}_{0})=h_{1}(\mathbf{x},\mathbf{x}_{0}),\text{ }%
\forall \mathbf{x}=\left( a+c,y\right) \in \Gamma _{0},\text{ }y\in
(b-c,b+c), \\ 
\forall \mathbf{x}_{0}\in E_{B}\left( a,b\right) ,%
\end{array}%
\right.  \label{3.2}
\end{equation}%
\begin{equation}
r\left( \mathbf{x}\right) =\frac{\Delta \left( \sqrt{\sigma \left( \mathbf{x}%
\right) }\right) }{\sqrt{\sigma \left( \mathbf{x}\right) }}.  \label{3.3}
\end{equation}%
Let $\mathbf{x}_{0}=\mathbf{x}_{0}(\theta ),$ where the angle $\theta \in %
\left[ 0,2\pi \right] .$ Next, we set 
\begin{equation}
\psi (\mathbf{x},\theta )=\ln V(\mathbf{x},\theta )  \label{3.30}
\end{equation}%
and substitute this in (\ref{2.6}). We obtain%
\begin{equation}
\Delta \psi +\left( \nabla \psi \right) ^{2}+r\left( \mathbf{x}\right) =0,%
\mbox{ }\mathbf{x}\in \Omega .  \label{3.4}
\end{equation}
Note that a logarithmic transformation for a CIP was first used in \cite[the
next line after (1.14)]{Klib95}. Recalling (\ref{2.01}) and differentiating (%
\ref{3.4}) with respect to $\theta $ we eliminate the unknown coefficient $%
r\left( x\right) .$ If the function $\psi \left( x,\theta \right) $ is
found, then by (\ref{3.4}) 
\begin{equation}
r\left( \mathbf{x}\right) =-\frac{1}{2\pi }\int\limits_{0}^{2\pi }\left(
\Delta \psi +\left( \nabla \psi \right) ^{2}\right) \left( \mathbf{x},\theta
\right) d\theta ,\mbox{ }x\in \Omega .  \label{3.5}
\end{equation}%
Given $r\left( x\right) ,$ the question of finding the target coefficient $%
\sigma \left( \mathbf{x}\right) $ is addressed below.

The above differentiation with respect to $\theta $ leads to:%
\begin{equation}
\Delta \psi _{\theta }+2\nabla \psi _{\theta }\nabla \psi =0,\mbox{ }x\in
\Omega .  \label{3.05}
\end{equation}%
However, (\ref{3.05}) is one equation with two unknown functions $\psi $ and 
$\psi _{\theta }.$To express $\psi $ via $\psi _{\theta },$ one might try to
use%
\[
\psi \left( \mathbf{x},\theta \right) =\int\limits_{\theta _{0}}^{\theta
}\psi _{\theta }\left( \mathbf{x},s\right) ds+\psi \left( \mathbf{x},\theta
_{0}\right) 
\]%
for a certain $\theta _{0}\in \left( 0,2\pi \right) .$ However, the function 
$\psi \left( \mathbf{x},\theta _{0}\right) $ is unknown for any value of $%
\theta _{0}.$ This causes the introduction of a viscosity term $-\varepsilon
\Delta \psi $ with a small parameter $\varepsilon \in \left( 0,1\right) $ 
\cite{EITIP2025}. This term complements (\ref{3.05}) with the following
equation 
\begin{equation}
-\varepsilon \Delta \psi +\Delta \psi _{\theta }+2\nabla \psi _{\theta
}\nabla \psi =0,~\mathbf{x}\in \Omega .  \label{3.06}
\end{equation}%
Denote 
\begin{equation}
\phi =\psi _{\theta }-\varepsilon \psi .  \label{3.50}
\end{equation}%
Then (\ref{3.05}) and (\ref{3.06}) lead to the following system of coupled
nonlinear PDEs for the functions $\psi _{\theta }$ and $\phi $, 
\begin{equation}
F_{1}\left( \psi _{\theta },\phi \right) \left( \theta \right) =\Delta \psi
_{\theta }+2\nabla \psi _{\theta }\nabla \left( \frac{\psi _{\theta }-\phi }{%
\varepsilon }\right) =0,\text{ }\mathbf{x}\in \Omega ,\text{ }\theta \in %
\left[ 0,2\pi \right] ,  \label{3.19}
\end{equation}%
\begin{equation}
F_{2}\left( \psi _{\theta },\phi \right) \left( \theta \right) =\Delta \phi
+2\nabla \psi _{\theta }\nabla \left( \frac{\psi _{\theta }-\phi }{%
\varepsilon }\right) =0,\text{ }\mathbf{x}\in \Omega ,\text{ }\theta \in %
\left[ 0,2\pi \right] .  \label{3.20}
\end{equation}%
This system is subject to the overdetermined boundary conditions derived
from the data (\ref{3.2}) 
\begin{equation}
\psi _{\theta }\left( \mathbf{x},\theta \right) \mid _{\mathbf{x}\in
\partial \Omega }=\partial _{\theta }s_{0}\left( \mathbf{x},\theta \right) ,%
\mbox{ }\left( \mathbf{x},\theta \right) \in \partial \Omega \times \left[
0,2\pi \right] ,  \label{3.21}
\end{equation}%
\begin{equation}
\partial _{x}\left( \partial _{\theta }\psi \left( \mathbf{x},\theta \right)
\right) \mid _{\mathbf{x}\in \Gamma _{0}}=\partial _{\theta }s_{1}\left( 
\mathbf{x},\theta \right) ,\mbox{ }\forall \left( \mathbf{x},\theta \right)
\in \Gamma _{0}\times \left[ 0,2\pi \right] ,  \label{3.22}
\end{equation}%
\begin{equation}
\phi \left( \mathbf{x},\theta \right) \mid _{\mathbf{x}\in \partial \Omega
}=\partial _{\theta }s_{0}\left( \mathbf{x},\theta \right) -\varepsilon
s_{0}\left( \mathbf{x},\theta \right) ,\text{ }\forall \left( \mathbf{x}%
,\theta \right) \in \partial \Omega \times \left[ 0,2\pi \right] ,
\label{3.23}
\end{equation}%
\begin{equation}
\partial _{x}\phi \left( \mathbf{x},\theta \right) =\partial _{\theta
}s_{1}\left( \mathbf{x},\theta \right) -\varepsilon s_{1}\left( \mathbf{x}%
,\theta \right) ,\mbox{ }\forall \left( \mathbf{x},\theta \right) \in \Gamma
_{0}\times \left[ 0,2\pi \right] .  \label{3.24}
\end{equation}%
where $s_{0}$ and $s_{1}$ are derived from the data $h_{0}$ and $h_{1}$ (see
equations (3.9) and (3.10) in \cite{EITIP2025}). Introduce a new notation
for convenience 
\begin{equation}
(q,p)\left( \mathbf{x},\theta \right) =(\psi _{\theta },\phi )\left( \mathbf{%
x},\theta \right) .  \label{3.240}
\end{equation}%
The transformation procedure is complete.

\subsection{The Carleman estimate}

\label{sec:3.2}

The construction of the globally convergent convexification method for our
CIP is based on a Carleman estimate \cite{EITIP2025}. Denote 
\begin{equation}
\left. H_{0}^{2}\left( \Omega \right) =\left\{ u\in H^{2}\left( \Omega
\right) :u\mid _{\partial \Omega }=0,\partial _{x}u\mid _{\Gamma
_{0}}=0\right\} .\right.   \label{3.25}
\end{equation}%
Having (\ref{1.0}) and (\ref{2.06}), we define Carleman Weight Function
(CWF) as 
\begin{equation}
W_{\kappa }\left( \mathbf{x}\right) =\exp \left( 2\kappa x^{2}\right) ,
\label{3.26}
\end{equation}%
where $\kappa \geq 1$ is a parameter.

\textbf{Theorem 3.1} (Carleman estimate \cite[Theorem 10.3.1]{KL}).\emph{\
There exists a sufficiently large number }$\kappa _{0}=\kappa _{0}\left(
\Omega \right) \geq 1$\emph{\ and a number }$C=C\left( \Omega \right) >0,$ 
\emph{both numbers depending only on }$\Omega ,$\emph{\ such that with the
CWF }$W_{\kappa }\left( \mathbf{x}\right) $ \emph{as in (\ref{3.26}) the
following Carleman estimate holds:} 
\[
\int\limits_{\Omega }\left( \Delta u\right) ^{2}W_{\kappa }\left( \mathbf{x}%
\right) d\mathbf{x\geq }\frac{C}{\kappa }\int\limits_{\Omega }\left(
u_{xx}^{2}+u_{xy}^{2}+u_{yy}^{2}\right) W_{\kappa }\left( \mathbf{x}\right) d%
\mathbf{x+} 
\]%
\[
+C\kappa \int\limits_{\Omega }\left( \nabla u\right) ^{2}W_{\kappa }\left( 
\mathbf{x}\right) d\mathbf{x+}C\kappa ^{3}\int\limits_{\Omega
}u^{2}W_{\kappa }\left( \mathbf{x}\right) d\mathbf{x,} 
\]%
\[
\forall \kappa \geq \kappa _{0},\mbox{ }\forall u\in H_{0}^{2}\left( \Omega
\right) . 
\]

Everywhere below $C=C\left( \Omega \right) >0$ denotes different numbers
depending only on the domain $\Omega .$ We note that there is another CWF
for a general elliptic operator of the second order, see, e.g. \cite[section
2.4]{KL}, \cite[\S 1 of Chapter 4]{LRS}. However, that CWF depends on two
large parameters instead of one in (\ref{3.26}). The latter makes it hard to
use that CWF for computations.

\subsection{The globally strongly convex CWF-weighted Tikhonov-like
functional}

\label{sec:3.3}

The CWF (\ref{3.26}) is the weight function for our CWF-weighted
Tikhonov-like cost functional. We point out that even though $\kappa $ is
required to be sufficiently large, our rich computational experience of
working with the convexification method (see the above citations) shows that
actually $\kappa \in \left[ 1,5\right] $ is sufficient. In this regard,
there is a similarity here to an asymptotic theory. Indeed, typically such a
theory states that if a parameter $X$ is sufficiently large, then a certain
formula $Y$ is sufficiently accurate. However, in any specific computational
practice only numerical studies can establish which exactly values of $Y$
ensure a good accuracy of the formula $Y$.

Similarly with (\ref{3.25}) we define the space $H_{0}^{6}\left( \Omega
\right) $ as 
\begin{equation}
H_{0}^{6}\left( \Omega \right) =\left\{ u\in H^{6}\left( \Omega \right)
:u\mid _{\partial \Omega }=0,\partial _{x}u\mid _{\Gamma _{0}}=0\right\} .
\label{4.0}
\end{equation}%
Define spaces $H_{1,\theta }$ and $H_{1,2,\theta }$ depending on the
parameter $\theta \in \left[ 0,2\pi \right] ,$ 
\begin{equation}
\left. 
\begin{array}{c}
H_{1,\theta }=\left\{ 
\begin{array}{c}
q\left( \mathbf{x},\theta \right) :q\left( \mathbf{x},\theta \right) \in
H^{6}\left( \Omega \right) , \\ 
\left\Vert q\left( \mathbf{x},\theta \right) \right\Vert _{H_{1,\theta
}}^{2}=\left\Vert q\left( \mathbf{x},\theta \right) \right\Vert
_{H^{6}\left( \Omega \right) }^{2}<\infty ,\text{ }\forall \theta \in \left[
0,2\pi \right] ,%
\end{array}%
\right\} , \\ 
H_{1,\theta }^{0}=\left\{ q\left( \mathbf{x},\theta \right) :q\left( \mathbf{%
x},\theta \right) \in H_{0}^{6}\left( \Omega \right) \right\} , \\ 
H_{1,2,\theta }=H_{1,\theta }\times H_{1,\theta },\text{ }H_{1,2,\theta
}^{0}=H_{1,\theta }^{0}\times H_{1,\theta }^{0}, \\ 
\left\Vert \left( q,p\right) \left( \mathbf{x},\theta \right) \right\Vert
_{H_{1,2,\theta }}^{2}=\left\Vert q\left( \mathbf{x},\theta \right)
\right\Vert _{H_{1,\mathbf{\theta }}}^{2}+\left\Vert p\left( \mathbf{x}%
,\theta \right) \right\Vert _{H_{1,\theta }}^{2}.%
\end{array}%
\right.  \label{4.1}
\end{equation}%
We note that while in \cite[formula (4.1)]{EITIP2025} the norm $\left\Vert
q\left( x,\theta \right) \right\Vert _{H^{3}\left( \Omega \right) }$ was
used in the direct analog of (\ref{4.1}), here we use the norm $\left\Vert
q\left( x,\theta \right) \right\Vert _{H^{6}\left( \Omega \right) }.$ This
is because of the finite difference approximations of the second order
derivatives with the $O\left( h^{2}\right) -$accuracy, which we use below in
the discrete version of the convexification. Indeed, by Sobolev embedding
theorem 
\begin{equation}
\left. 
\begin{array}{c}
\left( q,p\right) \left( \mathbf{x},\theta \right) \in C^{4}\left( \overline{%
\Omega }\right) \times C^{4}\left( \overline{\Omega }\right) ,\text{ }%
\forall \left( q,p\right) \left( \mathbf{x},\theta \right) \in H_{1,2,\theta
}, \\ 
\left\Vert \left( q,p\right) \left( \mathbf{x},\theta \right) \right\Vert
_{C^{4}\left( \overline{\Omega }\right) \times C^{4}\left( \overline{\Omega }%
\right) }\leq C\left\Vert \left( q,p\right) \left( \mathbf{x},\theta \right)
\right\Vert _{H_{1,2,\theta }},\text{ } \\ 
\forall \left( q,p\right) \left( \mathbf{x},\theta \right) \in H_{1,2,\theta
},\text{ }\forall \theta \in \left[ 0,2\pi \right] .%
\end{array}%
\right.  \label{4.2}
\end{equation}

\textbf{Remark 3.1.} \emph{Regarding to the smoothness conditions in (\ref%
{4.0})-(\ref{4.2}), we note that, as a rule, extra smoothness requirements
are not considered as serious restrictions in the theory of CIPs, see, e.g. 
\cite{Nov,Rom}.}

Let $A>0$ be an arbitrary number, which we fix. Using notation (\ref{3.240}%
), we define the set $G_{\theta }\left( A\right) \subset $ $H_{1,2,\theta }$
as:%
\begin{equation}
G_{\theta }\left( A\right) =\left\{ 
\begin{array}{c}
\left( q,p\right) \left( \mathbf{x},\theta \right) \in H_{1,2,\theta
}:\left\Vert \left( q,p\right) \left( \mathbf{x},\theta \right) \right\Vert
_{H_{1,2,\theta }}<A, \\ 
q\mbox{ satisfies boundary conditions (\ref{3.21}), (\ref{3.22}),} \\ 
p\mbox{ satisfies boundary conditions (\ref{3.23}), (\ref{3.24}),} \\ 
\forall \theta \in \left[ 0,2\pi \right]%
\end{array}%
\right\} .  \label{4.3}
\end{equation}%
Hence, by (\ref{4.2}) and (\ref{4.3}) 
\begin{equation}
\left\Vert \left( q,p\right) \left( \mathbf{x},\theta \right) \right\Vert
_{C^{4}\left( \overline{\Omega }\right) \times C^{4}\left( \overline{\Omega }%
\right) }\leq CA,\forall \left( q,p\right) \left( \mathbf{x},\theta \right)
\in H_{1,2,\theta },\forall \theta \in \left[ 0,2\pi \right] .  \label{4.300}
\end{equation}

We define the functional 
\[
J_{\kappa ,\alpha }:\overline{G_{\theta }\left( A\right) }\rightarrow 
\mathbb{R},\text{ }\forall \theta \in \left[ 0,2\pi \right] 
\]%
as (see equations (4.1) and (4.2) in \cite{EITIP2025})%
\begin{equation}
J_{\kappa ,\alpha }\left( q,p\right) \left( \theta \right) =\sqrt{%
\varepsilon }\int\limits_{\Omega }\left[ \left( F_{1}\left( q,p\right)
\left( \mathbf{x},\theta \right) \right) ^{2}+\left( F_{2}\left( q,p\right)
\left( \mathbf{x},\theta \right) \right) ^{2}\right] W_{\kappa }\left( 
\mathbf{x}\right) d\mathbf{x}+  \label{4.4}
\end{equation}%
\[
+\alpha \left\Vert \left( q,p\right) \left( \mathbf{x},\theta \right)
\right\Vert _{H^{6}\left( \Omega \right) \times H^{6}\left( \Omega \right)
}^{2},\text{ }\forall \theta \in \left[ 0,2\pi \right] , 
\]%
where $\alpha \in (0,1)$ is the regularization parameter. We consider the
following minimization problem:

\textbf{Minimization Problem.} Minimize the functional $J_{\kappa ,\alpha
}\left( q,p\right) \left( \theta \right) $ on the set $\overline{G_{\theta
}\left( A\right) }$ for each $\theta \in \left[ 0,2\pi \right] .$

\textbf{Theorem 3.2} (global strong convexity) \cite{EITIP2025}. \emph{The
following hold true:}

\emph{1. For each }$\kappa >0$\emph{, for each }$\theta \in \left[ 0,2\pi %
\right] $\emph{\ and for each pair }$\left( q,p\right) \in \overline{%
G_{\theta }\left( A\right) }$\emph{\ the functional }$J_{\kappa ,\alpha
}\left( q,p\right) \left( \theta \right) $\emph{\ has the Fr\'{e}chet
derivative }%
\begin{equation}
J_{\kappa ,\alpha }^{\prime }\left( q,p\right) \left( \theta \right) \in
H_{1,2,\theta }^{0}\emph{,}\text{ }\forall \theta \in \left[ 0,2\pi \right] .
\label{5.01}
\end{equation}%
\emph{\ }

\emph{2. Let }$\kappa _{0}=\kappa _{0}\left( \Omega \right) \geq 1$\emph{\
be the number of Theorem 3.1. There exists a sufficiently large number }$%
\kappa _{1}=\kappa _{1}\left( A,\Omega ,\varepsilon \right) \geq \kappa _{0}$%
\emph{\ such that for each }$\kappa \geq \kappa _{1}$\emph{\ the functional }%
$J_{\kappa ,\alpha }\left( q,p\right) $\emph{\ is strongly convex on the set 
}$\overline{G_{\theta }\left( A\right) },$\emph{\ i.e. there exists a number 
}$C_{1}=C_{1}\left( A,\Omega ,\varepsilon \right) >0$\emph{\ such that}%
\[
J_{\kappa ,\alpha }\left( q_{2},p_{2}\right) \left( \theta \right)
-J_{\kappa ,\alpha }\left( q_{1},p_{1}\right) \left( \theta \right)
-J_{\kappa ,\alpha }^{\prime }\left( q_{1},p_{1}\right) \left( \theta
\right) \left( q_{2}-q_{1},p_{2}-p_{1}\right) \left( \mathbf{x},\theta
\right) \geq 
\]%
\begin{equation}
\geq C_{1}\exp \left( 2\kappa \left( a-c_{1}\right) ^{2}\right) \left\Vert
\left( q_{2}-q_{1},p_{2}-p_{1}\right) \left( \mathbf{x},\theta \right)
\right\Vert _{H^{2}\left( \Omega \right) \times H^{2}\left( \Omega \right)
}^{2}+  \label{5.3}
\end{equation}%
\[
+\alpha \left\Vert \left( q_{2}-q_{1},p_{2}-p_{1}\right) \left( \mathbf{x}%
,\theta \right) \right\Vert _{H_{1,2,\theta }}^{2}, 
\]%
\[
\forall \theta \in \left[ 0,2\pi \right] ,\mbox{ }\forall \left(
q_{1},p_{1}\right) \left( \mathbf{x},\theta \right) ,\left(
q_{2},p_{2}\right) \left( \mathbf{x},\theta \right) \in \overline{G_{\theta
}\left( A\right) },\mbox{ }\forall \kappa \geq \kappa _{1}. 
\]%
\emph{\ Numbers }$\kappa _{1}\left( A,\Omega ,\varepsilon \right) $\emph{\ \
and }$C_{1}\left( A,\Omega ,\varepsilon \right) $\emph{\ depend only on
listed parameters.}

\emph{3. For every }$\kappa \geq \kappa _{1}$\emph{\ and for every }$\theta
\in \left[ 0,2\pi \right] $ \emph{there exists unique minimizer }$\left(
q_{\min ,\kappa ,\alpha },p_{\min ,\kappa ,\alpha }\right) \left( \mathbf{x}%
,\theta \right) \in \overline{G_{\theta }\left( A\right) }$\emph{\ of the
functional }$J_{\kappa ,\alpha }\left( q,p\right) \left( \theta \right) $%
\emph{\ on the set }$\overline{G_{\theta }\left( A\right) }.$ \emph{In
addition,}%
\begin{equation}
J_{\kappa ,\alpha }^{\prime }\left( q_{\min ,\kappa ,\alpha },p_{\min
,\kappa ,\alpha }\right) \left( \theta \right) \left( \left( q_{\min ,\kappa
,\alpha }-q,p_{\min ,\kappa ,\alpha }-p\right) \left( \mathbf{x},\theta
\right) \right) \leq 0,  \label{5.4}
\end{equation}%
\[
\forall \theta \in \left[ 0,2\pi \right] ,\text{ }\forall \left( q,p\right)
\left( \mathbf{x},\theta \right) \in \overline{G_{\theta }\left( A\right) }. 
\]%
\emph{\ }

Everywhere below $C_{1}=C\left( A,\Omega ,\varepsilon \right) >0$ denotes
different numbers depending only on listed parameters. Suppose now that the
minimizer $\left( q_{\min ,\kappa ,\alpha },p_{\min ,\kappa ,\alpha }\right)
\left( \mathbf{x},\theta \right) $ of the functional $J_{\kappa ,\alpha
}\left( q,p\right) \left( \theta \right) $ on the set $\overline{G_{\theta
}\left( A\right) }$ is found. Next, keeping in mind (\ref{3.5}), (\ref{3.50}%
) and (\ref{3.240}), we set%
\begin{equation}
\partial _{\theta }\psi _{\min ,\kappa ,\alpha }\left( \mathbf{x},\theta
\right) =q_{\min ,\kappa ,\alpha }\left( \mathbf{x},\theta \right) ,\phi
_{\min ,\kappa ,\alpha }\left( \mathbf{x},\theta \right) =p_{\min ,\kappa
,\alpha }\left( \mathbf{x},\theta \right) ,  \label{5.5}
\end{equation}%
\begin{equation}
\psi _{\min ,\kappa ,\alpha }\left( \mathbf{x},\theta \right) =\frac{%
\partial _{\theta }\psi _{\min ,\kappa ,\alpha }\left( \mathbf{x},\theta
\right) -\phi _{\min ,\kappa ,\alpha }\left( \mathbf{x},\theta \right) }{%
\varepsilon },  \label{5.6}
\end{equation}%
\begin{equation}
r_{\kappa ,\alpha }\left( \mathbf{x}\right) =-\frac{1}{2\pi }%
\int\limits_{0}^{2\pi }\left( \Delta \psi _{\min ,\kappa ,\alpha }+\left(
\nabla \psi _{\min ,\kappa ,\alpha }\right) ^{2}\right) \left( \mathbf{x}%
,\theta \right) d\theta ,\mbox{
}\mathbf{x}\in \Omega .  \label{5.7}
\end{equation}%
Denote 
\begin{equation}
w\left( \mathbf{x}\right) =\sqrt{\sigma \left( \mathbf{x}\right) }.
\label{5.8}
\end{equation}%
By (\ref{3.3}) we now need to solve the following elliptic equation:%
\begin{equation}
\Delta w-r_{\kappa ,\alpha }\left( \mathbf{x}\right) w=0,\text{ }\mathbf{x}%
\in \Omega .  \label{5.9}
\end{equation}%
It follows from (\ref{1-1}), (\ref{2-1}) and (\ref{5.8}) that this equation
should be complemented by two boundary conditions%
\begin{equation}
w\mid _{\partial \Omega }=1,\partial _{\nu }w\mid _{\partial \Omega }=0,
\label{5.10}
\end{equation}%
where $\partial _{\nu }$ is the normal derivative at $\partial \Omega .$
Next, we assign the reconstructed function $\sigma _{\kappa ,\alpha }\left( 
\mathbf{x}\right) $ as: 
\begin{equation}
\sigma _{\kappa ,\alpha }\left( \mathbf{x}\right) =w^{2}\left( \mathbf{x}%
\right) .  \label{5.11}
\end{equation}

Problem (\ref{5.9}), (\ref{5.10}) has an over-determination in the boundary
conditions. Hence, this problem was solved in \cite{EITIP2025} by the
so-called Quasi-Reversibility Method, see \cite[formula (4.17) and section
4.3]{KL}. Specifically, a slightly modified version of \cite{KT3} of this
method was used. We do not discuss details here for brevity. We still
mention that the convergence of this method was proven in \cite{KL,KT3}
using Carleman estimates. More precisely, the following functional $I\left(
w\right) $ was minimized in \cite{EITIP2025} on the set of functions $w\in
H^{2}\left( \Omega \right) $ satisfying boundary conditions (\ref{5.10}) 
\begin{equation}
I\left( w\right) =\int\limits_{\Omega }\left( \Delta w-r_{\kappa ,\alpha
}\left( \mathbf{x}\right) w\right) ^{2}d\mathbf{x}.  \label{5.110}
\end{equation}

\subsection{The accuracy of the minimizer}

\label{sec:3.4}

The minimizer $\left( q_{\min ,\kappa ,\alpha },p_{\min ,\kappa ,\alpha
}\right) \left( \mathbf{x},\theta \right) $ of the functional $J_{\kappa
,\alpha }\left( q,p\right) \left( \theta \right) $ is called
\textquotedblleft regularized solution" in the theory of Ill-Posed Problems 
\cite{T}. Following one of conventional postulates of this theory, we assume
the existence of the true solution $\sigma ^{\ast }\left( \mathbf{x}\right) $
of our CIP satisfying conditions (\ref{1-1}), (\ref{2-1}) with the noiseless
data $h_{0}^{\ast }\left( \mathbf{x,}\theta \right) ,h_{1}^{\ast }\left( 
\mathbf{x,}\theta \right) $ in (\ref{2.6}). An important question is the
question of an estimate of the accuracy of the regularized solution, i.e. an
estimate of a norm of the difference between $\left( q_{\min ,\kappa ,\alpha
},p_{\min ,\kappa ,\alpha }\right) \left( \mathbf{x,}\theta \right) $ and $%
\left( q^{\ast },p^{\ast }\right) \left( \mathbf{x,}\theta \right) ,$ where
the last pair corresponds to $\sigma ^{\ast }\left( \mathbf{x}\right) .$ As
it follows from the discussion of the previous section, below Theorem 3.2,
this estimate, in turn leads to an estimate of a certain norm of the
difference $\sigma _{\kappa ,\alpha }\left( \mathbf{x}\right) -\sigma ^{\ast
}\left( \mathbf{x}\right) .$ Such an estimate was derived in \cite[Theorem
5.3]{EITIP2025}. To simplify the presentation, that estimate was derived in 
\cite[Theorem 5.3]{EITIP2025} for a norm of the difference $r_{\kappa
,\alpha }\left( \mathbf{x}\right) -r^{\ast }\left( \mathbf{x}\right) ,$
where $r^{\ast }\left( \mathbf{x}\right) $ is generated by $\sigma ^{\ast
}\left( \mathbf{x}\right) $ as in (\ref{3.3}).

While the case of noisy data was considered in \cite[Theorem 5.3]{EITIP2025}%
, we simplify the presentation here by assuming that our data (\ref{2.6})
are noiseless. This is because we only outline here the results of \cite%
{EITIP2025} without going to some details. Thus, we estimate here the
difference $r_{\kappa ,\alpha }\left( \mathbf{x}\right) -r^{\ast }\left( 
\mathbf{x}\right) $, assuming that 
\begin{equation}
h_{0}\left( \mathbf{x,}\theta \right) =h_{0}^{\ast }\left( \mathbf{x,}\theta
\right) ,\text{ }h_{1}\left( \mathbf{x,}\theta \right) =h_{1}^{\ast }\left( 
\mathbf{x,}\theta \right) .  \label{5.12}
\end{equation}

Since we now work in the framework of the viscosity solution, then, based on
(\ref{3.05}), (\ref{3.06}) and (\ref{3.240}), we assume that the pair $%
(q^{\ast },p^{\ast })\left( \mathbf{x},\theta \right) =(\psi _{\theta
}^{\ast },\phi ^{\ast })\left( \mathbf{x},\theta \right) $ satisfies
equations (\ref{3.19}) and (\ref{3.20}), i.e. 
\begin{equation}
F_{1}(q^{\ast },p^{\ast })\left( \mathbf{x},\theta \right) =F_{2}(q^{\ast
},p^{\ast })\left( \mathbf{x},\theta \right) =0,\text{ }\forall \mathbf{x}%
\in \Omega ,\text{ }\forall \theta \in \left[ 0,2\pi \right] .  \label{5.120}
\end{equation}

\textbf{Theorem 3.3:}

1. \emph{Let the function }$\sigma ^{\ast }\left( \mathbf{x}\right) $\emph{\
be the exact solution of our CIP satisfying conditions (\ref{1-1}), (\ref%
{2-1}) and with the noiseless data }$h_{0}^{\ast }\left( \mathbf{x,}\theta
\right) ,h_{1}^{\ast }\left( \mathbf{x,}\theta \right) $\emph{. Assume that
we have noiseless data in (\ref{2.6}), i.e. assume that (\ref{5.12}) holds.}

2. \emph{Also, let the function }$r^{\ast }\left( \mathbf{x}\right) $\emph{\
be linked with the function }$\sigma ^{\ast }\left( \mathbf{x}\right) $\emph{%
\ via the direct analog of formula (\ref{3.3}). Let }$(q^{\ast },p^{\ast
})\left( \mathbf{x},\theta \right) $ \emph{be the pair of functions\ \
generated by the function }$\sigma ^{\ast }\left( \mathbf{x}\right) .$\ 
\emph{Let }$G_{\theta }\left( A\right) $\emph{\ be the set of vector
functions defined in (\ref{4.3}). Assume that }%
\begin{equation}
(q^{\ast },p^{\ast })\left( \mathbf{x},\theta \right) \in G_{\theta }\left(
A\right) .  \label{5.121}
\end{equation}

3. \emph{Let }$\kappa _{1}$\emph{\ be the number of Theorem 3.2 and let }$%
\kappa \geq \kappa _{1}.$\emph{\ Let }

$\left( q_{\min ,\kappa ,\alpha },p_{\min ,\kappa ,\alpha }\right) \left( 
\mathbf{x},\theta \right) \in \overline{G_{\theta }\left( A\right) }$\emph{\
be the minimizer of the functional }$J_{\kappa ,\alpha }\left( q,p\right)
\left( \theta \right) $\emph{\ on the set }$\overline{G_{\theta }\left(
A\right) },$ \emph{which was found in Theorem 3.2. Let the function }$%
r_{\kappa ,\alpha }\left( x\right) $\emph{\ being found via (\ref{5.5})- (%
\ref{5.7}) be the corresponding approximation for function }$r^{\ast }\left( 
\mathbf{x}\right) .$\emph{\ Then the following accuracy estimates are valid:}%
\begin{equation}
\left\Vert \left( q_{\min ,\kappa ,\alpha }-q^{\ast }\right) \left( \mathbf{x%
},\theta \right) \right\Vert _{H^{2}\left( \Omega \right) }+\left\Vert
\left( p_{\min ,\kappa ,\alpha }-p^{\ast }\right) \left( \mathbf{x,}\theta
\right) \right\Vert _{H^{2}\left( \Omega \right) }\leq  \label{5.13}
\end{equation}%
\[
\leq C_{1}\sqrt{\alpha }\exp \left[ -\kappa \left( a-c\right) ^{2}\right] ,%
\text{ }\forall \theta \in \left[ 0,2\pi \right] , 
\]%
\begin{equation}
\left\Vert r_{\kappa ,\alpha }-r^{\ast }\right\Vert _{L_{2}\left( \Omega
\right) }\leq C_{1}\sqrt{\alpha }\exp \left[ -\kappa \left( a-c\right) ^{2}%
\right] .  \label{5.14}
\end{equation}

\textbf{Proof. }We provide this proof here since it is different from the
one of \cite[Theorem 5.3]{EITIP2025}. Recall (\ref{3.240}). The functional $%
J_{\kappa ,\alpha }\left( q,p\right) \left( \theta \right) $ in (\ref{4.4})
is the sum of two parts. The first part contains differential operators $%
F_{1}\left( q,p\right) \left( \theta \right) $ and $F_{2}\left( q,p\right)
\left( \theta \right) $ in (\ref{3.19}) and (\ref{3.20}), which are
generated by our CIP. And the second part contains only the regularization
term $\alpha \left\Vert \left( q,p\right) \left( \mathbf{x},\theta \right)
\right\Vert _{H_{1,2,\theta }}^{2},$ which is not directly linked to our
CIP. For each $\theta \in \left[ 0,2\pi \right] $ this can be written as:%
\begin{equation}
J_{1,\kappa ,\alpha }\left( q,p\right) \left( \theta \right) =\sqrt{%
\varepsilon }\int\limits_{\Omega }\left[ \left( F_{1}\left( q,p\right)
\left( \mathbf{x},\theta \right) \right) ^{2}+\left( F_{2}\left( q,p\right)
\left( \mathbf{x},\theta \right) \right) ^{2}\right] W_{\kappa }\left( 
\mathbf{x}\right) d\mathbf{x,}  \label{5.15}
\end{equation}%
\begin{equation}
J_{2,\kappa ,\alpha }\left( q,p\right) \left( \theta \right) =\alpha
\left\Vert \left( q,p\right) \left( \mathbf{x},\theta \right) \right\Vert
_{H_{1,2,\theta }}^{2},  \label{5.16}
\end{equation}%
\begin{equation}
J_{\kappa ,\alpha }\left( q,p\right) \left( \theta \right) =J_{1,\kappa
,\alpha }\left( q,p\right) \left( \theta \right) +J_{2,\kappa ,\alpha
}\left( q,p\right) \left( \theta \right) .  \label{5.17}
\end{equation}%
By (\ref{4.1}), (\ref{4.3}), (\ref{5.120}), (\ref{5.121}) and (\ref{5.15})-(%
\ref{5.17}) 
\begin{equation}
J_{\kappa ,\alpha }\left( q^{\ast },p^{\ast }\right) \left( \theta \right)
=\alpha \left\Vert \left( q^{\ast },p^{\ast }\right) \left( \mathbf{x}%
,\theta \right) \right\Vert _{H^{6}\left( \Omega \right) \times H^{6}\left(
\Omega \right) }^{2}\leq \alpha A,\text{ }\forall \theta \in \left[ 0,2\pi %
\right] .  \label{5.18}
\end{equation}%
Next, by (\ref{5.3})%
\[
J_{\kappa ,\alpha }\left( q^{\ast },p^{\ast }\right) \left( \theta \right)
-J_{\kappa ,\alpha }\left( q_{\min ,\kappa ,\alpha },p_{\min ,\kappa ,\alpha
}\right) \left( \theta \right) - 
\]%
\begin{equation}
-J_{\kappa ,\alpha }^{\prime }\left( \left( q_{\min ,\kappa ,\alpha
},p_{\min ,\kappa ,\alpha }\right) \right) \left( \theta \right) \left(
\left( q^{\ast }-q_{\min ,\kappa ,\alpha },p^{\ast }-p_{\min ,\kappa ,\alpha
}\right) \left( \mathbf{x},\theta \right) \right) \geq  \label{5.19}
\end{equation}%
\[
\geq C_{1}\exp \left( 2\kappa \left( a-c\right) ^{2}\right) \left\Vert
\left( q^{\ast }-q_{\min ,\kappa ,\alpha },p^{\ast }-p_{\min ,\kappa ,\alpha
}\right) \left( \mathbf{x},\theta \right) \right\Vert _{H^{2}\left( \Omega
\right) \times H^{2}\left( \Omega \right) }^{2}, 
\]%
\[
\forall \theta \in \left[ 0,2\pi \right] . 
\]

Using (\ref{5.4}), we obtain%
\begin{equation}
-J_{\kappa ,\alpha }^{\prime }\left( \left( q_{\min ,\kappa ,\alpha
},p_{\min ,\kappa ,\alpha }\right) \right) \left( \theta \right) \left(
q^{\ast }-q_{\min ,\kappa ,\alpha },p^{\ast }-p_{\min ,\kappa ,\alpha
}\right) \left( \mathbf{x},\theta \right) \leq 0.  \label{5.20}
\end{equation}%
In addition, 
\begin{equation}
-J_{\kappa ,\alpha }\left( \left( q_{\min ,\kappa ,\alpha },p_{\min ,\kappa
,\alpha }\right) \right) \left( \theta \right) \leq 0.  \label{5.21}
\end{equation}%
Hence, using (\ref{5.18})-(\ref{5.21}), we obtain 
\[
\left\Vert \left( q^{\ast }-q_{\min ,\kappa ,\alpha },p^{\ast }-p_{\min
,\kappa ,\alpha }\right) \left( \mathbf{x},\theta \right) \right\Vert
_{H^{2}\left( \Omega \right) \times H^{2}\left( \Omega \right) }^{2}\leq
C_{1}\alpha \exp \left( -2\kappa \left( a-c\right) ^{2}\right) , 
\]%
which implies the first target accuracy estimate (\ref{5.13}). The second
target accuracy estimate (\ref{5.14}) follows from (\ref{3.5}), (\ref{3.50}%
), (\ref{3.240}), (\ref{5.7}) and (\ref{5.13}). $\square $

\section{The Semi-Discrete Analogue of the Functional $J_{\protect\kappa ,%
\protect\alpha }$}

\label{sec:4}

We call the case of this section \textquotedblleft semi-discrete" because we
consider a discrete analog of the functional $J_{1,\kappa ,\alpha }\left(
q,p\right) \left( \theta \right) $ in (\ref{5.15}), whereas we consider the
functional $J_{2,\kappa ,\alpha }\left( q,p\right) \left( \theta \right) $
in (\ref{5.16}) in its original continuous form. We need the material of
this section since we apply the convexification method on a coarse grid on
the first stage of our numerical procedure: due to a slow performance of
this method on a fine grid.

\subsection{Discretization}

\label{sec:4.1}

Recall definition (\ref{1.0}) of the domain $\Omega $ and consider the
following grid $\left\{ \left( x_{i},y_{j}\right) \right\}
_{i,j=0}^{n+1}\subset \overline{\Omega }$ in $\overline{\Omega }$ with the
grid step size $h\in \left( 0,1\right) $ and $n>4$%
\begin{equation}
\left. 
\begin{array}{c}
a-c=x_{0}<x_{1}<\ldots <x_{n}<x_{n+1}=a+c,\quad x_{i}-x_{i-1}=h, \\ 
b-c=y_{0}<y_{1}<\ldots <y_{n}<y_{n+1}=b+c,\quad y_{j}-y_{j-1}=h.%
\end{array}%
\right.  \label{6.00}
\end{equation}%
Everywhere below $O\left( h^{k}\right) $ denotes different functions such
that 
\begin{equation}
\left\vert O\left( h^{k}\right) \right\vert \leq C_{1}h^{k},k>0.
\label{6.01}
\end{equation}%
Recall that $\Gamma _{0}\subset \partial \Omega $ was defined in (\ref{2.06}%
). Denote 
\begin{equation}
\left. 
\begin{array}{c}
\Omega ^{h}=\left\{ \left( x_{i},y_{j}\right) \right\} _{i,j=1}^{n}, \\ 
\overline{\Omega ^{h}}=\left\{ \left( x_{i},y_{j}\right) \right\}
_{i,j=0}^{n+1}, \\ 
\partial \Omega ^{h}=\overline{\Omega ^{h}}\diagdown \Omega ^{h}, \\ 
\Gamma _{0}^{h}=\left\{ \left( a+c,y_{j}\right) ,\text{ }j=0,...,n+1\right\}
.%
\end{array}%
\right.  \label{6.1}
\end{equation}

Denote the running point in $\overline{\Omega ^{h}}$ as $\mathbf{x}%
_{i,j}^{h}=\left( x_{i},y_{j}\right) ,$ $\forall i,j=0,...,n+1.$ Any
function $z\left( \mathbf{x}\right) \in C\left( \overline{\Omega }\right) $
generates the discrete function $z^{h}\left( \mathbf{x}_{i,j}^{h}\right) $
defined on the grid (\ref{6.00}). Hence, 
\begin{equation}
z^{h}\left( \mathbf{x}_{i,j}^{h}\right) =z\left( x_{i},y_{j}\right)
,i,j=0,...,n+1.  \label{6.001}
\end{equation}

Let the function $f\left( \mathbf{x}\right) \in C^{4}\left( \overline{\Omega 
}\right) .$ We define its first derivative with respect to $x$ in finite
differences at the points $\left( x_{i},y_{j}\right) \in \Omega ^{h},$ $%
i,j=1,...,n$ as: 
\begin{equation}
\partial _{x}^{h}f_{i,j}=\partial _{x}^{h}f^{h}\left( x_{i},y_{j}\right) =%
\frac{f\left( x_{i+1},y_{j}\right) -f\left( x_{i-1},y_{j}\right) }{2h}.
\label{6.2}
\end{equation}%
Next, its second derivative with respect to $x$ in finite differences at the
points $\left( x_{i},y_{j}\right) \in \Omega ^{h},$ $i,j=1,...,n$ $\ $is
defined as: 
\begin{equation}
\partial _{x}^{h,2}f_{i,j}=\partial _{x}^{h,2}f^{h}\left( x_{i},y_{j}\right)
=\frac{f\left( x_{i+1},y_{j}\right) -2f\left( x_{i},y_{j}\right) +f\left(
x_{i-1},y_{j}\right) }{h^{2}}.  \label{6.4}
\end{equation}%
In addition, because of the Neumann boundary conditions (\ref{3.22}) and (%
\ref{3.24}) at $\Gamma _{0},$ we define its finite difference approximation
at points $f\left( x_{n+1},y_{k}\right) =f\left( a+c,y_{k}\right) $ as: 
\begin{equation}
\left. 
\begin{array}{c}
\partial _{x}^{h}f_{n+1,j}=\partial _{x}^{h}f^{h}\left( a+c,y_{j}\right) =
\\ 
=\left[ 3f\left( a+c,y_{j}\right) -4f\left( x_{n},y_{j}\right) +f\left(
x_{n-1},y_{j}\right) \right] /\left( 2h\right) , \\ 
\left( a+c,y_{j}\right) \in \Gamma _{0}^{h},\text{ }j=0,...,n+1.%
\end{array}%
\right.  \label{6.6}
\end{equation}%
Note that%
\begin{equation}
\left. 
\begin{array}{c}
\partial _{x}^{h}f_{i,j}^{h}=f_{x}\left( x_{i},y_{j}\right) +O\left(
h^{2}\right) ,\text{ }\partial _{x}^{h,2}f_{i,j}=f_{xx}\left(
x_{i},y_{j}\right) +O\left( h^{2}\right) , \\ 
i,j=1,...,n, \\ 
\partial _{x}^{h}f_{n+1,k}=\partial _{x}^{h}f^{h}\left( a+c,y_{k}\right)
=f_{x}\left( a+c,y_{k}\right) +O\left( h^{2}\right) ,k=0,...,n+1, \\ 
\left\vert O\left( h^{2}\right) \right\vert \leq C\left\Vert f\right\Vert
_{C^{4}\left( \overline{\Omega }\right) }h^{2}, \\ 
\forall f\in C^{4}\left( \overline{\Omega }\right) .%
\end{array}%
\right.  \label{6.60}
\end{equation}%
First and second finite difference derivatives $\partial _{x}f_{i,j}^{h}$
and $\partial _{y}^{2}f_{i,j}^{h}$ with respect to $y$ as well as the mixed
derivative $\partial _{xy}^{h,2}f_{i,j}^{h}$ at the points $\left(
x_{i},y_{j}\right) \in \Omega ^{h}$ are defined similarly, and analogs of
formulas (\ref{6.60}) hold. Next, we consider finite difference analogs of
the Laplace operator and the gradient vector, 
\begin{equation}
\Delta ^{h}f_{i,j}=\partial _{x}^{h,2}f_{i,j}+\partial
_{y}^{h,2}f_{i,j},~\nabla ^{h}f_{ij}=\left( \partial
_{x}^{h}f_{i,j},\partial _{y}^{h}f_{i,j}\right) .  \label{6.7}
\end{equation}%
Given the above constructions, Proposition 3.1 follows immediately from (\ref%
{6.001})-(\ref{6.7}) and Taylor formula.

\textbf{Proposition 4.1.} \emph{For each function }$f\in C^{4}\left( 
\overline{\Omega }\right) $\emph{\ the following accuracy estimates are
valid for points }$\mathbf{x}_{i,j}^{h}\in \Omega ^{h}$ \emph{as} $%
h\rightarrow 0^{+}$%
\[
\left. 
\begin{array}{c}
\partial _{x}^{h,2}f_{i,j}=f_{xx}\left( x_{i},y_{j}\right) +O\left(
h^{2}\right) ,\text{ }\partial _{y}^{h,2}f_{i,j}=f_{yy}\left(
x_{i},y_{j}\right) +O\left( h^{2}\right) , \\ 
\partial _{xy}^{h,2}f_{i,j}^{h}=f_{xy}\left( x_{i},y_{j}\right) +O\left(
h^{2}\right) , \\ 
\nabla ^{h}f_{i,j}=\nabla f\left( x_{i},y_{j}\right) +O\left( h^{2}\right) ,
\\ 
\partial _{x}f_{n+1,k}^{h}=f_{x}\left( a+c,x_{2,k}\right) +O\left(
h^{2}\right) , \\ 
\left\vert O\left( h^{2}\right) \right\vert \leq C\left\Vert f\right\Vert
_{C^{4}\left( \overline{\Omega }\right) }h^{2}, \\ 
i,j=1,...,n;\text{ }k=0,...,n+1.%
\end{array}%
\right. 
\]%
\emph{\ }

We have for any function $f\in C^{1}\left( \overline{\Omega }\right) $: 
\begin{equation}
\left. \int\limits_{\Omega }f\left( \mathbf{x}\right) dx\mathbf{=}%
\int\limits_{\Omega ,rt}f\left( \mathbf{x}_{ij}^{h}\right) d\mathbf{x}%
_{ij}^{h}\mathbf{+}O\left( h\right) ,\text{ }\left\vert O\left( h\right)
\right\vert \leq C\left\Vert f\right\Vert _{C^{1}\left( \overline{\Omega }%
\right) }h,\right.  \label{6.9}
\end{equation}%
where the symbol 
\[
\int\limits_{\Omega ,rt}f\left( \mathbf{x}_{ij}^{h}\right) d\mathbf{x}%
_{ij}^{h} 
\]%
denotes the approximation of the regular integral in (\ref{6.9}) by the
rectangle method, i.e. 
\begin{equation}
\int\limits_{\Omega ,rt}f\left( \mathbf{x}_{ij}^{h}\right) d\mathbf{x}%
_{ij}^{h}=h^{2}\sum_{i,j=0}^{n+1}f^{2}\left( x_{i},y_{j}\right) ,\text{ }%
\forall f\in C^{1}\left( \overline{\Omega }\right) .  \label{6.90}
\end{equation}%
Even though the rectangle method works for a more general set of functions $%
f\left( \mathbf{x}\right) ,$ it is convenient for the goal of this paper to
restrict our attention in (\ref{6.90}) only to functions $f\in C^{1}\left( 
\overline{\Omega }\right) .$

\subsection{Some function spaces}

\label{sec:4.2}

First, we introduce some discrete analogs of two function spaces $%
L_{2}\left( \Omega \right) $ and $H^{2}\left( \Omega \right) .$We set for
all functions $f\left( \mathbf{x}\right) $ defined at all points of the
domain $\Omega :$%
\begin{equation}
L_{2}^{h}\left( \Omega ^{h}\right) =\left\{ f^{h}:\left\Vert
f^{h}\right\Vert _{L_{2}^{h}\left( \Omega ^{h}\right)
}^{2}=\int\limits_{\Omega ,rt}f^{2}\left( \mathbf{x}_{ij}^{h}\right) d%
\mathbf{x}_{ij}^{h}=h^{2}\sum\limits_{i,j=0}^{n+1}f^{2}\left(
x_{i},y_{j}\right) \right\} .  \label{6.10}
\end{equation}%
\begin{equation}
H^{2,h}\left( \Omega ^{h}\right) =\left\{ 
\begin{array}{c}
f^{h}:\left\Vert f^{h}\right\Vert _{H^{2,h}\left( \Omega ^{h}\right)
}^{2}=\left\Vert f^{h}\right\Vert _{L_{2}^{h}\left( \Omega ^{h}\right) }^{2}+
\\ 
+h^{2}\sum\limits_{i,j=1}^{n}\left[ \left( \partial _{x}^{h}f_{i,j}\right)
^{2}+\left( \partial _{y}^{h}f_{i,j}\right) ^{2}\right] + \\ 
+h^{2}\sum\limits_{i,j=1}^{n}\left[ \left( \partial _{x}^{h,2}f_{i,j}\right)
^{2}+\left( \partial _{y}^{h,2}f_{i,j}\right) ^{2}+\left( \partial
_{yy}^{2}f_{i,j}^{h}\right) ^{2}\right]%
\end{array}%
\right\} .  \label{6.11}
\end{equation}

\ Clearly, sums in (\ref{6.11}) are approximate values of integrals over $%
\Omega $ of corresponding functions calculated by the rectangle rule (\ref%
{6.90}). Hence, by Proposition 4.1 and (\ref{6.9})-(\ref{6.12}) 
\begin{equation}
\left. 
\begin{array}{c}
\left\Vert f^{h}\right\Vert _{L_{2}^{h}\left( \Omega ^{h}\right)
}^{2}=\left\Vert f\right\Vert _{L_{2}\left( \Omega \right) }^{2}+O\left(
h\right) \mbox{ as }h\rightarrow 0^{+},\forall f\in C^{1}\left( \overline{%
\Omega }\right) , \\ 
\left\vert O\left( h\right) \right\vert \leq C\left\Vert f\right\Vert
_{C^{1}\left( \overline{\Omega }\right) }^{2}h.%
\end{array}%
\right.  \label{6.12}
\end{equation}%
\begin{equation}
\left. 
\begin{array}{c}
\left\Vert f^{h}\right\Vert _{H^{2,h}\left( \Omega ^{h}\right)
}^{2}=\left\Vert f\right\Vert _{H^{2}\left( \Omega \right) }^{2}+O\left(
h\right) \mbox{ as }h\rightarrow 0^{+},\text{ }\forall f\in C^{4}\left( 
\overline{\Omega }\right) , \\ 
\left\vert O\left( h\right) \right\vert \leq C\left\Vert f\right\Vert
_{C^{4}\left( \overline{\Omega }\right) }^{2}h.%
\end{array}%
\right.  \label{6.13}
\end{equation}

Recalling that the spaces $H_{1,\theta }$ and $H_{1,2,\theta }$ are defined
in (\ref{4.1}), introduce six more function spaces: 
\begin{equation}
H_{1,2,\theta ,fd}=\left\{ 
\begin{array}{c}
\left( q^{h},p^{h}\right) \left( \mathbf{x}_{i,j}^{h},\theta \right)
:q^{h},p^{h}\in H^{2,h}\left( \Omega ^{h}\right) , \\ 
\left\Vert \left( q^{h},p^{h}\right) \left( \mathbf{x}_{i,j}^{h},\theta
\right) \right\Vert _{H_{1,2,\theta ,fd}}^{2}= \\ 
=\left\Vert q^{h}\left( \mathbf{x}_{i,j}^{h},\theta \right) \right\Vert
_{H^{2,h}\left( \Omega ^{h}\right) }^{2}+\left\Vert p^{h}\left( \mathbf{x}%
_{i,j}^{h},\theta \right) \right\Vert _{H^{2,h}\left( \Omega ^{h}\right)
}^{2}<\infty , \\ 
\forall \theta \in \left[ 0,2\pi \right] ,%
\end{array}%
\right\} ,  \label{6.141}
\end{equation}%
\begin{equation}
\left. H_{1,2,2,\theta ,fd}=\left\{ 
\begin{array}{c}
\left( q^{h},p^{h},q,p\right) : \\ 
\left( q^{h},p^{h}\right) \in H_{1,2,\theta ,fd},\text{ }\left( q,p\right)
\in H_{1,2,\theta }, \\ 
\left( q^{h},p^{h}\right) \left( x_{i},y_{j},\theta \right) =\left(
q,p\right) \left( x_{i},y_{j},\theta \right) , \\ 
i,j=0,...,n+1, \\ 
\left\Vert \left( q^{h},p^{h},q,p\right) \right\Vert _{H_{1,2,2,\theta
,fd}}^{2}= \\ 
=\left\Vert \left( q^{h},p^{h}\right) \right\Vert _{H_{1,2,\theta
,fd}}^{2}+\left\Vert \left( q,p\right) \right\Vert _{H_{1,2,\theta }}^{2},
\\ 
\forall \theta \in \left[ 0,2\pi \right] ,%
\end{array}%
\right\} \right.  \label{6.142}
\end{equation}%
\begin{equation}
H_{1,2,\theta ,fd}^{0}=\left( q^{h},p^{h}\right) :\left\{ 
\begin{array}{c}
q^{h}\left( \mathbf{x}_{i,j}^{h},\theta \right) \mid _{\partial \Omega
^{h}}=p^{h}\left( \mathbf{x}_{i,j}^{h},\theta \right) \mid _{\partial \Omega
^{h}}=0, \\ 
i,j=0,...,n+1, \\ 
\partial _{x}q_{n+1,k}^{h}=\partial _{x}p_{n+1,k}^{h}=0,k=0,...,n+1, \\ 
\mbox{see (\ref{6.6}) for }\partial _{x}q_{n+1,k}^{h},\partial
_{x}p_{n+1,k}^{h}, \\ 
\left\Vert \cdot \right\Vert _{H_{1,2,\theta ,fd}^{0}}=\left\Vert \cdot
\right\Vert _{H_{1,2,\theta ,fd}}%
\end{array}%
\right\} ,  \label{6.143}
\end{equation}%
\begin{equation}
\left. H_{1,2,2,\theta ,fd}^{0}=\left\{ 
\begin{array}{c}
\left( q^{h},p^{h},q,p\right) \in H_{1,2,2,\theta ,fd}: \\ 
\left( q^{h},p^{h}\right) \in H_{1,2,\theta ,fd}^{0},\left( q,p\right) \in
H_{1,2,\theta }^{0}, \\ 
\left\Vert \cdot \right\Vert _{H_{1,2,2,\theta ,fd}^{0}}=\left\Vert \cdot
\right\Vert _{H_{1,2,2,\theta ,fd}},\text{ }\forall \theta \in \left[ 0,2\pi %
\right]%
\end{array}%
\right\} ,\right.  \label{6.145}
\end{equation}%
\begin{equation}
H_{1,2,\theta }^{2}=\left\{ 
\begin{array}{c}
\left( q,p\right) \left( \mathbf{x},\theta \right) \in H^{2}\left( \Omega
\right) \times H^{2}\left( \Omega \right) : \\ 
\left\Vert \left( q,p\right) \left( \mathbf{x},\theta \right) \right\Vert
_{H_{1,2,\theta }^{2}}^{2}=\left\Vert q\left( \mathbf{x},\theta \right)
\right\Vert _{H^{2}\left( \Omega \right) }^{2}+\left\Vert p\left( \mathbf{x}%
,\theta \right) \right\Vert _{H^{2}\left( \Omega \right) }^{2}, \\ 
\forall \theta \in \left[ 0,2\pi \right] .%
\end{array}%
\right\}  \label{6.200}
\end{equation}%
\begin{equation}
\left. H_{1,2,2,\theta ,fd}^{2}=\left\{ 
\begin{array}{c}
\left( q^{h},p^{h},q,p\right) : \\ 
\left( q^{h},p^{h}\right) \in H_{1,2,\theta ,fd},\text{ }\left( q,p\right)
\in H_{1,2,\theta }^{2}, \\ 
\left\Vert \left( q^{h},p^{h},q,p\right) \right\Vert _{H_{1,2,2,\theta
,fd}^{2}}= \\ 
=\left\Vert \left( q^{h},p^{h}\right) \right\Vert _{H_{1,2,\theta
,fd}}^{2}+\left\Vert \left( q,p\right) \right\Vert _{H_{1,2,\theta }^{2}}^{2}%
\text{ }\forall \theta \in \left[ 0,2\pi \right] .%
\end{array}%
\right\} ,\right.  \label{6.1460}
\end{equation}

Recall that the space $H_{1,2,\theta }$ was defined in (\ref{4.1}). Keeping
in mind (\ref{6.1}), we define now the semi-discrete analog $G_{\theta
}^{h}\left( A\right) $ of the set $G_{\theta }\left( A\right) $ in (\ref{4.3}%
) as:

\textbf{Definition. }$G_{\theta }^{h}\left( A\right) $\emph{\ is the set of
all vector functions }$\left( q^{h},p^{h},q,p\right) \in H_{1,2,2,\theta
,fd} $\emph{\ such that }$\left( q,p\right) \in G_{\theta }\left( A\right) .$%
\emph{\ Boundary conditions (\ref{3.22}), (\ref{3.24}) for }

$\left( q^{h},p^{h}\right) $ at $\left( x_{n+1},y_{j}\right) =\left(
a+c,y_{j}\right) $\emph{\ are understood in terms of (\ref{6.6}).}

Hence, by (\ref{6.12}), (\ref{6.13}), (\ref{6.141}), (\ref{6.200}) and (\ref%
{6.1460}) 
\begin{equation}
G_{\theta }^{h}\left( A\right) \subset H_{1,2,2,\theta ,fd},  \label{6.146}
\end{equation}%
\begin{equation}
\left. 
\begin{array}{c}
\left\Vert \left( q^{h},p^{h}\right) \right\Vert _{H_{1,2,\theta
,fd}}^{2}\left( \theta \right) =\left\Vert \left( q,p\right) \right\Vert
_{H_{1,2,\theta }^{2}}^{2}\left( \theta \right) +O\left( h\right) , \\ 
\forall \left( q^{h},p^{h},q,p\right) \in G_{\theta }^{h}\left( A\right) ,%
\text{ }\forall \theta \in \left[ 0,2\pi \right] .%
\end{array}%
\right.  \label{6.0146}
\end{equation}

Thus, using (\ref{3.21})-(\ref{3.24}), (\ref{3.240}), (\ref{4.3}), (\ref%
{6.001}), (\ref{6.1}), (\ref{6.6}), Definition, (\ref{6.141}), (\ref{6.142})
and (\ref{6.146}), we obtain%
\begin{equation}
\left. 
\begin{array}{c}
q^{h}\left( x_{i},y_{j},\theta \right) \mid _{\left( x_{i},y_{j}\right) \in
\partial \Omega ^{h}}=\partial _{\theta }s_{0}\left( x_{i},y_{j},\theta
\right) , \\ 
\partial _{x}^{h}q^{h}\left( x_{n+1},y_{j},\theta \right) =\partial _{\theta
}s_{1}\left( x_{n+1},y_{j},\theta \right) ,\text{ }\left(
x_{n+1},y_{j}\right) \in \Gamma _{0}, \\ 
p^{h}\left( x_{i},y_{j},\theta \right) \mid _{\left( x_{i},y_{j}\right) \in
\partial \Omega ^{h}}=\left( \partial _{\theta }s_{0}-\varepsilon
s_{0}\right) \left( x_{i},y_{j},\theta \right) , \\ 
\partial _{x}^{h}p^{h}\left( x_{n+1},y_{j},\theta \right) =\left( \partial
_{\theta }s_{1}-\varepsilon s_{1}\right) \left( x_{n+1},y_{j},\theta \right)
,\left( x_{n+1},y_{j}\right) \in \Gamma _{0}, \\ 
i,j=0,...,n+1, \\ 
\forall \left( q^{h},p^{h},q,p\right) \in G_{\theta }^{h}\left( A\right) ,%
\text{ }\forall \theta \in \left[ 0,2\pi \right] ,%
\end{array}%
\right.  \label{6.014}
\end{equation}%
where $\partial _{x}^{h}q^{h}\left( x_{n+1},y_{j},\theta \right) $ and $%
\partial _{x}^{h}p^{h}\left( x_{n+1},y_{j},\theta \right) $ are understood
in terms of (\ref{6.6}).

\subsection{Semi-discretization of the functional $J_{\protect\kappa ,%
\protect\alpha }$}

\label{sec:4.3}

First, using (\ref{6.141}), we consider the discrete analog $J_{1,h,\kappa
,\alpha }$ of the functional $J_{1,\kappa ,\alpha }$ in (\ref{5.15}),%
\[
J_{1,h,\kappa ,\alpha }\left( q^{h},p^{h}\right) \left( \theta \right) = 
\]%
\begin{equation}
=\sqrt{\varepsilon }\int\limits_{\Omega ,rt}\left[ \left( F_{1}^{h}\left(
q^{h},p^{h}\right) \left( \mathbf{x}_{i,j}^{h},\theta \right) \right)
^{2}+\left( F_{2}^{h}\left( q^{h},p^{h}\right) \left( \mathbf{x}%
_{i,j}^{h},\theta \right) \right) ^{2}\right] W_{\kappa }^{h}\left( \mathbf{x%
}_{i,j}^{h}\right) d\mathbf{x}_{i,j}^{h}\mathbf{,}  \label{6.15}
\end{equation}%
\[
\forall \left( q^{h},p^{h}\right) \left( \mathbf{x}_{i,j}^{h},\theta \right)
\in H_{1,2,\theta ,fd},\text{ }\forall \theta \in \left[ 0,2\pi \right] . 
\]%
Operators $F_{1}^{h},F_{2}^{h}$ in (\ref{6.15}) are finite difference
versions of operators $F_{1},F_{2}$ in (\ref{3.19}), (\ref{3.20}). Next, the
functional 
\begin{equation}
\left. J_{2,\kappa ,\alpha }\left( q,p\right) \left( \theta \right) \text{
remains the same as in (\ref{5.16}).}\right.  \label{6.150}
\end{equation}%
Now, using (\ref{6.142}), (\ref{6.146}), (\ref{6.15}) and (\ref{6.150}) and
Definition, we define the semi-discrete analog of the functional $J_{\kappa
,\alpha }\left( q,p\right) \left( \theta \right) $ similarly with (\ref{5.17}%
) as:%
\begin{equation}
\left. 
\begin{array}{c}
J_{h,\kappa ,\alpha }\left( q^{h},p^{h},q,p\right) \left( \theta \right)
=J_{1,h,\kappa ,\alpha }\left( q^{h},p^{h}\right) \left( \theta \right)
+J_{2,\kappa ,\alpha }\left( q,p\right) \left( \theta \right) , \\ 
\forall \left( q^{h},p^{h},q,p\right) \in G_{\theta }^{h}\left( A\right) ,%
\text{ }\forall \theta \in \left[ 0,2\pi \right] .%
\end{array}%
\right.  \label{6.16}
\end{equation}

\section{Theorems About the Semi-Discrete Functional $J_{h,\protect\kappa ,%
\protect\alpha }$}

\label{sec:5}

\subsection{The $h-$strong convexity of $J_{h,\protect\kappa ,\protect\alpha %
}$}

\label{sec:5.1}

Below $\left\{ ,\right\} $ and $\left[ ,\right] $ are scalar products in $%
H_{1,2,2,\theta ,fd}$ \ and $H_{1,2,\theta }$ respectively.

\textbf{Theorem 5.1} (the $h-$strong convexity). \emph{The following hold
true:}

\emph{1. For each }$\kappa >0$\emph{, for each }$\theta \in \left[ 0,2\pi %
\right] $\emph{\ and for each vector function }$\left(
q^{h},p^{h},q,p\right) \in \overline{G_{\theta }^{h}\left( A\right) }$\emph{%
\ the functional }$J_{h,\kappa ,\alpha }\left( q^{h},p^{h},q,p\right) \left(
\theta \right) $\emph{\ has the Fr\'{e}chet derivative }%
\begin{equation}
J_{h,\kappa ,\alpha }^{\prime }\left( q^{h},p^{h},q,p\right) \left( \theta
\right) \in H_{1,2,2,\theta ,fd}^{0},  \label{6.17}
\end{equation}%
\emph{where }$H_{1,2,2,\theta ,fd}^{0}$\emph{\ is the function space defined
in (\ref{6.143}).\ Furthermore,}%
\[
J_{h,\kappa ,\alpha }^{\prime }\left( q_{1}^{h},p_{1}^{h},q_{1},p_{1}\right)
\left( \theta \right) \left(
q_{2}^{h}-q_{1}^{h},p_{2}^{h}-p_{1}^{h},q_{2}-q_{1},p_{2}-p_{1}\right)
\left( x_{i,j}^{h},\theta \right) = 
\]%
\begin{equation}
=J_{\kappa ,\alpha }^{\prime }\left( q_{1},p_{1}\right) \left( \theta
\right) \left( q_{2}-q_{1},p_{2}-p_{1}\right) +O\left( h\right) \exp \left(
2\kappa \left( a+c\right) ^{2}\right) ,  \label{6.170}
\end{equation}%
\[
\forall \left( q_{1}^{h},p_{1}^{h},q_{1},p_{1}\right) ,\left(
q_{2}^{h},p_{2}^{h},q_{2},p_{2}\right) \mathbf{\in }G_{\theta }^{h}\left(
A\right) ,\text{ }\forall \theta \in \left[ 0,2\pi \right] , 
\]%
\emph{where} $J_{\kappa ,\alpha }^{\prime }\left( q,p\right) \left( \theta
\right) $\emph{\ is the Fr\'{e}chet derivative of the functional }$J_{\kappa
,\alpha }\left( q,p\right) \left( \theta \right) ,$\emph{\ which was found
in Theorem 3.2.}

\emph{2. Let }$\kappa _{0}=\kappa _{0}\left( \Omega \right) \geq 1$\emph{\
and }$\kappa _{1}=\kappa _{1}\left( A,\Omega ,\varepsilon \right) \geq
\kappa _{0}$ \emph{be the numbers of Theorems 3.1 and 3.2 respectively. Then
there exists a sufficiently large number }$\kappa _{2}=\kappa _{2}\left(
A,\Omega ,\varepsilon \right) \geq \kappa _{1}$\emph{\ such that for each }$%
\kappa \geq \kappa _{2}$\emph{\ and for each }$\theta \in \left[ 0,2\pi %
\right] $ \emph{the functional }$J_{h,\kappa ,\alpha }\left(
q^{h},p^{h},q,p\right) \left( \theta \right) $\emph{\ is }$h-$\emph{strongly
convex on the set }$\overline{G_{\theta }^{h}\left( A\right) }$,\emph{\ i.e.}%
\[
J_{h,\kappa ,\alpha }\left( q_{2}^{h},p_{2}^{h},q_{2},p_{2}\right) \left(
\theta \right) -J_{h,\kappa ,\alpha }\left(
q_{1}^{h},p_{1}^{h},q_{1},p_{1}\right) \left( \theta \right) - 
\]%
\[
-J_{h,\kappa ,\alpha }^{\prime }\left(
q_{1}^{h},p_{1}^{h},q_{1},p_{1}\right) \left( \theta \right) \left(
q_{2}^{h}-q_{1}^{h},p_{2}^{h}-p_{1}^{h},q_{2}-q_{1},p_{2}-p_{1}\right) \geq 
\]%
\begin{equation}
\geq C_{1}\exp \left( 2\kappa \left( a-c\right) ^{2}\right) \left\Vert
\left( q_{2}^{h}-q_{1}^{h},p_{2}^{h}-p_{1}^{h}\right) \right\Vert
_{H_{1,2,\theta ,fd}}^{2}+  \label{6.18}
\end{equation}%
\[
+\alpha \left\Vert \left( q_{2}-q_{1},p_{2}-p_{1}\right) \right\Vert
_{H_{1,2,\theta }}^{2}-C_{1}h\exp \left( 2\kappa \left( a+c\right)
^{2}\right) , 
\]%
\[
\forall \left( q_{1}^{h},p_{1}^{h},q_{1},p_{1}\right) ,\left(
q_{2}^{h},p_{2}^{h},q_{2},p_{2}\right) \in \overline{G_{\theta }^{h}\left(
A\right) }, 
\]%
\[
\forall \theta \in \left[ 0,2\pi \right] ,\text{ }\forall \kappa \geq \kappa
_{2}. 
\]%
\emph{\ }

\emph{3. Assume that there exists a pair of numbers }$\left( \kappa
_{3},\theta _{0}\right) $\emph{\ satisfying }$\kappa _{3}\geq \kappa _{2}$%
\emph{, }$\theta _{0}\in \left[ 0,2\pi \right] $\emph{\ and such that there
exists a minimizer }%
\begin{equation}
X_{\min }=\left( \widetilde{q}_{\min ,\kappa _{3},\alpha }^{h},\widetilde{p}%
_{\min ,\kappa _{3},\alpha }^{h},\widetilde{q}_{\min ,\kappa _{3},\alpha },%
\widetilde{p}_{\min ,\kappa _{3},\alpha }\right) \in \overline{G_{\theta
_{0}}^{h}\left( A\right) }  \label{6.41}
\end{equation}%
\emph{of the functional\ }$J_{h,\kappa _{3},\alpha }\left(
q^{h},p^{h},q,p\right) \left( \theta _{0}\right) $\emph{\ on the set }$%
\overline{G_{\theta _{0}}^{h}\left( A\right) }.$ \emph{In other words,
assume that}%
\begin{equation}
\left. 
\begin{array}{c}
J_{h,\kappa _{3},\alpha }\left( X_{\min }\right) \left( \theta _{0}\right)
\leq J_{h,\kappa ,\alpha }\left( q^{h},p^{h},q,p\right) \left( \theta
_{0}\right) , \\ 
\forall \left( q^{h},p^{h},q,p\right) \left( \theta _{0}\right) \in 
\overline{G_{\theta _{0}}^{h}\left( A\right) }.%
\end{array}%
\right.  \label{6.019}
\end{equation}%
\emph{Then\ }%
\begin{equation}
\hspace{-1cm}\left. 
\begin{array}{c}
J_{h,\kappa _{3},\alpha }^{\prime }\left( X_{\min }\right) \left( X_{\min
}-\left( q^{h},p^{h},q,p\right) \right) \left( \theta _{0}\right) \leq 0, \\ 
\forall \left( q^{h},p^{h},q,p\right) \in \overline{G_{\theta
_{0}}^{h}\left( A\right) }.%
\end{array}%
\right.  \label{6.19}
\end{equation}

\textbf{Remark 5.1.} Note that while the strong convexity estimate (%
\ref{5.3}) of Theorem 3.2 for the continuous functional $J_{\kappa ,\alpha }
$ has only two positive terms in the right hand side of (\ref{5.3}%
), a negative term $-C_{1}h\exp \left( 2\kappa \left( a+c\right)
^{2}\right) $ is added in the right hand side of (\ref{6.18}). This
is exactly what the term \textquotedblleft $h-$strongly convexity"
means. Fortunately this term is small, as long as $h$ is small.
This allows us to estimate some distances, as stated in Theorems 5.2 and
5.3. On the other hand, because of the presence of this negative term, we
cannot prove the existence and uniqueness of the minimizer of the functional 
$J_{h,\kappa _{3},\alpha }$ on the set $\overline{G_{\theta
_{0}}^{h}\left( A\right) },$\ which is unlike Theorem 3.2. Rather,
we only assume its existence, see (\ref{6.41}). Still, it follows from
Theorem 5.2 that even if there exist several minimizers of the functional $%
J_{h,\kappa _{3},\alpha }$ on the set $\overline{G_{\theta
_{0}}^{h}\left( A\right) },$ they are close to each other as long
as $\sqrt{\alpha }+\sqrt{h}$ is sufficiently small. And the latter
is sufficient for the computational purpose. In this regard we refer to
subsections 1.2 and 1.3.

\textbf{Proof of Theorem 5.1.} For every $\theta \in \left[ 0,2\pi \right] $
define the set\textbf{\ }$G_{\theta ,0}^{h}\left( 2A\right) $ similarly with
the set $G_{\theta }^{h}\left( A\right) $ in Definition. More precisely,%
\begin{equation}
\left. 
\begin{array}{c}
G_{\theta ,0}^{h}\left( 2A\right) \text{ is defined the same way as }%
G_{\theta }^{h}\left( A\right) \text{ in Definition,} \\ 
\text{except that }A\text{ is replaced with }2A, \\ 
\text{and the right hand sides of boundary conditions} \\ 
\text{(\ref{3.21})-(\ref{3.24}) are replaced with zeros.}%
\end{array}%
\right.  \label{6.20}
\end{equation}%
Then (\ref{6.141})-(\ref{6.146}) and (\ref{6.20}) imply 
\begin{equation}
\overline{G_{\theta ,0}^{h}\left( 2A\right) }\subset H_{1,2,2,\theta
,fd}^{0}.  \label{6.21}
\end{equation}

Let\textbf{\ }$\left( q_{1}^{h},p_{1}^{h},q_{1},p_{1}\right) \mathbf{\in }%
\overline{G_{\theta }^{h}\left( A\right) }$ and $\left(
q_{2}^{h},p_{2}^{h},q_{2},p_{2}\right) \mathbf{\in }\overline{G_{\theta
}^{h}\left( A\right) }$ be two arbitrary pairs. Denote 
\begin{equation}
\left( m^{h},z^{h},m,z\right) =\left( q_{2}^{h},p_{2}^{h},q_{2},p_{2}\right)
-\left( q_{1}^{h},p_{1}^{h},q_{1},p_{1}\right) .  \label{6.22}
\end{equation}%
Then, using (\ref{6.20}), we obtain 
\begin{equation}
\left( m^{h},z^{h},m,z\right) \in \overline{G_{\theta ,0}^{h}\left(
2A\right) }.  \label{6.23}
\end{equation}

By (\ref{3.19}), (\ref{3.240}) and (\ref{6.7}) \textbf{\ }%
\begin{equation}
F_{1}^{h}\left( q^{h},p^{h}\right) \left( x_{i},y_{j},\theta \right) =\Delta
^{h}q_{i,j}\left( \theta \right) +2\nabla ^{h}q_{i,j}\left( \theta \right)
\nabla ^{h}\left( \frac{q_{i,j}-p_{i,j}}{\varepsilon }\right) \left( \theta
\right) .  \label{6.230}
\end{equation}%
Hence, using (\ref{6.20}), we obtain%
\[
\left. 
\begin{array}{c}
F_{1}^{h}\left( \left( q_{2}^{h},p_{2}^{h}\right) \left( x_{i},y_{j},\theta
\right) \right) =F_{1}^{h}\left( \left(
q_{1}^{h}+m^{h},p_{1}^{h}+z^{h}\right) \left( x_{i},y_{j},\theta \right)
\right) = \\ 
=F_{1}^{h}\left( \left( q_{1}^{h},p_{1}^{h}\right) \left( x_{i},y_{j},\theta
\right) \right) + \\ 
+\left( \Delta ^{h}m_{i,j}\left( \theta \right) +2\nabla ^{h}q_{1,i,j}\left(
\theta \right) \nabla ^{h}\left( \left( m_{i,j}-z_{i,j}\right) /\varepsilon
\right) \right) + \\ 
+2\nabla ^{h}m_{i,j}\left( \theta \right) \nabla ^{h}\left( \left(
q_{1,i,j}-p_{1,i,j}\right) /\varepsilon \right) \left( \theta \right) + \\ 
+2\nabla ^{h}m_{i,j}\left( \theta \right) \nabla ^{h}\left( \left(
m_{i,j}-z_{i,j}\right) /\varepsilon \right) \left( \theta \right) .%
\end{array}%
\right. 
\]%
Hence, 
\begin{equation}
\left. 
\begin{array}{c}
\left[ F_{1}^{h}\left( \left( q_{2}^{h},p_{2}^{h}\right) \left(
x_{i},y_{j},\theta \right) \right) \right] ^{2}-\left[ F_{1}^{h}\left(
\left( q_{1}^{h},p_{1}^{h}\right) \left( x_{i},y_{j},\theta \right) \right) %
\right] ^{2}= \\ 
=L_{1,\mbox{lin},h}\left( \left( m^{h},z^{h}\right) \left(
x_{i},y_{j},\theta \right) \right) + \\ 
+L_{1,\mbox{nonlin},h}\left( \left( m^{h},z^{h}\right) \left(
x_{i},y_{j},\theta \right) \right) .%
\end{array}%
\right.  \label{6.24}
\end{equation}%
In (\ref{6.24}) terms $L_{1,\mbox{lin},h}\left( \left( m^{h},z^{h}\right)
\left( x_{i},y_{j},\theta \right) \right) $ and

$L_{1,\mbox{nonlin},h}\left( \left( m^{h},z^{h}\right) \left(
x_{i},y_{j},\theta \right) \right) $ are linear and nonlinear ones
respectively with respect to $\left( m^{h},z^{h}\right) \left(
x_{i},y_{j},\theta \right) .$ More precisely, 
\begin{equation}
\left. 
\begin{array}{c}
L_{1,\mbox{lin},h}\left( \left( m^{h},z^{h}\right) \left( x_{i},y_{j},\theta
\right) \right) =2F_{1}^{h}\left( \left( q_{1}^{h},p_{1}^{h}\right) \left(
x_{i},y_{j},\theta \right) \right) \times \\ 
\times \left[ 
\begin{array}{c}
\Delta ^{h}m_{i,j}\left( \theta \right) +2\nabla ^{h}q_{1,i,j}\left( \theta
\right) \nabla ^{h}\left( \left( m_{i,j}-z_{i,j}\right) /\varepsilon \right)
\left( \theta \right) + \\ 
+2\nabla ^{h}m_{i,j}\left( \theta \right) \nabla ^{h}\left( \left(
q_{1,i,j}-p_{1,i,j}\right) /\varepsilon \right) \left( \theta \right)%
\end{array}%
\right] .%
\end{array}%
\right.  \label{6.25}
\end{equation}%
\begin{equation}
\left. 
\begin{array}{c}
L_{1,\mbox{nonlin},h}\left( \left( m^{h},z^{h}\right) \left(
x_{i},y_{j},\theta \right) \right) = \\ 
=\left[ 
\begin{array}{c}
\left( \Delta ^{h}m_{i,j}\left( \theta \right) +2\nabla ^{h}q_{1,i,j}\left(
\theta \right) \right) \nabla ^{h}\left( \left( m_{i,j}-z_{i,j}\right)
/\varepsilon \right) \left( \theta \right) + \\ 
+2\nabla ^{h}m_{i,j}\left( \theta \right) \nabla ^{h}\left( \left(
q_{1,i,j}-p_{1,i,j}\right) /\varepsilon \right) \left( \theta \right) + \\ 
+2\nabla ^{h}m_{i,j}\left( \theta \right) \nabla ^{h}\left( \left(
m_{i,j}-z_{i,j}\right) /\varepsilon \right) \left( \theta \right)%
\end{array}%
\right] ^{2}+ \\ 
+4F_{1}^{h}\left( \left( q_{1}^{h},p_{1}^{h}\right) \left(
x_{i},y_{j},\theta \right) \right) \cdot \nabla ^{h}m_{i,j}\left( \theta
\right) \nabla ^{h}\left( \left( m_{i,j}-z_{i,j}\right) /\varepsilon \right)
\left( \theta \right) .%
\end{array}%
\right.  \label{6.26}
\end{equation}%
Using (\ref{6.90}), (\ref{6.20})-(\ref{6.23}) and (\ref{6.26}), we obtain%
\[
\sqrt{\varepsilon }\int\limits_{\Omega ,rt}\left\vert L_{1,\mbox{nonlin}%
,h}\left( \left( m^{h},z^{h}\right) \left( \mathbf{x}_{i,j}^{h},\theta
\right) \right) \right\vert W_{\kappa }^{h}\left( \mathbf{x}%
_{i,j}^{h}\right) d\mathbf{x}_{i,j}^{h}\leq 
\]%
\begin{equation}
\leq C_{1}\exp \left( 2\kappa \left( a+c\right) ^{2}\right) \left\Vert
\left( m^{h},z^{h}\right) \left( \mathbf{x}_{i,j}^{h},\theta \right)
\right\Vert _{H_{1,2,\theta ,fd}}^{2}.  \label{6.27}
\end{equation}

Obviously an expression, which is similar with the one in (\ref{6.24}), is
valid for the term $\left[ F_{2}^{h}\left( \left( q_{2}^{h},p_{2}^{h}\right)
\left( x_{i},y_{j},\theta \right) \right) \right] ^{2}-\left[
F_{2}^{h}\left( q_{1}^{h},p_{1}^{h}\right) \left( x_{i},y_{j},\theta \right) %
\right] ^{2}$. In this case

$L_{1,\mbox{lin},h}\left( \left( m^{h},z^{h}\right) \left(
x_{i},y_{j},\theta \right) \right) $ and $L_{1,\mbox{nonlin},h}\left( \left(
m^{h},z^{h}\right) \left( x_{i},y_{j},\theta \right) \right) $ are

replaced with slightly different operators $L_{2,\mbox{lin},h}\left( \left(
m^{h},z^{h}\right) \left( x_{i},y_{j},\theta \right) \right) $ and

$L_{2,\mbox{nonlin},h}\left( \left( m^{h},z^{h}\right) \left(
x_{i},y_{j},\theta \right) \right) $ respectively, and analogs of formulas

(\ref{6.25})-(\ref{6.27}) are valid.

Thus, (\ref{6.15}), (\ref{6.150}), (\ref{6.23}) and (\ref{6.25})-(\ref{6.27}%
) imply%
\[
J_{h,\kappa ,\alpha }\left( q_{2}^{h},p_{2}^{h},q_{2},p_{2}\right) \left(
\theta \right) -J_{h,\kappa ,\alpha }\left(
q_{1}^{h},p_{1}^{h},q_{1},p_{1}\right) \left( \theta \right) = 
\]%
\[
=\sqrt{\varepsilon }\int\limits_{\Omega _{rt}}L_{1,\mbox{lin},h}\left(
\left( m^{h},z^{h}\right) \left( \mathbf{x}_{i,j}^{h},\theta \right) \right)
W_{\kappa }^{h}\left( \mathbf{x}_{i,j}^{h}\right) d\mathbf{x}_{i,j}^{h}+ 
\]%
\[
+\sqrt{\varepsilon }\int\limits_{\Omega _{rt}}L_{2,\mbox{lin},h}\left(
\left( m^{h},z^{h}\right) \left( \mathbf{x}_{i,j}^{h},\theta \right) \right)
W_{\kappa }^{h}\left( \mathbf{x}_{i,j}^{h}\right) d\mathbf{x}_{i,j}^{h}+ 
\]%
\begin{equation}
+2\left[ \left( \left( q_{1},p_{1}\right) ,\left( m,z\right) \right) \left( 
\mathbf{x},\theta \right) \right] +  \label{6.28}
\end{equation}%
\[
+\sqrt{\varepsilon }\int\limits_{\Omega _{rt}}L_{1,\mbox{nonlin},h}\left(
\left( m^{h},z^{h}\right) \left( \mathbf{x}_{i,j}^{h},\theta \right) \right)
W_{\kappa }^{h}\left( \mathbf{x}_{i,j}^{h}\right) d\mathbf{x}_{i,j}^{h}+ 
\]%
\[
+\sqrt{\varepsilon }\int\limits_{\Omega _{rt}}L_{2,\mbox{nonlin},h}\left(
\left( m^{h},z^{h}\right) \left( \mathbf{x}_{i,j}^{h},\theta \right) \right)
W_{\kappa }^{h}\left( \mathbf{x}_{i,j}^{h}\right) d\mathbf{x}_{i,j}^{h}+ 
\]%
\[
+\left\Vert \left( m,z\right) \left( \mathbf{x},\theta \right) \right\Vert
_{H_{1,2,\theta }}^{2},\text{ }\forall \theta \in \left[ 0,2\pi \right]
,\forall \kappa \geq \kappa _{2}. 
\]

Keeping in mind (\ref{6.141})-(\ref{6.145}) and (\ref{6.20})-(\ref{6.23}),
consider the functional%
\[
P^{h}\left( \theta \right) :H_{1,2,2,\theta ,fd}^{0}\rightarrow \mathbb{R}, 
\]%
\[
P\left( t_{1}^{h},t_{2}^{h},t_{1},t_{2}\right) \left( \theta \right) =\sqrt{%
\varepsilon }\int\limits_{\Omega _{rt}}L_{1,\mbox{lin},h}\left( \left(
t_{1}^{h},t_{2}^{h}\right) \left( \mathbf{x}_{i,j}^{h},\theta \right)
\right) W_{\kappa }^{h}\left( \mathbf{x}_{i,j}^{h}\right) d\mathbf{x}%
_{i,j}^{h}+ 
\]%
\begin{equation}
+\sqrt{\varepsilon }\int\limits_{\Omega _{rt}}L_{2,\mbox{lin},h}\left(
\left( t_{1}^{h},t_{2}^{h}\right) \left( \mathbf{x}_{i,j}^{h},\theta \right)
\right) W_{\kappa }^{h}\left( \mathbf{x}_{i,j}^{h}\right) d\mathbf{x}%
_{i,j}^{h}+  \label{6.29}
\end{equation}%
\[
+2\left[ \left( q_{1},p_{1}\right) ,\left( t_{1},t_{2}\right) \right] , 
\]%
\[
\forall \left( t_{1}^{h},t_{2}^{h},t_{1},t_{2}\right) \in H_{1,2,2,\theta
,fd}^{0},\text{ }\forall \theta \in \left[ 0,2\pi \right] . 
\]%
This is a bounded linear functional defined on the Hilbert space $%
H_{1,2,2,\theta ,fd}^{0}.$ Hence, by Riesz theorem for each $\theta \in %
\left[ 0,2\pi \right] $ there exists unique point $\widehat{P}^{h}\left(
\theta \right) \in H_{1,2,2,\theta ,fd}^{0}$ \ such that%
\begin{equation}
P^{h}\left( t_{1}^{h},t_{2}^{h},t_{1},t_{2}\right) \left( \theta \right)
=\left\{ \widehat{P}^{h}\left( \theta \right) ,\left(
t_{1}^{h},t_{2}^{h},t_{1},t_{2}\right) \right\} ,  \label{6.30}
\end{equation}%
\[
\forall \left( t_{1}^{h},t_{2}^{h},t_{1},t_{2}\right) \in H_{1,2,2,\theta
,fd}^{0},\text{ }\forall \theta \in \left[ 0,2\pi \right] . 
\]%
It follows from (\ref{4.1}), (\ref{6.141})-(\ref{6.145}), (\ref{6.22}) and (%
\ref{6.27})-(\ref{6.30}) that%
\[
\left. 
\begin{array}{c}
\left\vert 
\begin{array}{c}
J_{h,\kappa ,\alpha }\left( q_{2}^{h},p_{2}^{h},q_{2},p_{2}\right) \left(
\theta \right) -J_{h,\kappa ,\alpha }\left(
q_{1}^{h},p_{1}^{h},q_{1},p_{1}\right) \left( \theta \right) - \\ 
-\left\{ \widehat{P}\left( \theta \right) ,\left(
q_{2}^{h}-q_{1}^{h},p_{2}^{h}-p_{1}^{h},q_{2}-q_{1},p_{2}-p_{1}\right)
\left( \mathbf{x}_{i,j}^{h},\theta \right) \right\}%
\end{array}%
\right\vert \leq \\ 
\leq C_{1}\left\Vert \left(
q_{2}^{h}-q_{1}^{h},p_{2}^{h}-p_{1}^{h},q_{2}-q_{1},p_{2}-p_{1}\right)
\right\Vert _{H_{1,2,2,\theta ,fd}}^{2},\text{ }\forall \theta \in \left[
0,2\pi \right] .%
\end{array}%
\right. 
\]%
Hence, $\widehat{P}^{h}\left( \theta \right) =J_{h,\kappa ,\alpha }^{\prime
}\left( q_{1}^{h},p_{1}^{h},q_{1},p_{1}\right) \left( \theta \right) \in
H_{1,2,2,\theta ,fd}^{0}$ \ is the Fr\'{e}chet derivative of the functional%
\emph{\ }$J_{h,\kappa ,\alpha }\left( q_{1}^{h},p_{1}^{h},q_{1},p_{1}\right)
\left( \theta \right) $ at an arbitrary point $\left(
q_{1}^{h},p_{1}^{h},q_{1},p_{1}\right) \left( \theta \right) \mathbf{\in }%
\overline{G_{\theta }^{h}\left( A\right) },$ which proves (\ref{6.17}).

We now prove (\ref{6.170}). Using (\ref{3.26}) and (\ref{6.90}), we obtain%
\[
\sqrt{\varepsilon }\int\limits_{\Omega _{rt}}L_{1,\mbox{lin},h}\left( \left(
m^{h},z^{h}\right) \left( \mathbf{x}_{i,j}^{h},\theta \right) \right)
W_{\kappa }^{h}\left( \mathbf{x}_{i,j}^{h}\right) d\mathbf{x}_{i,j}^{h}= 
\]%
\begin{equation}
=\sqrt{\varepsilon }h^{2}\sum_{i,j=1}^{n}2F_{1}\left( \left(
q_{1},p_{1}\right) \left( x_{1,i},x_{2,j},\theta \right) \right) \times
\label{6.32}
\end{equation}%
\[
\times \left[ 
\begin{array}{c}
\Delta ^{h}m+2\nabla ^{h}q_{1}\nabla ^{h}\left( \left( m-z\right)
/\varepsilon \right) + \\ 
+2\nabla ^{h}m\nabla ^{h}\left( \left( q_{1}-p_{1}\right) /\varepsilon
\right) \left( \theta \right)%
\end{array}%
\right] \left( x_{i},y_{j},\theta \right) W_{\kappa }\left( x_{i}\right) . 
\]%
Hence, (\ref{6.01}), (\ref{6.9}), Proposition 4.1 and (\ref{6.32}) imply%
\[
\sqrt{\varepsilon }\int\limits_{\Omega _{rt}}L_{1,\mbox{lin},h}\left( \left(
m^{h},z^{h}\right) \left( \mathbf{x}_{i,j}^{h},\theta \right) \right)
W_{\kappa }^{h}\left( \mathbf{x}_{i,j}^{h}\right) d\mathbf{x}_{i,j}^{h}= 
\]%
\begin{equation}
=\sqrt{\varepsilon }\int\limits_{\Omega }2F_{1}\left( \left(
q_{1},p_{1}\right) \left( \mathbf{x},\theta \right) \right) \times
\label{6.33}
\end{equation}%
\[
\times \left[ \Delta m+2\nabla q_{1}\nabla \left( \left( m-z\right)
/\varepsilon \right) +2\nabla m\nabla \left( \left( q_{1}-p_{1}\right)
/\varepsilon \right) \right] \left( \mathbf{x},\theta \right) W_{\kappa
}\left( \mathbf{x}\right) d\mathbf{x+} 
\]%
\[
+O\left( h\right) \exp \left( 2\kappa \left( a+c\right) ^{2}\right) ,\text{ }%
\forall \theta \in \left[ 0,2\pi \right] . 
\]%
A similar formula and in a similar way can also be obtained for 
\begin{equation}
\sqrt{\varepsilon }\int\limits_{\Omega _{rt}}L_{2,\mbox{lin},h}\left( \left(
m^{h},z^{h}\right) \left( \mathbf{x}_{i,j}^{h},\theta \right) \right)
W_{\kappa }^{h}\left( \mathbf{x}_{i,j}^{h}\right) d\mathbf{x}_{i,j}^{h}.
\label{6.34}
\end{equation}%
Comparison of (\ref{6.33}) and its analog for (\ref{6.34}) with the
analogous formula (6.3) for the continuous case in \cite{EITIP2025} shows
that (\ref{6.170}) holds true.

Acting similarly, we obtain%
\begin{equation}
\left. 
\begin{array}{c}
J_{h,\kappa ,\alpha }\left( q^{h},p^{h},q,p\right) \left( \theta \right)
=J_{\kappa ,\alpha }\left( q,p\right) \left( \theta \right) +O\left(
h\right) \exp \left( 2\kappa \left( a+c\right) ^{2}\right) , \\ 
\forall \left( q^{h},p^{h},q,p\right) \in G_{\theta }^{h}\left( A\right) ,%
\text{ }\forall \theta \in \left[ 0,2\pi \right] .%
\end{array}%
\right.  \label{6.35}
\end{equation}%
Comparing (\ref{6.170}) and (\ref{6.35}) with the continuous case of (\ref%
{5.3}), we obtain the desired $h-$strong convexity property (\ref{6.18}).

We now prove (\ref{6.19}). Suppose that the opposite is true. In other
words, assume that there exists a point $\left( q^{h},p^{h},q,p\right) \in
G_{\theta _{0}}^{h}\left( A\right) $ such that 
\begin{equation}
\hspace{-1cm}\left. J_{h,\kappa _{3},\alpha }^{\prime }\left( X_{\min
}\right) \left( X_{\min }-\left( q^{h},p^{h},q,p\right) \right) \left(
\theta _{0}\right) >0.\right.  \label{6.36}
\end{equation}%
Denote%
\begin{equation}
Y=\left( q^{h},p^{h},q,p\right) \in G_{\theta _{0}}^{h}\left( A\right) .
\label{6.37}
\end{equation}%
Then (\ref{6.41}), (\ref{6.36}) and (\ref{6.37}) imply%
\begin{equation}
J_{h,\kappa ,\alpha }^{\prime }\left( X_{\min }\right) \left( X_{\min
}-Y\right) >0.  \label{6.38}
\end{equation}%
Consider the interval of the straight line connecting points $\widetilde{Z}%
_{\min }$ and $Y$. This interval is%
\[
I=\left\{ \left( 1-\beta \right) X_{\min }+\beta Y,\beta \in \left[ 0,1%
\right] \right\} . 
\]
Since the set $\overline{G_{\theta _{0}}^{h}\left( A\right) }$ is convex,
then $I\in \overline{G_{\theta _{0}}^{h}\left( A\right) }.$ This is
equivalent with 
\[
X_{\min }+\beta \left( Y-X_{\min }\right) \in \overline{G_{\theta
_{0}}^{h}\left( A\right) },\text{ }\forall \beta \in \left[ 0,1\right] . 
\]%
We have for sufficiently small $\beta >0:$%
\begin{equation}
\left. 
\begin{array}{c}
J_{h,\kappa _{3},\alpha }\left( X_{\min }+\beta \left( Y-X_{\min }\right)
\right) \left( \theta _{0}\right) =J_{h,\kappa _{3},\alpha }\left( X_{\min
}\right) \left( \theta _{0}\right) + \\ 
+J_{h,\kappa _{3},\alpha }^{\prime }\left( X_{\min }\right) \left( \beta
\left( Y-X_{\min }\right) \right) \left( \theta _{0}\right) +o(1),\text{ as }%
\beta \rightarrow 0^{+}.%
\end{array}%
\right.  \label{6.39}
\end{equation}%
Using (\ref{6.38}), we obtain 
\[
J_{h,\kappa _{3},\alpha }^{\prime }\left( X_{\min }\right) \left( \beta
\left( Y-X_{\min }\right) \right) \left( \theta _{0}\right) =\beta
J_{h,\kappa _{3},\alpha }^{\prime }\left( X_{\min }\right) \left( Y-X_{\min
}\right) \left( \theta _{0}\right) <0. 
\]%
Hence, (\ref{6.39}) implies that for sufficiently small values of $\beta >0$ 
\begin{equation}
J_{h,\kappa _{3},\alpha }\left( X_{\min }+\beta \left( Y-X_{\min }\right)
\right) \left( \theta _{0}\right) <J_{h,\kappa _{3},\alpha }\left( X_{\min
}\right) \left( \theta _{0}\right) .  \label{6.40}
\end{equation}%
However, (\ref{6.40}) contradicts to (\ref{6.019}).$\ \ \square $

\subsection{Accuracy estimates}

\label{sec:5.2}

\textbf{Theorem 5.2.}\emph{\ Assume that conditions of Theorems 3.3 and 5.1
hold. In particular, let }$\kappa _{3}$ \emph{and }$\theta _{0}$\emph{\ be
the pair of numbers of \ item 3 of Theorem 5.1 and let }

$\left( q_{\min ,\kappa _{3},\alpha },p_{\min ,\kappa _{3},\alpha }\right)
\left( \mathbf{x},\theta _{0}\right) \in \overline{G_{\theta _{0}}\left(
A\right) }$\emph{\ be the corresponding minimizer of the functional }$%
J_{\kappa _{3},\alpha }\left( q,p\right) \left( \theta \right) $\emph{\ on
the set }$\overline{G_{\theta _{0}}\left( A\right) }$\emph{. Let }%
\begin{equation}
Z_{\min }=\left( q_{\min ,\kappa _{3},\alpha }^{h},p_{\min ,\kappa
_{3},\alpha }^{h},q_{\min ,\kappa _{3},\alpha },p_{\min ,\kappa _{3},\alpha
}\right) \in \overline{G_{\theta _{0}}^{h}\left( A\right) }  \label{6.041}
\end{equation}%
\emph{be the corresponding point of the space }$H_{1,2,2,\theta _{0},fd}.$ 
\emph{Also, let }$X_{\min }\in \overline{G_{\theta _{0}}^{h}\left( A\right) }
$\emph{\ be the minimizer (\ref{6.41}) of the functional\ }$J_{h,\kappa
_{3},\alpha }\left( q^{h},p^{h},q,p\right) \left( \theta _{0}\right) $\emph{%
\ on the set }$\overline{G_{\theta _{0}}^{h}\left( A\right) },$ \emph{which
was mentioned in item 3 of Theorem 5.1. Denote}%
\begin{equation}
Z^{\ast }=\left( q^{\ast h},p^{\ast h},q^{\ast },p^{\ast }\right) \in
G_{\theta _{0}}^{h}\left( A\right) ,  \label{6.410}
\end{equation}%
\emph{see item 2 of Theorem 3.3 and, in particular, (\ref{5.121}) about the
pair }$(q^{\ast },p^{\ast })\left( \mathbf{x},\theta _{0}\right) \in
G_{\theta _{0}}\left( A\right) .$ \emph{Then the following accuracy
estimates hold}%
\begin{equation}
\left\Vert Z_{\min }-X_{\min }\right\Vert _{H_{1,2,2,\theta
_{0},fd}^{2}}\leq C_{1}\left( \sqrt{\alpha }+\sqrt{h}\right) \exp \left(
\kappa _{3}\left( a+c\right) ^{2}\right) ,  \label{6.42}
\end{equation}%
\begin{equation}
\left\Vert Z^{\ast }-X_{\min }\right\Vert _{H_{1,2,2,\theta
_{0},fd}^{2}}\leq C_{1}\left( \sqrt{\alpha }+\sqrt{h}\right) \exp \left(
\kappa _{3}\left( a+c\right) ^{2}\right) .  \label{6.420}
\end{equation}

\textbf{Remark 5.2}.  Rephrasing one of statements of subsection
1.3, we say now that if $\sqrt{\alpha }+\sqrt{h}$ is sufficiently
small, then estimates (\ref{6.42}) and (\ref{6.420}) imply that even if the
functional $J_{h,\kappa _{3},\alpha }\left( q^{h},p^{h},q,p\right) \left(
\theta _{0}\right) $ has several minimizers on the set $\overline{%
G_{\theta _{0}}^{h}\left( A\right) }$ , still all these minimizers
are close to each other, as well as to the points $Z_{\min }$ and $%
Z$$^{\ast }.$ The latter is good for computations.

\textbf{Proof of Theorem 5.2.} Using (\ref{6.1460}), (\ref{6.146}), (\ref%
{6.41}), (\ref{6.18}) and (\ref{6.041}), we obtain%
\begin{equation}
\left. 
\begin{array}{c}
J_{h,\kappa _{3},\alpha }\left( Z_{\min }\right) \left( \theta _{0}\right)
-J_{h,\kappa _{3},\alpha }\left( X_{\min }\right) \left( \theta _{0}\right) -
\\ 
-J_{h,\kappa _{3},\alpha }^{\prime }\left( X_{\min }\right) \left( Z_{\min
}-X_{\min }\right) \left( \theta _{0}\right) \geq \\ 
\geq C_{1}\exp \left( 2\kappa _{3}\left( a-c\right) ^{2}\right) \left\Vert
Z_{\min }-X_{\min }\right\Vert _{H_{1,2,2,\theta _{0},fd}^{2}}^{2}- \\ 
-C_{1}h\exp \left( 2\kappa _{3}\left( a+c\right) ^{2}\right) .%
\end{array}%
\right.  \label{6.43}
\end{equation}%
Since by (\ref{6.19}) 
\[
-J_{h,\kappa ,\alpha }\left( X_{\min }\right) \left( \theta _{0}\right)
-J_{h,\kappa ,\alpha }^{\prime }\left( X_{\min }\right) \left( \theta
_{0}\right) \left( Z_{\min }-X_{\min }\right) \leq 0, 
\]%
then (\ref{6.43}) implies%
\begin{equation}
\left. 
\begin{array}{c}
J_{h,\kappa ,\alpha }\left( Z_{\min }\right) \left( \theta _{0}\right) \geq
\\ 
\geq C_{1}\exp \left( 2\kappa _{3}\left( a-c\right) ^{2}\right) \left\Vert
Z_{\min }-X_{\min }\right\Vert _{H_{1,2,2,\theta _{0},fd}^{2}}^{2}- \\ 
-C_{1}h\exp \left( 2\kappa _{3}\left( a+c\right) ^{2}\right) .%
\end{array}%
\right.  \label{6.44}
\end{equation}%
By (\ref{6.35}) and (\ref{6.041}) 
\[
\left. J_{h,\kappa _{3},\alpha }\left( Z_{\min }\right) \left( \theta
_{0}\right) =J_{\kappa _{3},\alpha }\left( q_{\min ,\kappa _{3},\alpha
},p_{\min ,\kappa _{3},\alpha }\right) \left( \theta _{0}\right) +O\left(
h\right) \exp \left( 2\kappa _{3}\left( a+c\right) ^{2}\right) .\right. 
\]%
Substituting this in (\ref{6.44}), we obtain%
\begin{equation}
\left. 
\begin{array}{c}
J_{\kappa _{3},\alpha }\left( q_{\min ,\kappa _{3},\alpha },p_{\min ,\kappa
_{3},\alpha }\right) \left( \theta _{0}\right) \geq \\ 
\geq C_{1}\exp \left( 2\kappa _{3}\left( a-c\right) ^{2}\right) \left\Vert
Z_{\min }-X_{\min }\right\Vert _{H_{1,2,2,\theta _{0},fd}^{2}}^{2}- \\ 
-C_{1}h\exp \left( 2\kappa _{3}\left( a+c\right) ^{2}\right) .%
\end{array}%
\right.  \label{6.45}
\end{equation}

It follows from (\ref{3.19}), (\ref{3.20}), (\ref{3.240}) and (\ref{4.4})
that 
\begin{equation}
\left. 
\begin{array}{c}
J_{\kappa ,\alpha }\left( q+\widehat{q},p+\widehat{p}\right) \left( \theta
\right) \leq C_{1}J_{\kappa ,\alpha }\left( q,p\right) \left( \theta \right)
+C_{1}J_{\kappa ,\alpha }\left( \widehat{q},\widehat{p}\right) \left( \theta
\right) , \\ 
\forall \left( q,p\right) ,\left( \widehat{q},\widehat{p}\right) \in 
\overline{G_{\theta }\left( A\right) },\text{ }\forall \kappa >0,\text{ }%
\forall \theta \in \left[ 0,2\pi \right] .%
\end{array}%
\right.  \label{6.450}
\end{equation}%
Hence, in (\ref{6.45})%
\begin{equation}
\left. 
\begin{array}{c}
J_{\kappa _{3},\alpha }\left( q_{\min ,\kappa _{3},\alpha },p_{\min ,\kappa
_{3},\alpha }\right) \left( \theta _{0}\right) = \\ 
=J_{\kappa _{3},\alpha }\left( \left( q_{\min ,\kappa _{3},\alpha },p_{\min
,\kappa _{3},\alpha }\right) -\left( q^{\ast },p^{\ast }\right) +\left(
q^{\ast },p^{\ast }\right) \right) \left( \theta _{0}\right) \leq \\ 
\leq C_{1}J_{\kappa _{3},\alpha }\left( \left( q_{\min ,\kappa _{3},\alpha
}-q^{\ast },p_{\min ,\kappa _{3},\alpha }-p^{\ast }\right) \right) \left(
\theta _{0}\right) + \\ 
+C_{1}J_{\kappa _{3},\alpha }\left( q^{\ast },p^{\ast }\right) \left( \theta
_{0}\right) .%
\end{array}%
\right.  \label{6.451}
\end{equation}%
Hence, (\ref{5.18}) and (\ref{6.451}) imply%
\begin{equation}
\left. 
\begin{array}{c}
J_{\kappa _{3},\alpha }\left( q_{\min ,\kappa _{3},\alpha },p_{\min ,\kappa
_{3},\alpha }\right) \left( \theta _{0}\right) \leq \\ 
\leq C_{1}J_{\kappa _{3},\alpha }\left( \left( q_{\min ,\kappa _{3},\alpha
}-q^{\ast },p_{\min ,\kappa _{3},\alpha }-p^{\ast }\right) \right) \left(
\theta _{0}\right) +C_{1}\alpha .%
\end{array}%
\right.  \label{6.452}
\end{equation}%
Next, using again (\ref{3.19}), (\ref{3.20}), (\ref{3.240}) and (\ref{4.4}),
we obtain%
\begin{equation}
\left. J_{\kappa _{3},\alpha }\left( \left( q_{\min ,\kappa _{3},\alpha
}-q^{\ast },p_{\min ,\kappa _{3},\alpha }-p^{\ast }\right) \right) \left(
\theta _{0}\right) \leq \right.  \label{6.453}
\end{equation}%
\[
\leq C_{1}\exp \left( 2\kappa _{3}\left( a+c\right) ^{2}\right) \left\Vert
\left( q_{\min ,\kappa _{3},\alpha }-q^{\ast },p_{\min ,\kappa _{3},\alpha
}-p^{\ast }\right) \right\Vert _{H_{1,2,\theta _{0}}^{2}}^{2}+C_{1}\alpha . 
\]%
Hence, by (\ref{5.13}) and (\ref{6.451})-(\ref{6.453})%
\begin{equation}
J_{\kappa _{3},\alpha }\left( q_{\min ,\kappa _{3},\alpha },p_{\min ,\kappa
_{3},\alpha }\right) \left( \theta _{0}\right) \leq C_{1}\alpha \exp \left(
2\kappa _{3}\left( a+c\right) ^{2}\right) .  \label{6.454}
\end{equation}%
Comparing (\ref{6.454}) with (\ref{6.45}), we obtain the first target
estimate (\ref{6.42}) of this theorem.

We now prove (\ref{6.420}). By triangle inequality and (\ref{6.42})%
\[
\left\Vert Z^{\ast }-X_{\min }\right\Vert _{H_{1,2,2,\theta
_{0},fd}^{2}}\leq \left\Vert Z^{\ast }-Z_{\min }\right\Vert
_{H_{1,2,2,\theta _{0},fd}^{2}}+\left\Vert Z_{\min }-X_{\min }\right\Vert
_{H_{1,2,2,\theta _{0},fd}^{2}}\leq 
\]%
\begin{equation}
\leq \left\Vert Z^{\ast }-Z_{\min }\right\Vert _{H_{1,2,2,\theta
_{0},fd}^{2}}+C_{1}\left( \sqrt{\alpha }+\sqrt{h}\right) \exp \left( \kappa
_{3}\left( a+c\right) ^{2}\right) .  \label{6.455}
\end{equation}%
Next, by (\ref{6.20}), (\ref{6.21}), (\ref{6.041}) and (\ref{6.410}) 
\[
\left( q^{\ast h}-q_{\min ,\kappa _{3},\alpha }^{h},p^{\ast h}-p_{\min
,\kappa _{3},\alpha }^{h},q^{\ast }-q_{\min ,\kappa _{3},\alpha },p^{\ast
}-p_{\min ,\kappa _{3},\alpha }\right) \in \overline{G_{\theta
_{0}}^{h}\left( 2A\right) }. 
\]%
Hence, using (\ref{6.146}), (\ref{6.0146}) and (\ref{6.041}), we obtain%
\begin{equation}
\left\Vert Z^{\ast }-Z_{\min }\right\Vert _{H_{1,2,2,\theta
_{0},,fd}^{2}}\leq \left\Vert \left( q^{\ast }-q_{\min ,\kappa _{3},\alpha
},p^{\ast }-p_{\min ,\kappa _{3},\alpha }\right) \right\Vert _{H_{1,2,\theta
_{0}}^{2}}+C_{1}\sqrt{h}.  \label{6.456}
\end{equation}%
Combining (\ref{6.456}) with (\ref{5.13}), we obtain%
\begin{equation}
\left\Vert Z^{\ast }-Z_{\min }\right\Vert _{H_{1,2,2,\theta
_{0},fd}^{2}}\leq C_{1}\sqrt{\alpha }\exp \left[ -\kappa \left( a-c\right)
^{2}\right] +C_{1}\sqrt{h}.  \label{6.457}
\end{equation}%
Combining (\ref{6.455}) with (\ref{6.457}), we obtain (\ref{6.420}), which
is the second target estimate of this theorem. $\square $

\subsection{Minimizing sequence}

\label{sec:5.3}

It is assumed in Theorem 5.2 that there exists such a pair of numbers $%
\kappa _{3}$ and\emph{\ }$\theta _{0}$ for which there exists a minimizer $%
Z_{\min }$ in (\ref{6.041}) of the functional $J_{h,\kappa ,\alpha }\left(
q^{h},p^{h},q,p\right) \left( \theta _{0}\right) $ on the set $\overline{%
G_{\theta _{0}}^{h}\left( A\right) }.$ We present in this section an analog
of Theorem 5.2 for the case when the assumption about the existence of a
minimizer is not imposed.

Let $\kappa _{2}=\kappa _{2}\left( A,\Omega ,\varepsilon \right) \geq 1$ be
the number of Theorem 5.1. In particular, unlike Theorem 5.2, where the
numbers $\kappa _{3}\geq \kappa _{2}$ and $\theta _{0}\in \left[ 0,2\pi %
\right] $ are fixed, we assume here that $\kappa \geq \kappa _{2}$ and $%
\theta \in \left[ 0,2\pi \right] $ are two arbitrary numbers. Denote 
\begin{equation}
m\left( h,\kappa ,\alpha ,\varepsilon ,A,\theta \right) =\inf_{\overline{%
G_{\theta }^{h}\left( A\right) }}J_{h,\kappa ,\alpha }\left(
q^{h},p^{h},q,p\right) ,\text{ }\forall \kappa \geq \kappa _{2}\text{, }%
\forall \theta \in \left[ 0,2\pi \right] .  \label{6.458}
\end{equation}%
It follows from (\ref{6.458}) that there exists a minimizing sequence $%
\left\{ \left( q_{n}^{h},p_{n}^{h},q_{n},p_{n}\right) \right\}
_{n=1}^{\infty }\subset \overline{G_{\theta }^{h}\left( A\right) }$ of the
functional $J_{h,\kappa ,\alpha }\left( q^{h},p^{h},q,p\right) $ such that 
\begin{equation}
\left. 
\begin{array}{c}
\lim_{n\rightarrow \infty }J_{h,\kappa ,\alpha }\left(
q_{n}^{h},p_{n}^{h},q_{n},p_{n}\right) =m\left( h,\kappa ,\alpha
,\varepsilon ,A,\theta \right) , \\ 
J_{h,\kappa ,\alpha }\left( q_{n}^{h},p_{n}^{h},q_{n},p_{n}\right) \geq
m\left( h,\kappa ,\alpha ,\varepsilon ,A,\theta \right) ,\text{ }\forall
n\geq 1.%
\end{array}%
\right.   \label{6.459}
\end{equation}%
Denote%
\begin{equation}
Z_{n}=\left( q_{n}^{h},p_{n}^{h},q_{n},p_{n}\right) \in \overline{G_{\theta
}^{h}\left( A\right) }.  \label{6.460}
\end{equation}

\textbf{Theorem 5.3.} \emph{Assume that (\ref{6.410}) is replaced with}%
\begin{equation}
Z^{\ast }=\left( q^{\ast h},p^{\ast h},q^{\ast },p^{\ast }\right) \in
G_{\theta }^{h}\left( A\right) ,\text{ }\forall \theta \in \left[ 0,2\pi %
\right] .  \label{1}
\end{equation}%
\emph{Let }$\kappa \geq \kappa _{2}$\emph{\ be an arbitrary number. Let }$%
\left( q_{\min ,\kappa ,\alpha },p_{\min ,\kappa ,\alpha }\right) \left( 
\mathbf{x},\theta \right) \in \overline{G_{\theta }\left( A\right) }$\emph{\
be the minimizer of the functional }$J_{\kappa ,\alpha }\left( q,p\right)
\left( \theta \right) $\emph{\ on the set }$\overline{G_{\theta }\left(
A\right) },$ \emph{which was found in Theorem 3.2. Let }$\left\{ \left(
q_{n}^{h},p_{n}^{h},q_{n},p_{n}\right) \right\} _{n=1}^{\infty }\subset 
\overline{G_{\theta }^{h}\left( A\right) }$ \emph{be the above minimizing
sequence of the functional }$J_{h,\kappa ,\alpha }\left(
q^{h},p^{h},q,p\right) $\emph{\ satisfying (\ref{6.459}). Then for each }$%
\theta \in \left[ 0,2\pi \right] $\emph{\ there exists such an integer }$%
N=N\left( h,\kappa ,\alpha ,\varepsilon ,A,\theta \right) >1$\emph{\
depending only on listed parameters that the following analogs of accuracy
estimates (\ref{6.42}) and (\ref{6.420}) hold }%
\begin{equation}
\left. 
\begin{array}{c}
\left\Vert q_{\min ,\kappa ,\alpha }-q_{n},p_{\min ,\kappa ,\alpha
}-p_{n}\right\Vert _{H_{1,2,\theta }^{2}}\leq \\ 
\leq C_{1}\left( \sqrt{\alpha }+\sqrt{h}\right) \exp \left( \kappa \left(
a+c\right) ^{2}\right) ,\text{ }\forall n\geq N,%
\end{array}%
\right.  \label{6.461}
\end{equation}%
\begin{equation}
\left. 
\begin{array}{c}
\left\Vert q^{\ast }-q_{n},p^{\ast }-p_{n}\right\Vert _{H_{1,2,\theta
}^{2}}\leq \\ 
\leq C_{1}\left( \sqrt{\alpha }+\sqrt{h}\right) \exp \left( \kappa \left(
a+c\right) ^{2}\right) ,\text{ }\forall n\geq N,\text{ }%
\end{array}%
\right.  \label{6.462}
\end{equation}%
\begin{equation}
\left. \left\Vert Z^{\ast }-Z_{n}\right\Vert _{H_{1,2,2,\theta ,fd}^{2}}\leq
C_{1}\left( \sqrt{\alpha }+\sqrt{h}\right) \exp \left( \kappa \left(
a+c\right) ^{2}\right) ,\text{ }\forall n\geq N.\right.  \label{6.4620}
\end{equation}

\textbf{Proof.} Using (\ref{6.18}), we obtain%
\[
J_{h,\kappa ,\alpha }\left( q_{n}^{h},p_{n}^{h},q_{n},p_{n}\right) \left(
\theta \right) -J_{h,\kappa ,\alpha }\left( q_{\min ,\kappa ,\alpha
}^{h},p_{\min ,\kappa ,\alpha }^{h},q_{\min ,\kappa ,\alpha },p_{\min
,\kappa ,\alpha }\right) \left( \theta \right) - 
\]%
\[
\left. 
\begin{array}{c}
-J_{h,\kappa ,\alpha }^{\prime }\left( q_{\min ,\kappa ,\alpha }^{h},p_{\min
,\kappa ,\alpha }^{h},q_{\min ,\kappa ,\alpha },p_{\min ,\kappa ,\alpha
}\right) \left( \theta \right) \\ 
\left( q_{n}^{h}-q_{\min ,\kappa ,\alpha }^{h},p_{n}^{h}-p_{\min ,\kappa
,\alpha }^{h},q_{n}-q_{\min ,\kappa ,\alpha },p_{n}-p_{\min ,\kappa ,\alpha
}\right) \geq%
\end{array}%
\right. 
\]%
\begin{equation}
\geq C_{1}\exp \left( 2\kappa \left( a-c\right) ^{2}\right) \left\Vert
\left( q_{n}^{h}-q_{\min ,\kappa ,\alpha }^{h},p_{n}^{h}-p_{\min ,\kappa
,\alpha }^{h}\right) \right\Vert _{H_{1,2,\theta ,fd}}^{2}+  \label{6.463}
\end{equation}%
\[
+\alpha \left\Vert \left( q_{n}-q_{\min ,\kappa ,\alpha },p_{n}-p_{\min
,\kappa ,\alpha }\right) \right\Vert _{H_{1,2,\theta }}^{2}-C_{1}h\exp
\left( 2\kappa \left( a+c\right) ^{2}\right) ,\text{ }\forall n\geq 1. 
\]%
Next, by (\ref{6.0146})%
\begin{equation}
\left. 
\begin{array}{c}
\left\Vert \left( q_{n}^{h}-q_{\min ,\kappa ,\alpha }^{h},p_{n}^{h}-p_{\min
,\kappa ,\alpha }^{h}\right) \right\Vert _{H_{1,2,\theta ,fd}}^{2}\geq \\ 
\geq \left\Vert \left( q_{n}-q_{\min ,\kappa ,\alpha },p_{n}-p_{\min ,\kappa
,\alpha }\right) \right\Vert _{H_{1,2,\theta }^{2}}^{2}-C_{1}h,\text{ }%
\forall n\geq 1.%
\end{array}%
\right.  \label{6.464}
\end{equation}%
Also, by (\ref{6.170})%
\[
\left. 
\begin{array}{c}
-J_{h,\kappa ,\alpha }^{\prime }\left( q_{\min ,\kappa ,\alpha }^{h},p_{\min
,\kappa ,\alpha }^{h},q_{\min ,\kappa ,\alpha },p_{\min ,\kappa ,\alpha
}\right) \left( \theta \right) \\ 
\left( q_{n}^{h}-q_{\min ,\kappa ,\alpha }^{h},p_{n}^{h}-p_{\min ,\kappa
,\alpha }^{h},q_{n}-q_{\min ,\kappa ,\alpha },p_{n}-p_{\min ,\kappa ,\alpha
}\right) \geq%
\end{array}%
\right. 
\]%
\begin{equation}
\geq -J_{\kappa ,\alpha }^{\prime }\left( q_{\min ,\kappa ,\alpha },p_{\min
,\kappa ,\alpha }\right) \left( \theta \right) \left( q_{n}-q_{\min ,\kappa
,\alpha },p_{n}-p_{\min ,\kappa ,\alpha }\right) -  \label{6.465}
\end{equation}%
\[
-C_{1}h\exp \left( 2\kappa \left( a+c\right) ^{2}\right) ,\text{ }\forall
n\geq 1. 
\]%
Using (\ref{6.463})-(\ref{6.465}), we obtain%
\begin{equation}
\left. 
\begin{array}{c}
J_{h,\kappa ,\alpha }\left( q_{n}^{h},p_{n}^{h},q_{n},p_{n}\right) \left(
\theta \right) - \\ 
-J_{h,\kappa ,\alpha }\left( q_{\min ,\kappa ,\alpha }^{h},p_{\min ,\kappa
,\alpha }^{h},q_{\min ,\kappa ,\alpha },p_{\min ,\kappa ,\alpha }\right)
\left( \theta \right) - \\ 
-J_{\kappa ,\alpha }^{\prime }\left( q_{\min ,\kappa ,\alpha },p_{\min
,\kappa ,\alpha }\right) \left( \theta \right) \left( q_{n}-q_{\min ,\kappa
,\alpha },p_{n}-p_{\min ,\kappa ,\alpha }\right) \geq \\ 
\geq C_{1}\exp \left( 2\kappa \left( a-c\right) ^{2}\right) \left\Vert
\left( q_{n}-q_{\min ,\kappa ,\alpha },p_{n}-p_{\min ,\kappa ,\alpha
}\right) \right\Vert _{H_{1,2,\theta }^{2}}^{2}- \\ 
-C_{1}h\exp \left( 2\kappa \left( a+c\right) ^{2}\right) ,\text{ }\forall
n\geq 1.%
\end{array}%
\right.  \label{6.466}
\end{equation}%
Next, by (\ref{5.4}) $-J_{\kappa ,\alpha }^{\prime }\left( q_{\min ,\kappa
,\alpha },p_{\min ,\kappa ,\alpha }\right) \left( \theta \right) \left(
q_{n}-q_{\min ,\kappa ,\alpha },p_{n}-p_{\min ,\kappa ,\alpha }\right) \leq
0.$ Since $-J_{h,\kappa ,\alpha }\left( q_{\min ,\kappa ,\alpha
}^{h},p_{\min ,\kappa ,\alpha }^{h},q_{\min ,\kappa ,\alpha },p_{\min
,\kappa ,\alpha }\right) \left( \theta \right) \leq 0$ as well, then (\ref%
{6.460}) and (\ref{6.466}) imply%
\begin{equation}
\left. 
\begin{array}{c}
J_{h,\kappa ,\alpha }\left( Z_{n}\right) \left( \theta \right) \geq \\ 
\geq C_{1}\exp \left( 2\kappa \left( a-c\right) ^{2}\right) \left\Vert
\left( q_{n}-q_{\min ,\kappa ,\alpha },p_{n}-p_{\min ,\kappa ,\alpha
}\right) \right\Vert _{H_{1,2,\theta }^{2}}^{2}- \\ 
-C_{1}h\exp \left( 2\kappa \left( a+c\right) ^{2}\right) ,\text{ }\forall
n\geq 1.%
\end{array}%
\right.  \label{6.467}
\end{equation}

Until now we have not used sufficiently large numbers $n$. We now start to
use them. By (\ref{6.458}), (\ref{6.459}) and (\ref{1}) there exists a
sufficiently large integer $N=N\left( h,\kappa ,\alpha ,\varepsilon
,A,\theta \right) >1$ depending only on listed parameters such that 
\[
J_{h,\kappa ,\alpha }\left( Z_{n}\right) \left( \theta \right) \leq
J_{h,\kappa ,\alpha }\left( Z^{\ast }\right) \left( \theta \right) ,\text{ }%
\forall n\geq N,\text{ }\forall \theta \in \left[ 0,2\pi \right] . 
\]%
Hence, (\ref{6.467}) implies%
\begin{equation}
\left. 
\begin{array}{c}
J_{h,\kappa ,\alpha }\left( Z^{\ast }\right) \left( \theta \right) \geq \\ 
\geq C_{1}\exp \left( 2\kappa \left( a-c\right) ^{2}\right) \left\Vert
\left( q_{n}-q_{\min ,\kappa ,\alpha },p_{n}-p_{\min ,\kappa ,\alpha
}\right) \right\Vert _{H_{1,2,\theta }^{2}}^{2}- \\ 
-C_{1}h\exp \left( 2\kappa \left( a+c\right) ^{2}\right) ,\text{ }\forall
n\geq N.%
\end{array}%
\right.  \label{6.468}
\end{equation}

Using (\ref{5.18}), Proposition 4.1, (\ref{6.9}), (\ref{6.90}), (\ref{6.0146}%
) and (\ref{1}), we obtain 
\[
\left. 
\begin{array}{c}
J_{h,\kappa ,\alpha }\left( Z^{\ast }\right) \left( \theta \right)
=J_{h,\kappa ,\alpha }\left( q^{\ast h},p^{\ast h},q^{\ast },p^{\ast
}\right) \left( \theta \right) = \\ 
=\sqrt{\varepsilon }\int\limits_{\Omega ,rt}\left[ \left( F_{1}^{h}\left(
q^{\ast h},p^{\ast h}\right) \left( \mathbf{x}_{i,j}^{h},\theta \right)
\right) ^{2}+\left( F_{2}^{h}\left( q^{\ast h},p^{\ast h}\right) \left( 
\mathbf{x}_{i,j}^{h},\theta \right) \right) ^{2}\right] W_{\kappa
}^{h}\left( \mathbf{x}_{i,j}^{h}\right) d\mathbf{x}_{i,j}^{h}+ \\ 
+\alpha \left\Vert \left( q^{\ast },p^{\ast }\right) \right\Vert
_{H_{1,2,\theta }}^{2}\leq \\ 
\leq J_{\kappa ,\alpha }\left( q^{\ast },p^{\ast }\right) +C_{1}h\exp \left(
2\kappa \left( a+c\right) ^{2}\right) \leq C_{1}\alpha +C_{1}h\exp \left(
2\kappa \left( a+c\right) ^{2}\right) .%
\end{array}%
\right. 
\]%
Combining this with (\ref{6.468}), we obtain (\ref{6.461}), which is the
first target estimate of this theorem. The second target estimate (\ref%
{6.462}) follows immediately from (\ref{5.13}), (\ref{6.461}) and triangle
inequality. The third target estimate (\ref{6.4620}) follows immediately
from (\ref{6.200})-(\ref{6.0146}) and (\ref{6.462}). $\square $

\section{Numerical Minimization of the Functional $J_{h,\protect\kappa ,%
\protect\alpha }\left( q^{h},p^{h},q,p\right) \left( \protect\theta \right) $
In (\protect\ref{6.16})}

\label{sec:6}

Our numerical setup for data generation via the solution of the forward
problem and the semi-discrete convexification method is identical to the one
used in \cite{EITIP2025}. We now briefly present it for the sake of
completeness. We refer to \cite{EITIP2025} for some details.

We again refer to Figure 1. For each $\theta \in \left[ 0,2\pi \right) $ and
for the corresponding position $\mathbf{x}_{0}=\mathbf{x}_{0}(\theta )\in
E_{B}\left( a,b\right) $ in (\ref{1.00}) the forward problem (\ref{2-2}) is
solved by the Finite Element Method (FEM) in the disk $P_{D}$ defined in (%
\ref{1.00}). The radius of this disk is $D=3$ and its center is at the point 
$\left( a,b\right) =(1.5,1.5)$. The circle $E_{B}\left( a,b\right) $ is now $%
E_{B}\left( 1.5,1.5\right) $ with its radius $B=2$. The spatial mesh size of
the FEM is 1/160. The square $\Omega $ in (\ref{1.0}) is now 
\begin{equation}
\Omega =\left\{ \mathbf{x}=\left( x,y\right) :1<x,y<2\right\} ,  \label{7.3}
\end{equation}%
i.e., $c=0.5$ in (\ref{1.0}). The source function $g(\mathbf{x}-\mathbf{x}%
_{0})$ is as in (\ref{2.02}) with $\xi =0.1.$ We use 199 discrete source
positions corresponding to $\theta _{n}=n\rho _{\theta }$, $n=1,2,\cdots
,199 $, $\rho _{\theta }=\pi /100$.

For the convexification method we have used a finite difference
scheme to discretize the Partial Differential Operators in (\ref{4.4}) in
the computational domain $\Omega $ in (\ref{7.3}). We have used a
two coarse grids: with $10\times 10$ of $20\times 20$
pixels, which means grid step sizes $h_{1}=0.1$ and $h_{2}=0.05$ 
respectively in (\ref{6.00}). On the other hand, the grid step
size for the convexification method on the fine grid in \cite[page 17]%
{EITIP2025} was $h_{3}=1/40=0.025$, i.e. it was two times less than 
$h_{2}$. This explains the $14.25=57/4$\ fold speed up in
computations on the grid with $h=h_{2}=0.05$, see (\ref{1.1}). We
have used parameters $\varepsilon ,\kappa $, $\alpha $\ in
(\ref{4.4}) the same as the ones in formula (8.6) of \cite{EITIP2025}: 

\begin{equation}
\alpha =0.01,\text{ }\varepsilon =0.0002,\kappa =3.  \label{8.01}
\end{equation}%
Since the minimization in the convexification method is carried
out on the set $\overline{G_{\theta }\left( A\right) }$ with
restrictions (see (\ref{4.3}) and Definition in subsection 4.2), then we
have used the MATLAB built-in function fminunc. The true solutions and the
output of the deep learning network are defined on a finer grid of $%
128\times 128$ pixels.

In the functional $J_{h,\kappa ,\alpha }\left(
q^{h},p^{h},q,p\right) \left( \theta \right) ,$ we write the
differential operators in the form of finite differences and minimize the
resulting discretized functionals with respect to the values of
corresponding functions at grid points. The norm $H^{6}\left( \Omega \right) 
$ is inconvenient for the numerical implementation.
Hence, we replace in (\ref{4.4}), (\ref{5.16}), (\ref{6.150}) and (\ref{6.16}%
) the regularization term $\alpha \left\Vert \left( q,p\right) \left( 
\mathbf{x},\theta \right) \right\Vert _{H_{1,2,\theta }}^{2}$ with 
$\alpha \left\Vert \left( q,p\right) \left( \mathbf{x},\theta \right)
\right\Vert _{H^{2}\left( \Omega \right) \times H^{2}\left( \Omega \right)
}^{2}.$ This replacement is similar with the one in \cite[page 17]%
{EITIP2025}. Since we actually work with not too many grid points and since
norms in finite dimensional spaces are equivalent, then this replacement has
not provided a negative impact on our numerical results. An extension of the
above theory on the fully discrete case is outside of the scope of this
publication. 

As the starting point for the iterations of fmincon we took $q=p=0.$%
 The iterations of fmincon stop when the condition 
\begin{equation}
\left\vert \nabla J_{h,\kappa ,\alpha }\left( q^{h},p^{h},q,p\right) \left(
\theta \right) \right\vert \leq 10^{-2}  \label{8.04}
\end{equation}%
is met. In fact, the reason of stopping criterion (\ref{8.04}) is
clear from Figure 2. We display in Figure 2. This figure explains the
stopping criterion (\ref{8.04}).

\begin{figure}[tbph]
\centering
\includegraphics[width=5in]{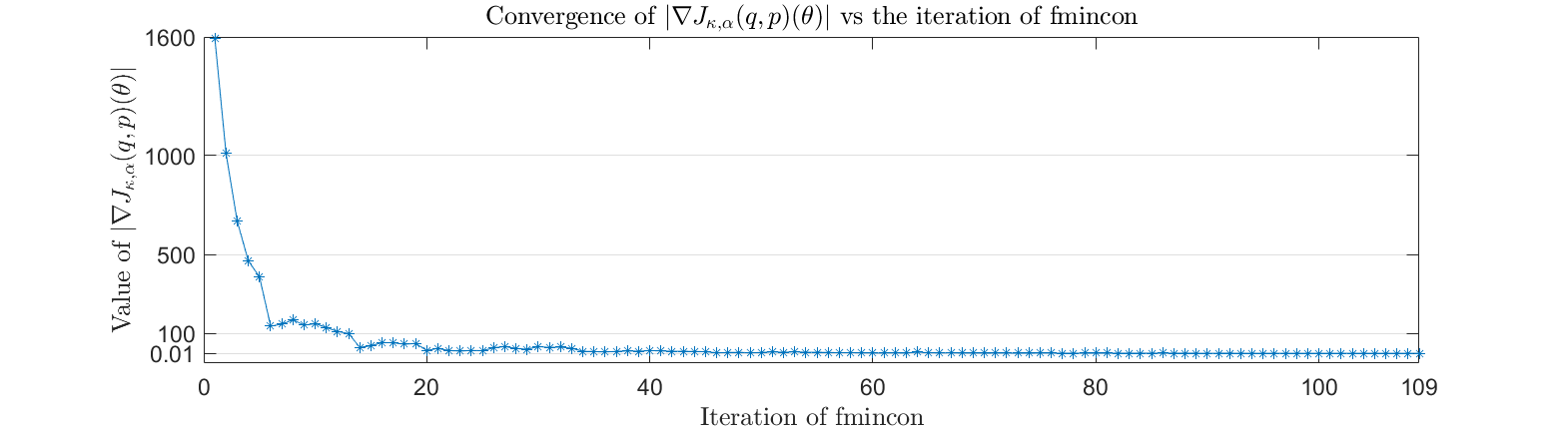}
\caption{Convergence behavior of the gradient $\left\vert \protect\nabla %
J_{h,\protect\kappa,\protect\alpha} \left(q^{h},p^{h},q,p\right)(\protect%
\theta)\right\vert$ versus the iterations of \textbf{fmincon}, for $%
h=h_{2}=0.05$, $\protect\alpha=0.01$, $\protect\varepsilon=0.0002$, $\protect%
\kappa=3$. This is for one of the test cases of Figure 8.  This figure explains our
stopping criterion: iterations of all cases of the convexification method on
the coarse grid are stopped as soon as the inequality $\left\vert \protect%
\nabla J_{h,\protect\kappa,\protect\alpha} \left(q^{h},p^{h},q,p\right)(%
\protect\theta)\right\vert \leq 1$ is met.}
\label{plot_grad}
\end{figure}



\bigskip



\section{Deep Learning With A Priory Accuracy Estimates}

\label{sec:7}

\subsection{Estimating the accuracy of the starting point of deep learning}

\label{sec:7.1}

We now estimate the accuracy of the starting point of iterations for
deep learning. We start from the training step of our deep learning
procedure. To simplify the presentation, we assume that conditions of
Theorem 5.2 hold. The case of Theorem 5.3 is considered similarly.

Let $\kappa _{3}$ and $\theta _{0}$ be two numbers of Theorem 5.2. Since we
now have only one value of the number $\theta =\theta _{0}\in \left[ 0,2\pi %
\right] $ instead of the whole interval $\theta \in \left[ 0,2\pi \right] ,$
then, recalling (\ref{3.50}), (\ref{3.240}), (\ref{5.5}), (\ref{5.6}) and
using (\ref{6.041}), we replace (\ref{5.7}) with%
\begin{equation}
r_{\kappa _{3},\alpha }\left( \mathbf{x}_{i,j}^{h}\right) =-\left( \Delta
^{h}\psi _{\min ,\kappa _{3},\alpha }+\left( \nabla \psi _{\min ,\kappa
_{3},\alpha }\right) ^{2}\right) \left( \mathbf{x}_{i,j}^{h},\theta
_{0}\right) .  \label{7.1}
\end{equation}%
Here the function $r_{\kappa _{3},\alpha }\left( \mathbf{x}_{i,j}^{h}\right) 
$ is our computational result obtained by the convexification method on a
coarse grid. Next, we interpolate formula (\ref{7.1}) on a fine grid in $%
\Omega $ and apply the procedure described in (\ref{5.8})-(\ref{5.110}). We
obtain this way an approximation $\sigma _{\kappa _{3},\alpha }\left( 
\mathbf{x}\right) $ of the true coefficient $\sigma ^{\ast }\left( \mathbf{x}%
\right) .$ However, it is hard to provide an analytical estimate of the
accuracy of that interpolation of (\ref{7.1}) on a fine grid. Hence,
we use the following assumption:

\textbf{Assumption 7.1.} The accuracy of the latter approximation of
the function $\sigma _{\kappa _{3},\alpha }\left( \mathbf{x}\right) $%
 is the same as the one in (\ref{6.420}), i.e. 
\begin{equation}
\left\Vert \sigma _{\kappa _{3},\alpha }-\sigma ^{\ast }\right\Vert
_{L_{2}\left( \Omega \right) }\leq C_{1}\left( \sqrt{\alpha }+\sqrt{h}%
\right) \exp \left( \kappa _{3}\left( a+c\right) ^{2}\right) .  \label{7.2}
\end{equation}

\subsection{A priory accuracy estimate}

\label{sec:7.2}

As it was stated in subsection 1.1, we use the solution obtained by the
semi-discrete version of the convexification method on a coarse grid as
input image for the DNN. In other words, we use the function $\sigma
_{\kappa _{3},\alpha }\left( \mathbf{x}\right) $ for this purpose. 
Therefore, if accepting Assumptions 1.1 and 7.1, then, rephrasing one of
sentences of subsection 1.4, a priory accuracy estimate (\ref{7.2}) holds
for our entire two-stage process.

 We note that multipliers of the term $\left( \sqrt{\alpha }+\sqrt{h%
}\right) $ in estimates of Theorems 5.2, 5.3 as well as in (\ref%
{7.2}) have large values for specific values of parameters listed in section
6. This means that for chosen values $\alpha =0.01,$ $h=0.05$ (see (%
\ref{8.01}) and (\ref{7.5})) the number in the right hand side of (\ref{7.2}%
) is large. Still, our argument here is that we consider those estimates as
asymptotic ones, i.e. sort of \textquotedblleft ideal" ones, assuming that $%
\sqrt{\alpha }+\sqrt{h}\rightarrow 0^{+},$ see {(\ref{1.3}). In
fact, such things often occur in computations if asymptotic theories are
used. Consider, for example an asymptotic theory, which claims that a
certain formula $X$\ is a quite accurate one if a certain parameter 
$Y$ is sufficiently large. However, the next natural question is:
What actually terms \textquotedblleft quite accurate" and \textquotedblleft
sufficiently large" mean in a specific computational setting? For example:
is $Y=1$ sufficiently large for $X$ to be quite accurate
or maybe $Y=100$ is the lowest appropriate value, which one should
use for that setting? The first author has a rich computational experience
with respect to the values of the parameter $\kappa $ in various
Carleman Weight Functions used in computations of the convexification
method. This experience tells us that in a specific computational setting at
hands, which always uses specific values of a number of parameters, only
numerical results can provide those values of $Y,$ which ensure a
good accuracy of the formula $X.$ Indeed, even though the
convexification theory always requires $\kappa >>1,$ see, e.g.
above theorems, in the computational practice $\kappa \in \left[ 1,5\right] $%
 has been consistently sufficient so far, see, e.g. \cite{KL}-\cite%
{TTTP}. In particular, in computations of this paper $\kappa =3,$
see (\ref{8.01}).}

Another possible argument here is that in real applications of our
CIP for EIT to, e.g. medicine, the values of corresponding parameters
forming the above mentioned multipliers of estimates of Theorems 5.2, 5.3
and (\ref{7.2}) might well be such that those multipliers are not large.
Still, the main point here is that the right hand sides of those estimates
are small, as long as the term $\left( \sqrt{\alpha }+\sqrt{h}\right) $ is
sufficiently small.

As to the computational studies of this paper, we have conducted
numerical experiments for those values of parameters, which we have selected
in the previous publication \cite[page 17]{EITIP2025}. Recomputing for
different values of those parameters would be a major time consuming effort,
which is outside of the scope of this publication. 

\subsection{Outline of the deep learning procedure for our case}

\label{sec:7.3}

We propose a supervised hybrid image-to-image model that
combines a Denoising Convolutional Neural Network (DnCNN) with a Residual
U-Net enhanced by Squeeze-and-Excitation (SE) attention blocks. This
architecture is designed for grayscale reconstruction tasks, such as
denoising Chinese-character images while preserving sharp strokes.

The learning procedure is supervised because each conductivity image
reconstructed by the semi-discrete convexification method on the coarse
grid is paired with its corresponding true conductivity image, which
serves as the ground-truth target. The trainable parameters of the network
are optimized by minimizing the discrepancy between the predicted and
ground-truth images using the combined MS-SSIM and $L^1$ loss function
defined in \eqref{3}.

\begin{figure}[htbp]
\centering
\includegraphics[width=0.9\linewidth]{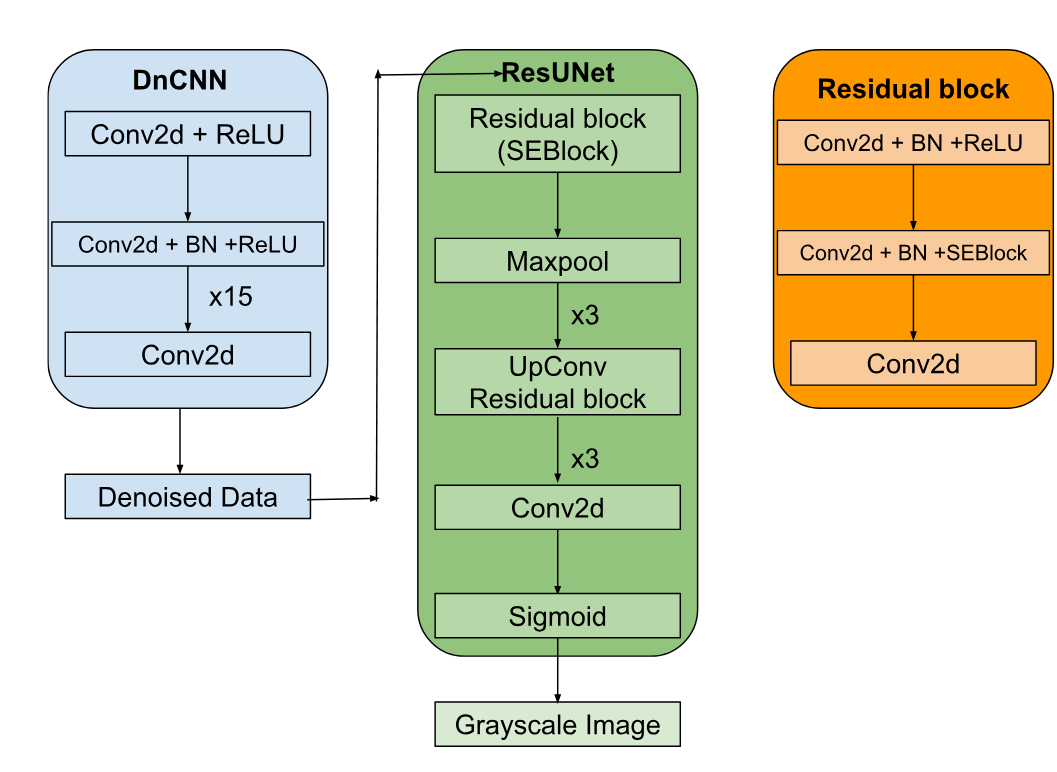}
\caption{Architecture of the proposed DnCNN + ResUNet model with SEBlock.
The grayscale input image is denoised by DnCNN and then passed to a
SE-enhanced ResUNet for reconstruction.}
\label{fig:cnn_layers}
\end{figure}

For this image-to-image translation task, we adopt the U-Net architecture
\cite{Unet15} (see Fig.~\ref{fig:cnn_layers}), which has demonstrated
remarkable success in biomedical image segmentation and other image
restoration tasks. The U-Net consists of a contracting path (encoder)
and an expansive path (decoder).

\clearpage

\subsection{The training procedure}

\label{sec:7.4} 

We first describe the training stage of our deep learning procedure. For each
training sample $i=1,\ldots,N$, we first apply the semi-discrete
convexification method on the coarse grid. The resulting reconstruction
$\sigma_{\mathrm{conv}}^{(i)}$ is then interpolated onto the fine grid in
$\Omega$ and used as the input image to the DNN. The corresponding true
conductivity image $\sigma^{(i)\ast}$ is used as the ground-truth target.

Thus, the training dataset consists of the paired images
\[
\bigl(\sigma_{\mathrm{conv}}^{(i)},\,\sigma^{(i)\ast}\bigr),
\qquad i=1,\ldots,N,
\]

where $\sigma_{\mathrm{conv}}^{(i)}$ is the reconstruction obtained by the
semi-discrete convexification method on the coarse grid and
$\sigma^{(i)\ast}$ is the corresponding true solution. For each training
sample, the convexification reconstruction is therefore not an initial guess
for the optimization of the neural-network parameters, but rather the input
image supplied to the network.

The network produces the predicted conductivity image
\[
\sigma_{\mathrm{pred}}^{(i)}
=
\text{U-Net}\bigl(\sigma_{\mathrm{conv}}^{(i)}\bigr),
\]
 
and its trainable parameters are optimized by minimizing the discrepancy
between $\sigma_{\mathrm{pred}}^{(i)}$ and the corresponding ground-truth
image $\sigma^{(i)\ast}$. We use the following loss function: 
\begin{equation}
\mathcal{L}_{\mathrm{total}}
=
\gamma \bigl(1-\mathrm{MS\text{-}SSIM}\bigr)
+
(1-\gamma)\,\mathrm{L1Loss},
\qquad
\gamma=0.84.
\label{3}
\end{equation}

Hence, during training, the procedure for each sample can be summarized as
\[
\sigma_{\mathrm{conv}}^{(i)}
\longrightarrow
\mathrm{DNN}
\longrightarrow
\sigma_{\mathrm{pred}}^{(i)},
\]
where $\sigma_{\mathrm{pred}}^{(i)}$ is compared with
$\sigma^{(i)\ast}$ through the loss function \eqref{3}. The minimization of
\eqref{3} is performed with respect to the trainable parameters of the neural
network.

For each training sample, the convexification reconstruction is obtained on
the coarse grid. As discussed above, this reconstruction satisfies the
corresponding a priori accuracy estimate of order 

\[
O\left(\sqrt{\alpha}+\sqrt{h}\right),
\qquad
\sqrt{\alpha}+\sqrt{h}\rightarrow0^{+},
\]
 
see \eqref{1.3}. Thus, the DNN is trained using convexification
reconstructions with the above a priori accuracy estimate as its input
images.

\subsection{The testing procedure}

\label{sec:7.5}

After the network has been trained, its parameters are fixed and no further
minimization of loss function \eqref{3} is performed. For each testing
sample, we again first apply the semi-discrete convexification method on the
coarse grid and obtain a new reconstruction
$\sigma_{\mathrm{conv}}(x)$. This reconstruction is interpolated onto the
fine grid and used directly as the input image to the already trained DNN.

Thus, the testing procedure has the form
\[
\sigma_{\mathrm{conv}}
\longrightarrow
\mathrm{trained\ DNN}
\longrightarrow
\sigma_{\mathrm{pred}}.
\]

The output $\sigma_{\mathrm{pred}}$ is the final refined reconstruction
produced by the trained network. Therefore, the convexification method is
used as the first stage for every sample in both the training and testing
procedures, while optimization of the neural-network parameters is performed
only during training.

\subsection{DnCNN denoising front-end}

\label{sec:7.5}

The first stage is a 17-layer DnCNN:

\begin{itemize}
\item \textbf{Input:} Single-channel ($1\times 256\times 256$) gray scale
image.

\item \textbf{Architecture:}

\begin{itemize}
\item Initial layer: $\text{Conv2D}(1\rightarrow 64, \text{kernel}=3) + 
\text{ReLU}$

\item 15 intermediate blocks: $\text{Conv2D}(64\rightarrow 64, \text{kernel}%
=3) + \text{BatchNorm} + \text{ReLU}$

\item Final layer: $\text{Conv2D}(64\rightarrow 1, \text{kernel}=3)$
\end{itemize}

\item \textbf{Output:} Residual learning output, $x - F(x)$.
\end{itemize}

\subsection{ResUNet with SE self-attention}

\label{sec:7.6}

The denoised image is passed to a Residual U-Net with SE attention blocks:

\begin{itemize}
\item \textbf{Encoder:}

\begin{itemize}
\item 3 stages with downsampling via $\text{MaxPool2D}(2)$

\item Residual blocks with output channels: $64 \rightarrow 128 \rightarrow
256$

\item Each residual block includes SE attention
\end{itemize}

\item \textbf{Bottleneck:} One residual block with 512 channels

\item \textbf{Decoder:}

\begin{itemize}
\item 3 stages with $\text{ConvTranspose2D}$ upsampling

\item Skip connections from encoder

\item Residual blocks with output channels: $256 \rightarrow 128 \rightarrow
64$
\end{itemize}

\item \textbf{Final Layer:} $\text{Conv2D}(64 \rightarrow 1) + \text{Sigmoid}
$
\end{itemize}

\subsection{Squeeze-and-excitation (SE) block}

\label{sec:7.7}

Each residual block is augmented with a channel-wise attention mechanism:

\begin{itemize}
\item Global average pooling across spatial dimensions

\item Two-layer fully-connected network: 
\[
\text{FC}(C \rightarrow C/r \rightarrow C) + \text{Sigmoid}, \quad r=16 
\]

\item The output attention weights modulate the original feature maps
channel-wise
\end{itemize}

\subsection{Loss function}

\label{sec:7.8}

We adopt a combined loss of Multi-Scale Structural Similarity (MS-SSIM) and
L1 loss (\ref{3}). This loss function encourages both: perceptual fidelity
and pixel-wise accuracy, especially for preserving the fine stroke
structures in characters we image.

\subsection{Dataset construction and partition}

\label{sec:dataset_construction}

The Chinese-character images used to construct the conductivity phantoms were
taken from the publicly available data set ``Handwritten Chinese Character
(Hanzi) Datasets,'' available at
\url{https://www.kaggle.com/datasets/pascalbliem/handwritten-chinese-character-hanzi-datasets/data}.
From this data set, we randomly selected $3{,}256$ samples.

The selected Chinese-character images determine the shapes of the conductivity
phantoms used in our numerical experiments. As described in Section~\ref{sec:8},
the conductivity is set to $\sigma=2$ inside the character shape and
$\sigma=1$ in the background.

For every selected sample, the deep-learning procedure uses a paired
input-target representation. The input image is the conductivity reconstruction
obtained by applying the semi-discrete convexification method on the coarse
grid, while the corresponding original Chinese-character conductivity image
is used as the ground-truth target. Thus, each sample has the form
\[
\left(
\sigma_{\mathrm{conv}}^{(i)},
\sigma^{(i)\ast}
\right),
\]
where $\sigma_{\mathrm{conv}}^{(i)}$ denotes the coarse-grid convexification
reconstruction and $\sigma^{(i)\ast}$ denotes the corresponding true
conductivity image. In the numerical implementation, reconstructed images and
their corresponding original images are paired according to their filenames.
All images are converted to grayscale, resized to a $256\times256$ pixel grid,
and converted to tensors before being supplied to the network.

The complete set of $3{,}256$ paired samples is randomly partitioned into
training, validation, and testing subsets using proportions of
$80\%$, $10\%$, and $10\%$, respectively. A fixed random seed of $42$ is used
to make this partition reproducible. With the integer partition used in our
implementation, this gives $2{,}604$ training samples, $325$ validation
samples, and $327$ testing samples. No individual image is shared between
these three subsets.

The partition is performed at the individual-image level rather than at the
Chinese-character-class level. Consequently, although the individual images
in the training, validation, and testing subsets are distinct, different
samples belonging to the same Chinese-character class may occur in different
subsets. Therefore, the held-out testing set measures generalization to unseen
images from the available data distribution rather than generalization to
completely unseen character classes.

No data augmentation is used in the reported experiments. In particular, no
random rotations, translations, flips, crops, geometric distortions, or other
synthetic transformations are applied during training. The only image
preprocessing consists of grayscale conversion, resizing to
$256\times256$, and conversion to tensor form.

\subsection{Training details}

\label{sec:7.9}

\begin{itemize}

\item \textbf{Dataset:} The construction, preprocessing, and partition of the
$3{,}256$ paired images are described in
Section~\ref{sec:dataset_construction}.

\item \textbf{Training/validation/testing split:}
$2{,}604$/$325$/$327$ samples, respectively, corresponding to the
$80\%$/$10\%$/$10\%$ partition described in
Section~\ref{sec:dataset_construction}.

\item \textbf{Image resolution:} $256\times256$ pixels.

\item \textbf{Optimizer:} AdamW
($\mathrm{lr}=10^{-3}$, weight decay $=10^{-5}$).

\item \textbf{Scheduler:} ReduceLROnPlateau
(factor $0.5$, patience $3$).

\item \textbf{Epochs:} $200$.

\item \textbf{Batch size:} $4$.

\item \textbf{Framework:} PyTorch with CUDA acceleration.

\end{itemize}

During training, the trainable parameters of the combined DnCNN--ResUNet
model are updated using AdamW with the initial learning rate and weight
decay specified above. The objective is the combined MS-SSIM and $L^1$
reconstruction loss in \eqref{3}. After each epoch, the loss and PSNR are
computed on the validation subset. The ReduceLROnPlateau scheduler monitors
the validation loss, uses a patience of three epochs, and reduces the
learning rate by a factor of $0.5$ when the scheduler's plateau criterion
is met. The validation samples are not used for gradient updates.

The model is trained for the prescribed $200$ epochs without early stopping.
The model parameters obtained at the end of the last epoch are saved and
used for testing, rather than selecting a checkpoint according to test-set
performance. The held-out test set is not used for gradient updates,
learning-rate adjustment, or checkpoint selection. During testing, the
saved model parameters remain fixed.

\subsection{Selection of loss function (\protect\ref{3})}

\label{sec:7.11}

We now explain the reason of the choice of the loss function in the form of (%
\ref{3}), in which a combined MS-SSIM and L1 loss is used. MS-SSIM has
assigned a higher weight to emphasize structural fidelity. This loss
formulation proved effective for reconstructing Chinese characters, as
MS-SSIM encourages multi-scale structural consistency, helping to preserve
stroke continuity and global character geometry, while the L1 term enforces
pixel-level accuracy and sharpness. The combination enables the model to
suppress mosaic artifacts and coherently reconstruct missing or degraded
regions, producing visually consistent characters that maintain both local
detail and overall structural integrity.

\subsection{Computational cost and implementation conditions}

\label{sec:computational_cost}

All coarse-grid and fine-grid convexification computations used in the timing
comparison were performed on the same computer and with the same code. The
computer was equipped with an Intel Core i7-9750H CPU (6 cores, 12 threads,
2.6 GHz) and 32 GB of RAM, running 64-bit Windows 10. Minimization in the
convexification procedure was carried out using the MATLAB built-in routine
\texttt{fminunc}. The grid step sizes in the timing comparison were
$h=0.05$ for the coarse grid of the present paper and $h=0.025$ for the
fine grid of \cite[p.~17]{EITIP2025}. The parameters $\alpha=0.01$,
$\epsilon=0.0002$, and $\kappa=3$ were identical in both computations.

On this machine, the fine-grid convexification computation required
approximately 57 minutes per image, whereas the coarse-grid computation
required approximately 4 minutes per image. The resulting speedup is
\[
\frac{57}{4}=14.25.
\]
This ratio refers exclusively to the convexification stage and excludes
both training and inference of the DNN, as well as computation of the
forward problem. In particular, it is not an end-to-end speedup of the
complete two-stage procedure.

The training of the DNN is a one-time offline cost. Training for 200 epochs
required approximately 20 hours on an NVIDIA GeForce RTX 3080 Ti GPU.
Once training is completed, the same trained network is used for the test
samples, and each refined reconstruction is produced in a single forward
pass without further optimization of the network parameters. The cost of
this inference step is separate from the convexification timing comparison;
no separate numerical measurement of inference time is reported here.

\subsection{Model summary table}

\label{sec:7.12}

\hbox{}\hfil\break

\begin{table}[!h]
\caption{Architecture Summary of Combined DnCNN + ResUNet Model.}
\label{tab:metrics}\centering
\begin{tabular}{|l|l|}
\hline
\textbf{Component} & \textbf{Details} \\ \hline
DnCNN & 17-layer residual denoiser (Conv-BN-ReLU) \\ 
ResUNet & 3-level encoder-decoder with skip connections \\ 
Residual Block & 2 conv layers + SE attention + identity skip \\ 
SE Block & Channel-wise self-attention via squeeze and excitation \\ 
Loss & $\gamma=0.84$: MS-SSIM + L1 \\ 
Output & Sigmoid to scale prediction to [0, 1] \\ \hline
\end{tabular}%
\end{table}

\subsection{Quantitative Evaluation and Robustness Analysis}

\label{sec:quantitative_robustness}

The preceding qualitative reconstructions demonstrate visually that the
proposed model can recover Chinese-character stroke structures from degraded
inputs. However, visual inspection alone does not establish whether the
reconstruction consistently improves across the test population. We
therefore performed a quantitative evaluation on the complete held-out test
set of 327 image pairs. The evaluation compares three conditions: the
degraded input without additional processing, a conventional U-Net baseline,
and the proposed ResUNet-based reconstruction model.

All methods were evaluated against the corresponding original character
images using peak signal-to-noise ratio (PSNR), multi-scale structural
similarity (MS-SSIM), mean absolute error (MAE), and Stroke Dice. PSNR
provides a conventional measure of pixel-level reconstruction fidelity,
whereas MS-SSIM measures structural similarity over multiple spatial scales.
Stroke Dice was additionally included because the primary objective of this
application is the recovery of the foreground stroke configuration that
defines the Chinese character rather than exact replication of every image
pixel. Higher PSNR, MS-SSIM, and Stroke Dice indicate better performance,
while lower MAE indicates better performance.

\begin{table}[htbp]
\caption{Quantitative reconstruction performance on the held-out test set.
Higher PSNR, MS-SSIM, and Stroke Dice are better, while lower MAE is better.}
\label{tab:test_metrics}\centering
{\small \ 
\resizebox{\linewidth}{!}{
    \begin{tabular}{lcccc}
        \hline
        Method &
        PSNR (dB) $\uparrow$ &
        MS-SSIM $\uparrow$ &
        MAE $\downarrow$ &
        Stroke Dice $\uparrow$ \\
        \hline
        Degraded Input
        & 11.0647
        & 0.6082
        & 0.2551
        & 0.5588 \\

        U-Net
        & 16.2677
        & 0.7878
        & 0.0693
        & 0.7234 \\

        Proposed ResUNet
        & \textbf{16.4068}
        & \textbf{0.7961}
        & \textbf{0.0672}
        & \textbf{0.7344} \\
        \hline
    \end{tabular}
    } }
\end{table}

As shown in Table~\ref{tab:test_metrics}, both neural reconstruction models
substantially improve the degraded inputs. Relative to the unprocessed
input, the U-Net increases the mean PSNR by approximately $5.20$ dB, while
the proposed model increases it by approximately $5.34$ dB. Similarly, the
U-Net increases MS-SSIM by approximately $0.1796$, whereas the proposed
model increases it by approximately $0.1879$. The proposed model also
achieves the highest mean Stroke Dice and the lowest MAE.

The improvement over U-Net on the original test distribution is modest in
magnitude, indicating that a conventional U-Net is already an effective
baseline when sufficient structural information remains in the degraded
input. Nevertheless, the improvement is consistent across a majority of the
held-out samples. The proposed model outperforms U-Net on approximately $%
61.2\%$ of test images according to PSNR, $67.9\%$ according to MS-SSIM, $%
71.3\%$ according to Stroke Dice, and $65.1\%$ according to MAE. The
stronger consistency observed for MS-SSIM and Stroke Dice is particularly
relevant to the objective of recovering the structural arrangement of
Chinese-character strokes.

\begin{figure}[htbp]
\centering
\includegraphics[width=\linewidth]{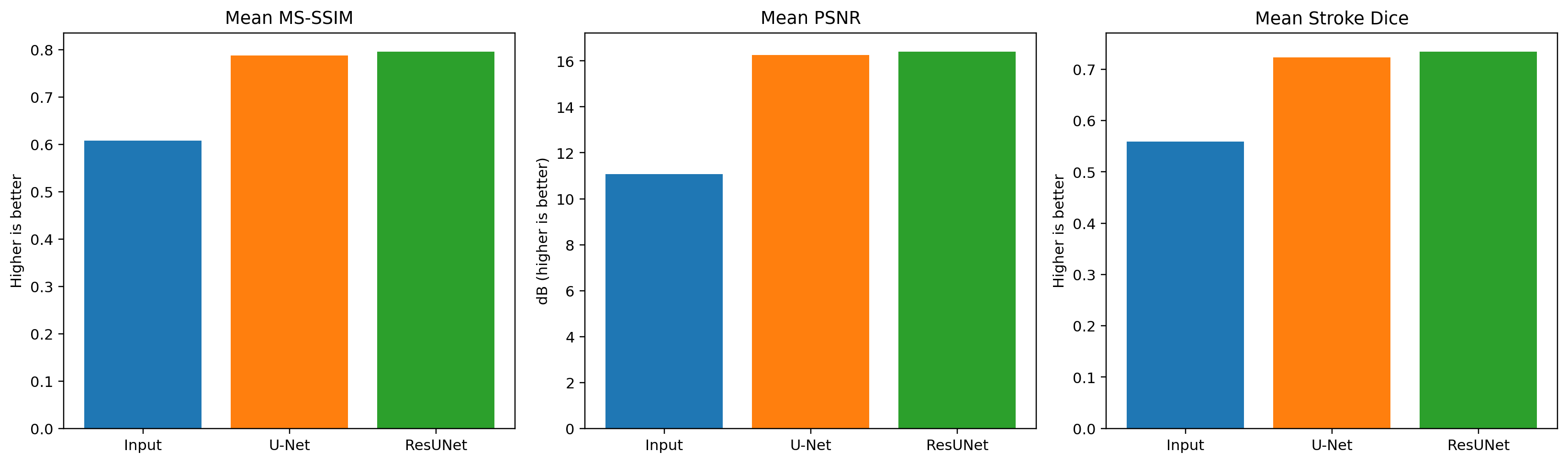}
\caption{Comparison of reconstruction metrics for the degraded input,
conventional U-Net, and the proposed ResUNet-based model on the held-out
test set.}
\label{fig:metric_comparison}
\end{figure}

Figure~\ref{fig:metric_comparison} provides a visual summary of the mean
performance of the degraded input, U-Net, and the proposed model.

\subsubsection{Robustness to Additional Input Corruption}

\label{sec:noise_robustness}

We note that the stability of the convexification method with respect 
to noise in the input data was numerically investigated in Tests 4 
and 5 of \cite{EITIP2025}. In those tests, noise was introduced into 
the input data of the convexification method, and the obtained 
numerical results demonstrated its stability with respect to 
such perturbations.

Motivated by these results, in the present work we also investigate 
the robustness of both learned reconstruction models with respect 
to additional noise. In addition to evaluating their performance 
on the original held-out test distribution, we examine their 
robustness under controlled additional corruption. More precisely, 
zero-mean Gaussian noise is added to the conductivity images 
reconstructed by the coarse-grid convexification method, which 
serve as the inputs to the trained networks, with noise level $\xi \in \{0.02,0.05,0.10\}$. The original test 
condition, without additional noise, corresponds to $\xi=0$. 
Neither model is retrained or fine-tuned for these additional 
noise levels. Instead, the same fixed checkpoints are used 
at each noise level.
\clearpage
\begin{table}[htbp]
\caption{Reconstruction performance under increasing Gaussian corruption.
Higher PSNR, MS-SSIM, and Stroke Dice are better, while lower MAE is better.}
\label{tab:noise_robustness}\centering
{\small \ 
\resizebox{\linewidth}{!}{
    \begin{tabular}{lccccc}
        \hline
        Model &
        Gaussian $\xi$ &
        PSNR (dB) $\uparrow$ &
        MS-SSIM $\uparrow$ &
        MAE $\downarrow$ &
        Stroke Dice $\uparrow$ \\
        \hline

        U-Net
        & 0.00
        & 16.2677
        & 0.7878
        & 0.0693
        & 0.7234 \\

        Proposed ResUNet
        & 0.00
        & \textbf{16.4068}
        & \textbf{0.7961}
        & \textbf{0.0672}
        & \textbf{0.7344} \\

        \hline

        U-Net
        & 0.02
        & \textbf{15.9585}
        & 0.7657
        & \textbf{0.0734}
        & 0.7184 \\

        Proposed ResUNet
        & 0.02
        & 15.7351
        & \textbf{0.7710}
        & 0.0739
        & \textbf{0.7277} \\

        \hline

        U-Net
        & 0.05
        & 12.8963
        & 0.6194
        & 0.1238
        & 0.6179 \\

        Proposed ResUNet
        & 0.05
        & \textbf{13.7949}
        & \textbf{0.6307}
        & \textbf{0.1049}
        & \textbf{0.6566} \\

        \hline

        U-Net
        & 0.10
        & 4.6051
        & 0.1485
        & 0.4899
        & 0.3410 \\

        Proposed ResUNet
        & 0.10
        & \textbf{13.1518}
        & \textbf{0.6020}
        & \textbf{0.1144}
        & \textbf{0.5805} \\

        \hline
    \end{tabular}
    } }
\end{table}

At the lowest additional corruption level, $\xi=0.02$, the two
architectures remain relatively close. U-Net achieves slightly higher PSNR
and slightly lower MAE, whereas the proposed model maintains higher MS-SSIM
and Stroke Dice. This indicates that, under mild additional corruption,
neither architecture uniformly dominates the other across all metrics.

As the corruption level increases, however, the difference becomes more
pronounced. At $\xi=0.05$, the proposed model outperforms U-Net on all
four metrics, achieving a PSNR of $13.7949$ dB compared with $12.8963$ dB,
MS-SSIM of $0.6307$ compared with $0.6194$, MAE of $0.1049$ compared with $%
0.1238$, and Stroke Dice of $0.6566$ compared with $0.6179$.

The largest difference occurs at $\xi=0.10$. Under this level of
corruption, the U-Net reconstruction deteriorates substantially, reaching
only $4.6051$ dB PSNR, $0.1485$ MS-SSIM, and $0.3410$ Stroke Dice. In
contrast, the proposed model retains $13.1518$ dB PSNR, $0.6020$ MS-SSIM,
and $0.5805$ Stroke Dice. Its MAE also remains substantially lower ($0.1144$
compared with $0.4899$ for U-Net).

Figure~\ref{fig:noise_robustness} illustrates the change in MS-SSIM as the
additional Gaussian corruption increases.

\clearpage

\begin{figure}[htbp]
\centering
\includegraphics[width=\linewidth]{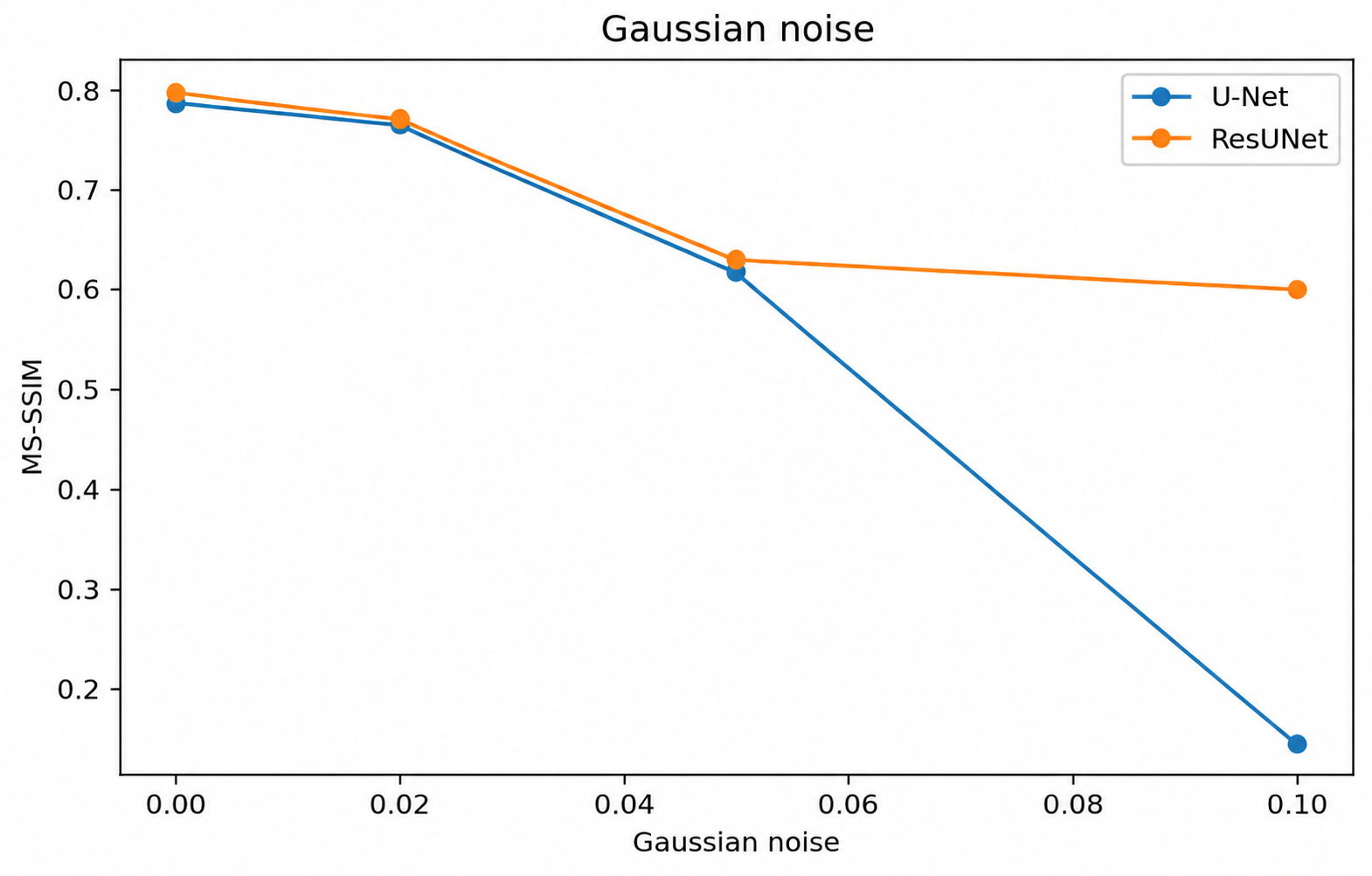}
\caption{MS-SSIM of U-Net and the proposed ResUNet-based model under
increasing Gaussian input corruption. The performance difference becomes
substantially larger at higher corruption levels.}
\label{fig:Gaussian Noise MS-SSIM Comparison}
\end{figure}

These results indicate that the additional components of the proposed
architecture provide their greatest advantage when the available input
information becomes substantially corrupted. Because the proposed
architecture contains multiple components beyond the conventional U-Net,
including the DnCNN preprocessing stage, residual connections, and
squeeze-and-excitation attention, the present comparison should be
interpreted as evidence for the robustness of the complete proposed
architecture rather than as an isolated ablation of any single component.

\subsubsection{Qualitative Comparison}

\label{sec:qualitative_comparison}

Representative examples from the held-out test set are shown in Figure~\ref%
{fig:matched_examples}. These examples supplement the quantitative analysis by illustrating how
differences in MS-SSIM and Stroke Dice correspond to visible differences
in the recovered stroke structures. MS-SSIM measures the structural
similarity between the reconstructed and true images at multiple spatial
scales, while Stroke Dice measures the spatial overlap between their
foreground character strokes, with $1$ indicating perfect overlap and
$0$ indicating no overlap. Higher values of both metrics indicate better
reconstruction quality. The examples are used only as qualitative
illustrations; the conclusions above are based on the complete held-out
test set rather than on manually selected reconstruction examples.

\clearpage

\begin{figure}[htbp]
\centering
\includegraphics[width=0.5\textwidth]{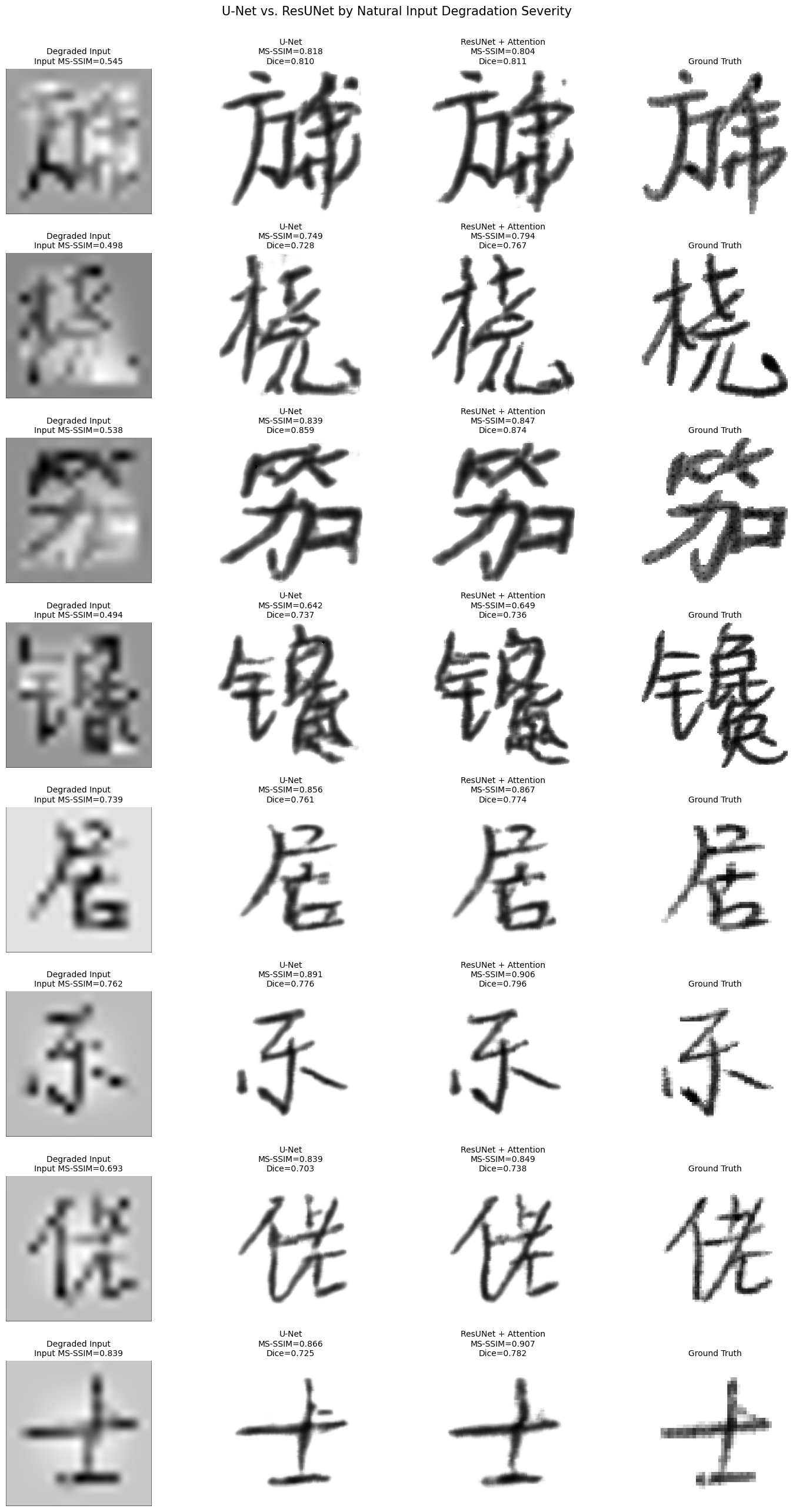}
\caption{Comparison of reconstruction results on the held-out test set.
Left column: images obtained by the semi-discrete convexification method.
Second column: reconstructions obtained using the baseline U-Net.
Third column: reconstructions obtained using the proposed ResUNet-based
model. Right column: the corresponding true images. The MS-SSIM and
Stroke Dice values measure structural similarity and stroke overlap with
the true images, respectively.}
\label{fig:matched_examples}
\end{figure}

\subsubsection{Scope of Out-of-Distribution Evaluation}

\label{sec:ood_scope}

The held-out test images used above are completely excluded from model
training and therefore provide an evaluation on unseen samples. However,
because these images originate from the same available dataset, we do not
characterize this experiment as a full out-of-distribution (OOD) evaluation.

The present dataset contains 3,256 paired images and consequently represents
only a limited subset of the substantially larger space of Chinese
characters, writing styles, fonts, and possible degradation conditions.
Constructing an artificial OOD split from such a limited corpus, for example
by withholding a small number of characters or styles, could make the
resulting performance highly dependent on the particular subset selected and
would not necessarily provide a representative measure of general OOD
behavior.

A comprehensive OOD evaluation would therefore be more appropriately
performed using a substantially larger and independently collected
Chinese-character dataset in which character identity, writing or font
style, and acquisition or degradation conditions can be explicitly separated
between training and evaluation data. Such a large-scale evaluation is
beyond the scope of the present dataset and is left for future work. The
controlled Gaussian-corruption experiment in Section~\ref%
{sec:noise_robustness} instead provides an initial assessment of the
robustness of the trained models when the input distribution is
systematically perturbed beyond the original test condition.

\section{\ Numerical Results}

\label{sec:8}

Our test phantoms are based on Chinese characters, which provide a
challenging test case due to their intricate shapes, sharp corners, and
combination of thick and thin strokes. The conductivity is set as $\sigma =2$
inside the character shape and $\sigma =1$ in the background. The images
shown here are representative samples drawn from the testing set, evaluated
after the training phase was completed. The testing procedure was arranged as follows. For each testing case, we first applied the above-described semi-discrete convexification
method on a coarse grid. Next, the resulting function 
$\sigma_{conv}\left(\mathbf{x}\right)$ was used as the input to the 
trained DNN. The final reconstructed image was then obtained by applying
the above-described DNN procedure.



We compare the reconstruction results obtained from two different spatial
discretizations of the coarse grids:

\begin{enumerate}
\item A too coarse grid with the step size%
\begin{equation}
h_{1}=0.1.  \label{7.4}
\end{equation}

\item A coarse grid with a finer grid with the step size 
\begin{equation}
h_{2}=0.05.  \label{7.5}
\end{equation}%
\newline
\end{enumerate}

\begin{figure}[]
\centering
\includegraphics[width=0.50\textwidth]{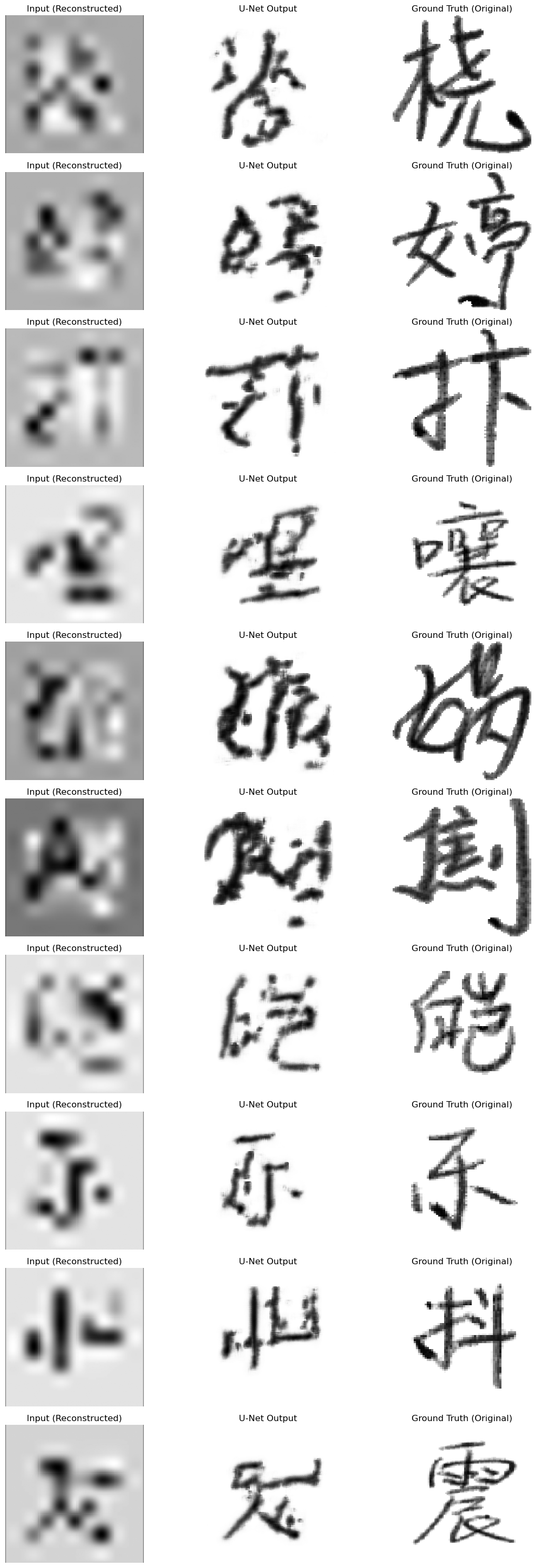}
\caption{Left column: image obtained by the semi-discrete version of the
convexification method on a too coarse grid with the grid step size $%
h_{1}=0.1$ as in (\protect\ref{7.4}). Right column: the true image. Middle
column: the reconstructed image after deep learning. These results are
visibly worse than those obtained with the a finer coarse grid step size $%
h_{2}=0.05$ as in (\protect\ref{7.5}) in the next Figure \protect\ref%
{fig:qualitative_results}. Apparently, the coarse grid step size $h_{1}=0.1$
is too large.}
\label{fig:poor_mesh_results}
\end{figure}

\begin{figure}[]
\centering
\includegraphics[width=0.50\textwidth]{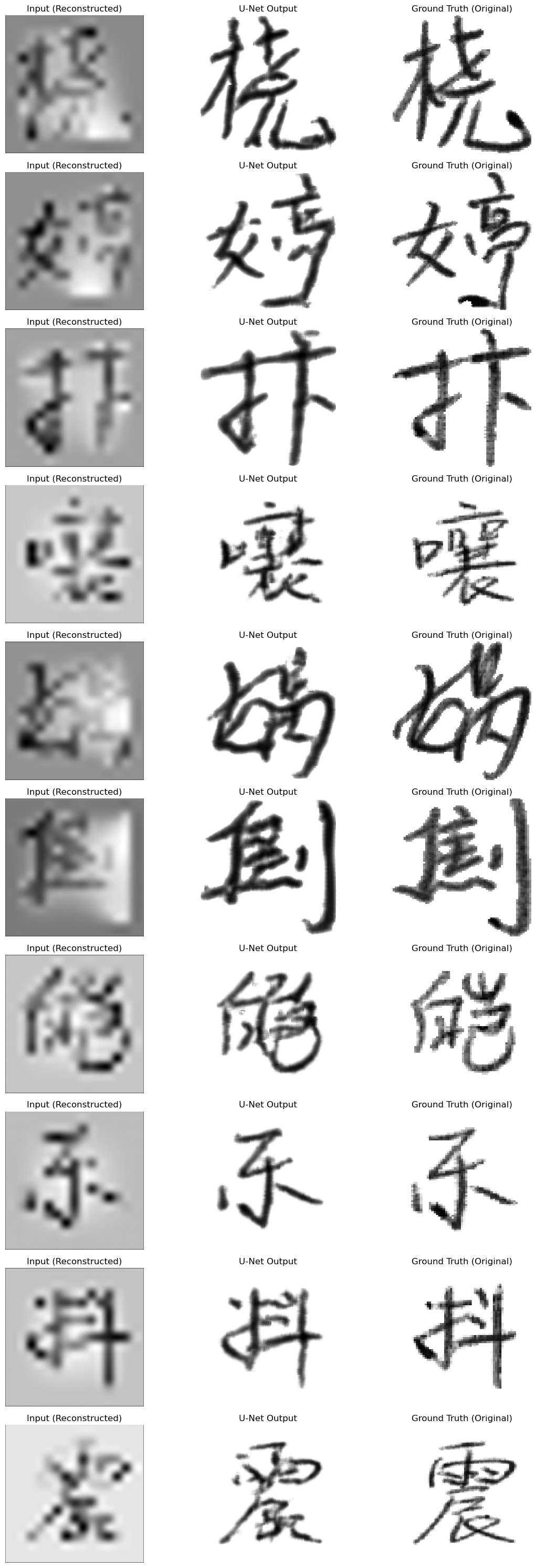}
\caption{ Left column: the image obtained by the semi-discrete version of
the convexification method on a coarse mesh with a finer coarse grid with
step size $h_{2}=0.05$ as in (\protect\ref{7.5}). The \emph{a priori}
accuracy estimate of the starting point is given in (\protect\ref{7.2}).
Right column: the true image. The same characters are used in this column as
in \protect\ref{fig:poor_mesh_results}. Middle column: the reconstructed
image after deep learning. One can observe an excellent reconstruction
accuracy in the middle column. Our algorithm recovers quite well intricate
structures and sharp contours.}
\label{fig:qualitative_results}
\end{figure}

\begin{figure}[]
\centering
\includegraphics[width=0.50\textwidth]{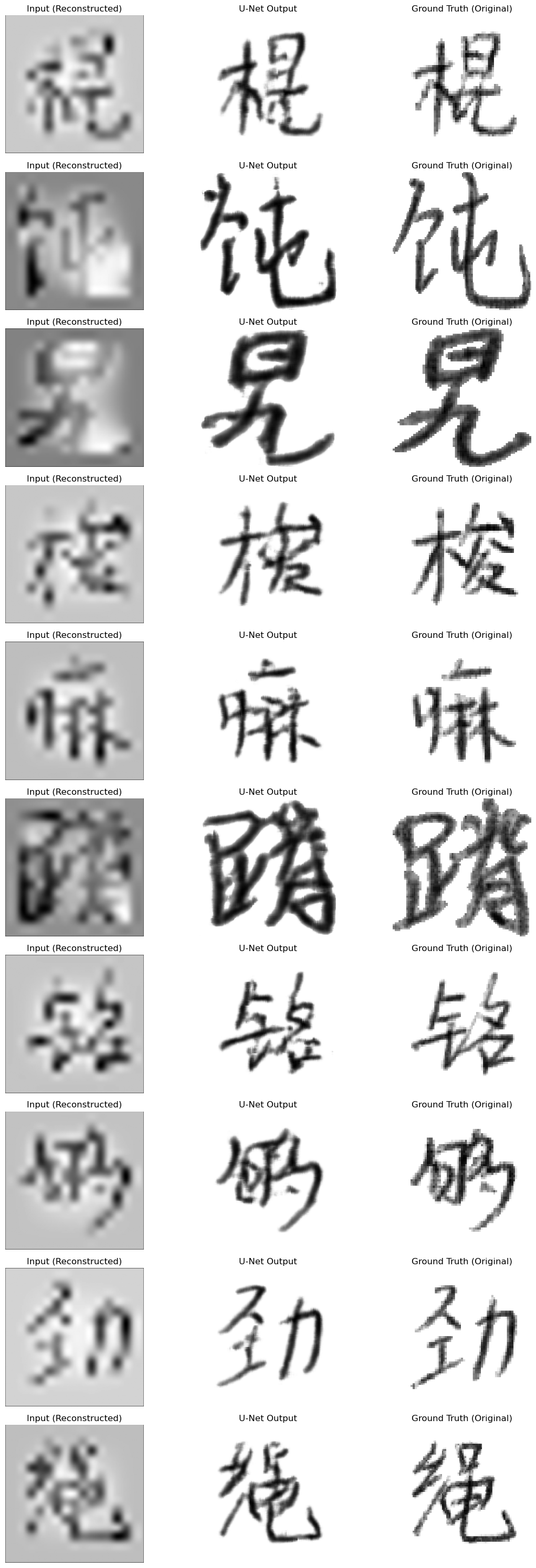}
\caption{This is the second series of images with the grid step size of a
finer coarse grid $h_{2}=0.05$ as in (\protect\ref{7.5}). Left column: image
obtained by the semi-discrete version of the convexification method with $%
h_{2}=0.05.$ Right column: the true image. Middle column: the reconstructed
image after deep learning. The same high reconstruction accuracy as the one
of Figure \protect\ref{fig:qualitative_results} is observed.}
\label{fig:anotherqualitative_results}
\end{figure}

Images on Figures \ref{fig:poor_mesh_results}, \ref{fig:qualitative_results}
and \ref{fig:anotherqualitative_results} are obtained on the testing step of
our deep learning procedure. Left columns display results of the performance
of the semi-discrete version of the convexification method with two
different grid step sizes (\ref{7.4}) and (\ref{7.5}). Images in the left
columns were obtained by the semi-discrete convexification method. They are
blurry since a coarse mesh with either of step sizes (\ref{7.4}) or (\ref%
{7.5}) was used.

\clearpage

\clearpage

\clearpage

\clearpage

The same characters in the right column are used in Figures \ref%
{fig:poor_mesh_results} and \ref{fig:qualitative_results}. This is done for
the sake of comparisons of the results of our deep learning procedure for
two different grid step sizes (\ref{7.4}) and (\ref{7.5}) of the
convexification method. A clear qualitative difference can
be observed between the results obtained from the two coarse-grid
discretizations of the semi-discrete convexification method, whose
reconstructions are used as input images to the DNN. In the case (\ref{7.5}%
) of a finer coarse grid, the output of our hybrid U-Net model closely
resembles the ground truth Chinese characters. Most of the fundamental
stroke structures such as horizontal strokes, vertical strokes,
right-falling strokes, and left-falling strokes are well preserved. In
particular, the network successfully identifies straight-line features as
straight lines and reconstructs the relative orientation and placement of
slanted strokes with reasonable accuracy. Moreover, the finer coarse grid
result demonstrates a strong ability to distinguish between dots and short
line segments, which is essential for the legibility of Chinese characters.

One can observe that the too coarse grid with $h_{1}=0.1$ as in (\ref{7.4})
produces outputs of Figure \ref{fig:poor_mesh_results} that are largely
illegible and fail to effectively resemble the original characters. Our
hybrid U-Net model struggles to differentiate between vertical strokes and
right-falling strokes, as well as between dots and line segments. This loss
of structural distinction significantly degrades readability and makes
character recognition nearly impossible in most cases.

In contrast to Figure \ref{fig:poor_mesh_results}, the primary limitation of
the computational result of Figures \ref{fig:qualitative_results} and \ref%
{fig:anotherqualitative_results} with $h_{2}=0.05$ as in (\ref{7.5}) appears
in characters with highly crowded structures, especially those involving
dense combinations of horizontal and vertical strokes forming square-like
patterns. In these cases, the proximity of strokes reduces the available
spatial resolution, making it more challenging for the network to accurately
separate individual components. Nevertheless, aside from these tightly
packed characters, the finer mesh consistently yields accurate and visually
coherent reconstructions. These observations highlight the critical
role of the grid step size. 

\section{Summary}

\label{sec:9}

We have proposed a new two-stage numerical procedure for Coefficient Inverse
Problems. On the first stage, a semi-discrete version of the convexification
numerical method is applied on a coarse grid. The reason why the coarse grid
is used is that the convexification performs slowly. On the second stage
deep learning refines the images obtained on the first stage. A
convergence analysis of that semi-discrete version was carried out.
Numerical experiments demonstrate both a high speed of computations and a
good accuracy of reconstructions of complicated media structures.

\textbf{Acknowledgments: }The work of MVK was partially supported by the
National Science Foundation grant DMS 2436227. The authors are grateful to
Dr. Jingzhi Li and Dr. Alexandre Timonov for useful discussions. Dr. Zhipeng
Yang has computed the above described numerical results for the
convexification method on sparce grids. He has obtained those results using
the software, which he has developed for the publication \cite{EITIP2025}.
The authors are grateful to Dr. Yang for his kindness by allowing the
authors to use his results in this publication.


\begin{thebibliography}{99}
\bibitem{Beilina1} L. Beilina, Domain decomposition finite element/finite
difference method for the conductivity reconstruction in a hyperbolic
equation, Commun. Nonlinear Sci. Numer. Simul. 37, 222--237, 2016.

\bibitem{Beilina2} L. Beilina and E. Lindstr\"{o}m, An adaptive finite
element/finite difference domain decomposition method for applications in
microwave imaging, Electronics, 11, 1359, 2022.

\bibitem{Borcea} L.~Borcea. \newblock Electrical impedance tomography. %
\newblock {\em Inverse problems}, 18(6):R99, 2002.

\bibitem{BukhKlib} A.~L. Bukhgeim and M.~V. Klibanov, Uniqueness in the
large of a class of multidimensional inverse problems, Soviet Math. Doklady,
17 (1981), pp.~244--247.

\bibitem{Chavent} G.~Chavent. 
\newblock {\em Nonlinear Least Squares for Inverse Problems: Theoretical
	Foundations and Step-by-Step Guide for Applications}. \newblock Springer
Science \& Business Media, Berlin, 2010.

\bibitem{Giorgi} G.~Giorgi, M.~Brignone, R.~Aramini, and M.~Piana. 
\newblock {Application of the inhomogeneous Lippmann--Schwinger equation to
	inverse scattering problems}. \newblock {\em SIAM J. Appl. Math.},
73:212--231, 2013.

\bibitem{Gonch1} A.~V. Goncharsky and S.~Y. Romanov. \newblock Iterative
methods for solving coefficient inverse problems of wave tomography in
models with attenuation. \newblock {\em Inverse Probl.}, 33:025003, 2017.

\bibitem{Gonch2} A.~V. Goncharsky, S.~Y. Romanov, and S.~Y. Seryozhnikov. %
\newblock On mathematical problems of two-coefficient inverse problems of
ultrasonic tomography. \newblock {\em Inverse Probl.}, 40:045026, 2024.

\bibitem{GT} G.~Gilbarg and N.~S. Trudinger. 
\newblock {\em Elliptic Partial
Differential Equations of Second Order}. \newblock Springer-Verlag, New
York, second edition, 1983.

\bibitem{Goh} S. M. Goh, J. C. Ye, et al. \newblock Deep learning for EIT: a
physics-informed approach. \newblock {\em Inverse Problems}, 35(10):105008,
2019.

\bibitem{Harrach1} B. Harrach, The Calderon problem with finitely many
unknowns is equivalent to convex semi definite optimization, SIAM J. Math.
Anal., 55, 5666--5684, 2023.

\bibitem{Harrach2} B. Harrah and H. Meftahi, A monotonicity-based
globalization of the level-set method for inclusion detection, Commun. Anal.
Comput. 5, 1--17, 2025.

\bibitem{Harrach3} B. Harrah and A. Brojatsch, On the required number of
electrodes for uniqueness and convex reformulation in an inverse coefficient
problem, Inverse Problems, 41: 105011, 2025.

\bibitem{Isakov} V.~Isakov, 
\newblock {\em Inverse Problems for Partial
Differential Equations}. \newblock Springer, New York, 2006.

\bibitem{Klib95} M.V. Klibanov and O.V. Ioussoupova, Uniform strict
convexity of a cost functional for three dimensional inverse scattering
problem, SIAM\ J.\ Mathematical Analysis, 26, 147-179, 1995.

\bibitem{Klib97} M.~V. Klibanov. \newblock Global convexity in a
three-dimensional inverse acoustic problem. 
\newblock {\em SIAM J. Math.
Anal.}, 28:1371--1388, 1997.

\bibitem{KT} M.~V. Klibanov and A. Timonov, \emph{Carleman Estimates for
Coefficient Inverse Problems and Numerical Applications}, VSP, Utrecht, 2004.

\bibitem{KEIT} M.~V. Klibanov, J.~Li, and W.~Zhang. 
\newblock {Electrical impedance tomography with restricted Dirichlet-to-Neumann
	map data}. \newblock {\em Inverse Probl.}, 35:035005, 2019.

\bibitem{KL} M.~V. Klibanov and J.~Li, 
\newblock {\em {Inverse Problems and Carleman Estimates: Global Uniqueness,
		Global Convergence and Experimental Data}}. \newblock De Gruyter, Berlin,
2021.

\bibitem{KlibHJ} M.~V. Klibanov, L.~H. Nguyen, and H.~V. Tran. 
\newblock {Numerical viscosity solutions to Hamilton-Jacobi equations via a
	Carleman estimate and the convexification method}. 
\newblock {\em Journal
of Computational Physics}, 451:110828, 2022.

\bibitem{KLZhyp2} M.V. Klibanov, J. Li and W. Zhang, A globally convergent
numerical method for a 3D coefficiet inverse problem for a wave-like
equation, \emph{SIAM J. Scientific Computing}, 44, A3341--A3365, 2022.

\bibitem{Ktransp} M.~V. Klibanov, J.~Li, and Z.~Yang, \newblock %
Convexification for the viscosity solution for a coefficient inverse problem
for the radiative transport equation, 39, 125002, 2023.

\bibitem{EITIP2025} M.~V. Klibanov, J.~Li, and Z.~Yang. \newblock %
Convexification with the viscosity term for electrical impedance tomography. %
\newblock {\em Inverse Problems}, 41, 065020, 2025.

\bibitem{TTTP} M.V. Klibanov, J. Li, V.G. Romanov and Z. Yang,
Convexification for the 3D problem of travel time tomography, \emph{SIAM J.
Scientific Computing}, 47, A1436-A1457, 2025.

\bibitem{KT3} M.~V. Klibanov and A. Timonov, Acoustic imaging via a
viscosity approximation of an elliptic system generated by the Lavrent'ev
integral operator, \emph{SIAM J. Imaging Sciences}, 18, 1002--1027, 2025.

\bibitem{LRS} M.M. Lavrent'ev, V.G. Romanov and S.P. Shishatskii, \emph{%
Ill-Posed Problems of Mathematical Physics and Analysis}, AMS, Providence:
RI, 1986.

\bibitem{LN} H.P.N. Le, T.~T. Le, and L.~H. Nguyen. 
\newblock {The Carleman
convexification method for Hamilton-Jacobi equations}. 
\newblock {\em
Computers} \& \emph{Mathematics with Applications}, 159, 173--185, 2024.

\bibitem{Lokk} C.~Lokker, E.~Bagheri, W.~Abdelkader, R.~Parrish, M.~Afzal,
T.~Navarro, C.~Cotoi, F.~Germini, L.~Linkins, R.~B.~Haynes, L.~Chu and
A.~Iorio, Deep learning to refine the identification of high-quality
clinical research articles from the biomedical literature: Performance
evaluation, \emph{J. Biomed. Inform.} 142, 104384, 2023.

\bibitem{Nov} R.G. Novikov, The $\overline{\partial }$- approach to
approximate inverse scattering at fixed energy in three dimensions,
International Math. Research Peports, 6 (2005), 287-349.

\bibitem{Rom} V.G. Romanov, \emph{Investigation Methods for Inverse Problems}%
, VSP, Utrecht, The Netherlands, 2002.

\bibitem{Unet15} O. Ronneberger, P. Fischer and T. Brox. \newblock U-net:
Convolutional networks for biomedical image segmentation. \newblock In \emph{%
International Conference on Medical image computing and computer-assisted
intervention}, pages 234--241, Springer, 2015.

\bibitem{Santos} T.B.R. Santos, R.M. Nakanishi, T.M.M. Oleg\'{a}rio, R.G.
Lima and J.L. Mueller, Resolution improvement and algorithmic dependence of
machine learning for post-processing respiratory EIT images, \emph{Applied
Mathematics for Modern Challenges}, 1, 21-38, 2023.

\bibitem{Scales} J.A. Scales, M.L. Smith and T.L. Fisher, Global
optimization methods for multimodal inverse problems, \emph{J. Comput. Phys.,%
} 103, 258-268, 1992.

\bibitem{T} A.N. Tikhonov, A.V. Goncharsky, V.V. Stepanov and A.G. Yagola, 
\emph{Numerical Methods for the Solution of Ill-Posed Problems}, Kluwer,
London, 1995.
\end{thebibliography}
\end{document}